
\input amstex
\documentstyle{amsppt}
\NoRunningHeads
\magnification=1200
\define\Naim{Na\"\i m }
\def\Cap{\text{\rom{Cap} }}
\def\ca{\text{\rom{cap }}}
\def\co{\text{co }}
\def\supp{\text{supp }}

\hyphenation{Gre-noble}

\topmatter
\title  Nonlinear equations
and weighted norm  inequalities\endtitle
\author N. J. Kalton and I. E. Verbitsky
\endauthor
\address Nigel J. Kalton, Department of Mathematics, University of
Missouri, Columbia, Missouri 65211\endaddress
\email nigel\@math.missouri.edu\endemail
\address Igor E. Verbitsky, Department of Mathematics,
 University of
Missouri, Columbia, Missouri 65211\endaddress
\email igor\@math.missouri.edu\endemail
\thanks The first
author was partially supported by NSF grant DMS-9500125,
 and the second by NSF grant
DMS-9401493  and the University of Missouri
Research Board grant  RB-96029.
\endthanks
\abstract We study connections between
the problem of the existence of positive solutions
for certain nonlinear  equations and weighted norm inequalities.
In particular, we obtain explicit criteria for the
solvability of the  Dirichlet problem
 $$\aligned -& \Delta u = v  \, u^q + w, \quad  u \ge 0 \quad
  \text {on} \quad
 \Omega, \\ &u = 0 \quad \text {on} \quad \partial \Omega,
\endaligned $$
  on a regular domain $\Omega$ in $\bold R^n$ in the ``superlinear
case'' $q > 1$.  The coefficients $v, w$
are arbitrary positive measurable functions (or measures) on $\Omega$.
We also consider more
 general nonlinear differential and integral equations,
 and study the  spaces of
coefficients and solutions
 naturally associated with these problems,
as well as the corresponding capacities.

Our characterizations of the existence of positive solutions
 take into account the interplay
between $v$, $w$, and the corresponding Green's kernel.
They are not only
sufficient, but also necessary, and are established
without any a priori regularity assumptions on $v$ and $w$; we also
obtain sharp two-sided estimates of solutions up to the boundary.
Some of our results are new even if $v \equiv 1$ and
 $\Omega$ is a ball or half-space.

 The corresponding weighted norm inequalities are
proved for integral operators with kernels satisfying a refined
version of the so-called $3 G$-inequality by an elementary
  ``integration by parts'' argument. This also gives a new
unified proof for  some classical
inequalities including
the Carleson measure theorem for Poisson integrals and  trace
inequalities for
Riesz potentials and Green potentials.
\endabstract
\toc
\head 1. Introduction\endhead
\head 2. Superlinear problems and related function spaces\endhead
\head 3. Quasi-metric kernels and infinitesimal inequalities\endhead
\head 4. Weighted norm inequalities and nonlinear integral
equations\endhead
\head 5. Capacitary inequalities and criteria of solvability \endhead
\head 6. Trace inequalities, Carleson measure theorems and nonlinear
convolution equations\endhead
\head 7. Existence of positive solutions for superlinear Dirichlet
problems\endhead
\endtoc
\endtopmatter
\document

 \head {1. Introduction} \endhead

The main goal of this paper is to obtain explicit
criteria for the existence of positive solutions
 for  a class of ``superlinear'' Dirichlet problems on a
regular
domain $\Omega
\subset \bold R^n.$  In particular  we are interested in the
solvability of the Dirichlet problem
$$-\Delta u=v(x) \, u^q+w(x)\tag 1.1$$
on $\Omega$ where $n\ge 3$ subject
to the conditions $u\ge 0$ and $u=0$ on $\partial \Omega$, when $q>1$
and $v,w$ are given positive measurable functions.  By a solution of
(1.1)
we understand (more precise
definitions  are given below) a nonnegative measurable function
$u$  satisfying
a.e. on $\Omega$ the equivalent integral equation
$u = G(v \, u^q) + G w$
which follows by applying
the corresponding Green's potential $G = (-\Delta)^{-1}$
to both sides of (1.1).
It is also of
interest to replace $v$ and $w$ by positive measures, so that the
equation becomes $$-\Delta u= \sigma \, u^q +\omega\tag 1.2 $$
where $\sigma$ and $\omega$ are arbitrary positive Borel measures
on $\Omega$, and consider a more general boundary condition
$u = \phi$ on $\partial \Omega$ for a nonnegative
measurable function $\phi$; then
solutions of (1.2) are understood in the
analogous sense
and are defined
$\sigma$-a.e. on $\Omega$.

We observe that in the case $q=1$ these equations turn into
the inhomogeneous Schr\"odinger equation with  potentials
$v$ and $\sigma$ respectively. Thus (1.1) and (1.2) may be
referred to as the $q$-Schr\"odinger equations, and the techniques
needed to treat them are reminiscent of those
employed in \cite{12}, \cite{13}, \cite{16}, \cite{17}, \cite{23},
\cite{35}, \cite{47}, etc.

Equations of this type
are widely used in differential geometry, physics, astronomy,
 and numerous
applied problems with ``nonlinear sources'' (heat
transfer,  fluid flow, control theory, etc.; see \cite{32} and the
bibliography therein). They fall into the class of equations
with convex nonlinearities which  generally are
known to be more difficult
to investigate than equations with concave operators (see \cite{29}).
 In a sense, (1.1) and (1.2) have become model problems for nonlinear
analysts,
especially
after  H. Brezis and L. Nirenberg  studied  in \cite{8}
the homogeneous  problem  related to (1.1) with
 $w \equiv 0$ and $v \equiv 1$
  in the difficult ``critical case''
$q= (n+2)/(n-2)$.

 The solvability of the inhomogeneous problems (1.1) and (1.2)
with  variable  coefficients  has
 been studied  extensively,  mostly
under strong additional
assumptions on the potentials and data. We  mention
sufficient conditions of
 solvability for a general class of boundary value problems
with unbounded
coefficients and
domains $\Omega$  established by
M. Schechter (see \cite{47}). However, the necessary
conditions for solvability (criteria of ``nonsolvabilty'') which match
 sufficiency results
are more difficult to obtain. This problem is solved in the present
paper for a wide
class of nonlinear differential and integral equations.

A starting point for us in studying (1.2)  was the following
sharp  criterion
for the existence of  solutions
 in the case $\sigma \equiv 1$  and $\omega$
compactly supported in $\Omega$
due to  D. R. Adams and M. Pierre \cite{3}.
(Note that as explained below the ``critical index'' for this problem
is $q= n/(n-2)$.)

 \proclaim{Theorem 1.1 \cite{3}} Suppose that $1 < q < \infty$,
$\sigma \equiv 1$,
and  $\omega$ is a compactly supported measure
on a bounded regular domain $\Omega$.
\newline (1)  If there exists a solution for
$(1.2)$ which vanishes at the boundary, then
$$|E|_\omega \le C \, \text{\rm Cap}_p (E), \quad E \Subset \Omega,
\tag 1.3 $$
where $C$ is independent of compact sets $E$ and
$\text{\rm Cap}_p (\cdot)$ is the capacity
 associated with the Sobolev space $W^{2,p} (\bold R^n)$,
$1/p + 1/q =1$.\newline
(2) Conversely, (1.2) has a solution with zero boundary values
if (1.3) holds with a small enough
 constant $C < C(q, n, \omega, \Omega)$.
\endproclaim

We observe that Theorem 1.1 provides not only sufficient  but also
necessary conditions for solvability
 in explicit geometric
terms.
The classes of measures characterized by estimates of type
(1.3) are well studied
in potential theory
 starting from the work
of V. G. Maz'ya \cite{34} on  spectral
problems for Schr\"odinger operators in the early sixties.
(We  cite the books \cite{2},  \cite{24}, \cite{35}, and
\cite{57}  as an invaluable   source
on  capacities, potential theory, related function spaces and
applications to partial differential equations.)
There  also are equivalent
alternative
 characterizations of these classes which do not involve
capacities and play an important role in the sequel (see \cite{27},
\cite{37}, \cite{44}, \cite{45}, \cite{46}, \cite{52}).

Note  that
if $1<q < n/(n-2)$, then
$\inf_E \, \text{Cap}_p (E) > 0$ for all
$E \subset \text{supp} \, \omega$ so that $q=n/(n-2)$
 is a critical index. In the more interesting case $q \ge n/(n-2)$
however,  for
compactly supported $\omega$, one deals essentially with
 $\Omega = \bold R^n$ where the capacities
 $\text{Cap}_p (\cdot)$ are applicable.
They are not adequate for the analysis of
solutions and coefficients up to the boundary.

The study of the boundary behavior for nonlinear problems
of this type requires new methods.
Even   in the case $\sigma \equiv 1$ on a bounded domain with
smooth boundary
the solvability problem
 for $(1.2)$  was open
for  noncompactly supported $\omega$. A solution to this and more
general  problems in both
capacitary and non-capacitary form  is given below.

We remark that similar questions
for the Dirichlet problem
$$\left \{
\aligned  & \Delta u= \sigma \,  u^q + \omega , \quad  u \ge 0 \quad
  \text {on} \quad
 \Omega,\\
 & u = \phi \quad \text {on} \quad \partial \Omega,
\endaligned \right.  \tag 1.4$$
with $\Delta$ in place of $- \Delta$
 are in the center of
the current work of probabilists. They consider (1.4)
in the  case $\sigma \equiv 1$  to study  the
 so-called
 superdiffusions (see \cite{30}, \cite{14}, \cite{15})
for $1< q \le 2$.
 Unfortunately, at the moment
no probabilistic models seem to be known for (1.2),
 or (1.4) with $q > 2$. Criteria of solvability
for (1.2) and (1.4),
 at least
 in the well studied case
$\sigma \equiv 1$ and $\text{supp} \, \omega \Subset \Omega$,
are known to be
different: solutions for (1.4) exist under much weaker assumptions
(absolute continuity of $\omega$ with respect to the corresponding
capacity) than  for (1.2).

 In order to
study the problems stated above
we first develop a general technique for studying
a class of nonlinear operator equations. We then
consider
certain nonlinear integral
 equations related to (1.1), (1.2), or more general differential
equations via Green's functions and the corresponding potential theory.
  This
approach which is applicable to many other similar
problems is developed in Sections 2-5.

 Before discussing our approach in detail
we would like to make some general
comments. We do not use any variational theory,
weighted Sobolev spaces,
 Calderon-Zygmund decompositions,  or maximal function
inequalities.
In this respect our
approach  resembles some ideas of the  {\it original\/} proof
 of T. Wolff's inequality which appeared in the context
of the spectral synthesis problem for Sobolev spaces in
\cite{23}. (Note that the subsequent
alternative proofs of Wolff's inequality
 due to Per Nilsson, J. L. Lewis, and D. R. Adams,
see  \cite{2}, p. 126, are not enough for our purposes.) We develop
a new discrete decomposition for integral operators with respect to
an arbitrary measure. This
leads to sharp estimates of the nonlinear iterations of
Green's potentials
with precise estimates of the constants
involved which is most important in this paper.

We also would like to point out
interesting connections (in the easier part which involves
 $L^p$-estimates for integral operators) with
  the Hilbert space proof of the Carleson measure theorem
 due to S. A. Vinogradov (see \cite{40}) and its weighted analogue
used recently by S. Treil  and A. Volberg \cite{50} in the vector valued
version of the Hunt-Muckenhoupt-Wheeden theorem
in  case $p=2$.

 However,
we avoid using dyadic decompositions,  sophisticated
capacitary inequalities, or
any known  tests for boundedness of integral operators
or quadratic forms.    Our
proofs  here are based on a  quite elementary
argument which  resembles the proof of Hardy's inequality and works
for all $1<p< \infty$.
In particular it yields a simple proof
of some deep classical results (Carleson measure inequality \cite{10},
trace inequalities of Maz'ya-D. Adams-Dahlberg  \cite{2}, \cite{35})
and their generalizations. Another application of these ideas to
the problem of $\Lambda_p$ sets for
Legendre and Jacobi polynomials is given in  \cite{26}
 where the underlying space is
assumed to be discrete.
Similar inequalities   also appear in the problem of solvability
for  multidimensional
 Riccati's equations \cite{22}, spectral
estimates for
Schr\"odinger operators and
multipliers of Sobolev spaces \cite{37}, etc.

We now turn to a more detailed discussion of our main results.
Let $X$ be a metric space and suppose $\sigma$ is a fixed positive
$\sigma-$finite Borel measure on $X.$  Suppose $K:X\times X\to
[0,\infty]$ is a Borel kernel function. We write
$$Kf(x)= K (f \, d \sigma) (x) = \int_X K(x,y) \, f(y) \, d\sigma(y)$$
for any measurable function $f \ge 0$. Occasionally we write $K^\sigma f$
rather than
$Kf$ to emphasize
the role of the underlying measure $\sigma$. In particular, if
$f \equiv 1$ we have
 $K \bold 1 = K^\sigma \bold 1 = K \sigma$. Similarly, for any Borel
measure $\omega$ on $X$, we set
$$K \omega (x) = K^\omega  \bold 1 (x) =
\int_X K(x,y) \, d \omega (y),$$
which hopefully does not lead to any confusion.

The problems (1.1),
 (1.2) and many
similar problems can be transformed into
``superlinear'' integral
equation of the type
$$ u(x) = K u^q(x)+f(x)\qquad (\sigma-\text{a.e.})\tag 1.5$$
where $f\ge 0$ is given and we require a nonnegative solution.
Let us denote by $S_{q,K}$ the set of all $f$ such that (1.5) has a
solution (i.e. there is a measurable  $u \ge 0$ satisfying (1.5)).
 This
problem
has been considered in several places in the literature
(cf. \cite{3},
 \cite{6}, \cite{29}, \cite{32}, \cite{52} ).
In the examples we wish to consider the kernel $K$ has an additional
property which we term the {\it quasi-metric assumption} i.e. $K$
is symmetric, $K(x,y)>0$ for all $x,y$, and there is a constant $\kappa$
such that for
$x,y,z\in X$ we have
 $$ \frac{1}{K(x,y)}\le
\kappa\left(\frac{1}{K(x,z)}+\frac{1}{K(z,y)}\right).\tag 1.6 $$
Then we can introduce a quasi-metric structure via
$\rho(x,y)=K(x,y)^{-1}.$  Note however that we do not require
$\rho(x,x)=0.$ If we define  the $\rho$-ball $B_a(x)=\{y:\rho(x,y)\le
a\}$
then we can write $$K^\sigma f(x)= K \nu (x) = \int_0^{\infty}
\frac{|B_t(x)|_{\nu}}{t^2}dt
\tag 1.7$$
where  $d \nu = f \, d \sigma$, $f \ge 0$.

We remark that in our main results we are able to avoid the usual
assumption
of the theory of homogeneous spaces  in the sense
of Coifman and Weiss
\cite{11} that there exists a doubling
measure with respect to $\rho$ on $X$.
Under that
assumption  our results are applicable to equations with
generalized fractional integral operators of \cite{20}.

Returning to (1.5) we seek to characterize those functions $f$ so that
for some $\epsilon>0$ we have $\epsilon f\in S_{q,K}$ i.e. there is a
solution to the equation $$ u=Ku^q+\epsilon f.$$
In this context it is natural to introduce the {\it solution space
$\Cal Z_{q,K}$} of all measurable functions $f$ such that for some
$\epsilon>0$ we have $\epsilon |f|\in S_{q,K}.$  If we have the
quasi-metric condition (1.6) then $\Cal Z_{q,K}$ either reduces to
$\{0\}$ or is a Banach function space on $(X,\sigma)$ with associated
norm $\|f\|_{\Cal Z}=\inf\{\lambda > 0: f\in\lambda S_{q,K}\}$,
and all solutions of the equation belong to $\Cal Z_{q,K}$.

In Section 2, which consists mainly of background material, we develop a
general theory of the solution space
$\Cal Z=
\Cal
Z_{q,T}$ associated to the equation $u=Tu^q+f$ where $T$ is an arbitrary
positive operator on the space of measurable functions. Most
 results on nonlinear  operator equations of this type we have found
 in the literature (see e.g.
\cite{29}, \cite{32})
 contain only sufficient conditions for solvability in
particular
 function spaces  which are not  intrinsically related
to the equation. To bridge this gap,  we use  techniques borrowed from
the Banach lattice theory
and in particular some facts from the Nikishin-Maurey theory
(see \cite{31}, \cite{33}, \cite{39}) to
identify $\Cal
Z_{q,T}$ as a Banach function space naturally associated with the
problem.

 Under mild assumptions,  we demonstrate a number of
different characterizations of this space.  For example we show that
$\Cal Z$ is invariant under the mapping $f \to \Cal Af=Tf^q$
and that $0\le f\in \Cal Z$ if and only if
$$\limsup_{n\to\infty} \, (\Cal
A^nf)^{1/q^n}\in L^{\infty} (\sigma).$$

We also show that earlier  results
 of P. Baras and M. Pierre \cite{6}  can be recast as
identifying
the predual space of $\Cal Z.$  The K\"othe dual $\Cal Z'$ of $\Cal Z$ is
shown to be the predual and its norm is given by the formula:
$$ \|g\|_{\Cal
Z'}=pq^{p-1}\inf\left\{\int\frac{h^p}{(T^*h)^{p-1}}d\sigma:\
 h \ge |g|\right\}. \tag 1.8$$
Here $T^*$ is the adjoint of the operator
 $T$, $1/p + 1/q =1$, and the equation $u = \Cal A u + f$
is solvable if and only if $f \ge 0$ lies in the unit ball
of $\Cal Z$.
This is quite similar to a  result proved by Baras and
Pierre, who consider more general convex functions than $x\to x^q.$
It
should be noted that our result applies to general positive operators
(not simply operators defined by kernels) and that our method of proof is
quite different from that of Baras and Pierre, depending only on quite
simple duality arguments.
Although this characterization of the predual space $\Cal Z'$ seems
to us both important and interesting, it is not used to establish the
main results of the paper.

At the same time in Section 2 we introduce the
Banach function space $\Cal W_{p, T}$ of $L^p$-weights associated
with the corresponding weighted  norm inequalities for future use.
Its $q$-convexification which contains $\Cal Z_{q, T}$
plays an important
role in the sequel and is intimately related to the general
theory of
$L^p$-capacities developed by N. G. Meyers (see \cite{2}).

In Section 3, we carry out
the study of the solution space further for the
case of an operator defined by a kernel $K$ satisfying the quasi-metric
assumption. Our main technical tool here
is a decomposition of the operator into its upper and lower part,
the latter being almost constant on quasi-metric balls (Proposition 3.4).
We show that if $\omega$ is a $\sigma-$finite Borel measure
on $X$ and
$$K\omega(x)=  K^\omega \bold 1 (x) = \int K(x,y) \, d\omega(y)$$
then $K\omega\in \Cal Z$ implies the inequality:
$$\sup\Sb x \in X\\a > 0\endSb \,
\left(\int_0^a\frac{|B_t(x)|_{\sigma}}{t^2}dt\right)
\left(\int_a^{\infty}\frac{|B_t(x)|_{\omega}}{t^2}dt\right)^{q/p}
\ < \infty
 \tag
1.9$$
which reflects the interplay between
 the local regularity of $\sigma$ and the behavior
of the ``tails'' of $K \omega$. It is reminiscent of the Kato class
for the potential $\sigma$ as $q \to 1$ (see \cite{13}).
We call (1.9) the
{\it infinitesimal inequality\/} because
of the method of the proof which boils down to  careful pointwise
estimates
of  $(\Cal A^nf)^{1/q^n}$ as $n \to \infty$; they involve  sharp
 constants depending on $q$ and $\kappa$.
 In particular if $f\in\Cal
Z$ then the infinitesimal inequality holds for the measure
$d\omega=f^q \, d\sigma.$

In Section 4, we then turn to the relationship between solutions of (1.2)
and weighted norm inequalities of the type:
$$ \int (K^\sigma g)^p \, d\omega \le C\int g^p \, d\sigma\tag 1.10$$
whenever $g \ge 0$
and $g\in L^p(\sigma).$
We prove that if $K\omega  \in\Cal Z$ then the
weighted norm inequality
(1.10) holds, and further that  (1.10) is also implied by the inequality
$$ \int (K^\sigma g)^p \, (K\omega)^q \, d \sigma
\le C\int g^p \, d\sigma.$$

In general $K\omega\in\Cal Z$ is equivalent to a sequence of weighted
norm inequalities
$$ \int (K^\sigma g)^p \,  d \omega_j
\le C_j \int g^p \, d\sigma, \quad j = 1, 2, \ldots,$$
where $d \omega_0 = d \omega$ and $d \omega_{j+1} = (K\omega_j)^q \, d \sigma$,
with a good control of the imbedding
constants: one should require
$\limsup_{j \to \infty} \, C_j^{1/q^j} < \infty$
 (Theorem 4.7); however this result is a little
unwieldy.

Our main result of the section (Theorems 4.8 and 4.9)
is that
$K\omega\in
\Cal Z$ if
and only if both the infinitesimal inequality (1.9) and the weighted
norm
inequality (1.10) hold, which is also equivalent to the pointwise
inequality $K(K \omega)^q \le C \, K \omega$.
Also in these results it is possible
to replace
(1.10) by its weak-type analogue, or even  the well-known
testing condition of Sawyer type (see \cite {44}, \cite {45}, \cite {46})
$$\int_{B} (K^\omega \chi_B)^q \, d \sigma \le C \, |B|_\omega
 \tag 1.11$$
for all quasi-metric
balls $B=B_a (x)$. This leads to the following theorem.

\proclaim{Theorem 1.2} Let $\omega, \sigma \in \Cal M_+(X)$ and let
$K f = K^\sigma f$.
Then the  following are  equivalent:
\newline  (1) $K \omega \in \Cal Z_{q, K}$ i.e. the equation
$u = K u^q + \epsilon \, K \omega$ has a solution for some $\epsilon >0$.
\newline  (2) $\omega$ satisfies both the infinitesimal inequality
(1.9) and weighted
norm inequality (1.10).
\newline (3) $\omega$
satisfies both the infinitesimal inequality
(1.9) and  testing inequality (1.11).\newline
 (4) There exists a constant
  $C$   so that   $K(K\omega)^q\le
C \, K\omega < \infty$ $\sigma$-a.e.\endproclaim

It is easy to see that if in (4) the constant
 $C \le q^{-1} p^{1-q}$, then the equation $u = K u^q +
  K \omega$
 has  a solution $u$
so that $ K \omega \le u \le
p \,  K \omega$. This remark means that Theorem 1.2 makes it possible
to characterize the existence of positive solutions in $L^r$ spaces
(or any other ideal space).

Similar results hold for all
$0\le f\in\Cal Z$.

\proclaim{ Theorem 1.3} Let $0 \le f < \infty$ $\sigma$-a.e.
 and $d \omega
= f^q \, d \sigma$. Then the equation $u = K u^q + \epsilon \, f$ has a
solution for some $\epsilon > 0$ if and only if
any one of statements (2)-(4) of
Theorem 1.2 holds.
\endproclaim
 It follows from the general estimates of Section 2
that  if  $d \omega
= f^q \, d \sigma$ and $K(K \omega)^q\le  q^{-q} p^{q(1-q)}
 \, K \omega < \infty$ as in statement (4) of Theorem 1.3
 then  the
equation $u = K u^q +  f$ has
 a solution so that $f + K f^q\le u \le f + p^q
\, K f^q$.

Note that in the proof of Theorems  1.2 and 1.3 we do not use any
known two weight theorems (cf. \cite{44}, \cite{45}); our proofs are
self-contained and shed new light on the role of weighted norm
inequalities and testing conditions of
Sawyer type.

In Section 5 we introduce  notions of capacity associated to the kernel
and then study conditions on the kernel $K$ and the measure $\sigma$
under which the statement
$K\omega\in\Cal Z$
becomes equivalent to the weighted norm inequality (1.10).  This requires
that (1.10) implies (1.9).  In fact we give two theorems of this
type.  In Theorems 5.6 and 5.9 we give conditions on the kernel $K$ so
that $K\omega\in\Cal Z$ if (1.10) holds or if a weaker capacitary
condition holds:
$$ |E|_{\omega}\le C \, \Cap E\tag 1.12$$
where $$\Cap E = \inf\left\{\int g^p \, d\sigma:\ 0\le g,\ K g\ge
\chi_E\right\}.$$  It is also possible to replace this condition by a
Sawyer-type testing condition. The hypothesis of Theorem
5.6 (which can
be shown to be necessary for the conclusion under certain mild
assumptions on
$K$) is that for some constant $C$ and every $x\in X$ and $a>0$ we have
$$ \int_0^a\frac{|B_t(x)|_{\sigma}}{t^2}dt \le
Ca^{q-1}\int_a^{\infty}\frac{|B_t(x)|_{\sigma}}{t^{1+q}}dt
< \infty.\tag 1.13$$
Roughly speaking this condition implies that the behavior of the kernel
at infinity dominates the behavior locally.  The hypothesis of Theorem
5.9 replaces this by the assumption that for some $C$ and every $x\in
X$,
$a>0$ we have both:
$$ \int_0^{2a}\frac{|B_t(x)|_\sigma}{t^2}dt\le
C\int_0^a\frac{|B_t(x)|_{\sigma}}{t^2}dt \tag 1.14$$ and
$$ \sup_{y \in B_a(x)} \, \int_0^{a}\frac{|B_t(y)|_\sigma}{t^2}dt\le
C\int_0^a\frac{|B_t(x)|_{\sigma}}{t^2}dt. \tag 1.15$$
Conditions (1.14)
and (1.15) essentially are assumptions that measure $\sigma$ is close to
being invariant for the kernel $K$. For convenience we
 state these results as
the following
theorem.

\proclaim{Theorem 1.4} Let  $K$  be  a  quasi-metric kernel.
Assume that there is a constant $C$ so that for every $x\in X$  and
$a>0,$   either (1.13),    or both (1.14) and (1.15) hold.
 Then $\Cal Z_{q, K}\neq \{0\}$  and
the  following  statements  are  equivalent:    \newline (1)
  $\omega$
satisfies  the
weighted   norm   inequality (1.10).\newline  (2)
$\omega$ satisfies  the capacity  condition $|E|_{\omega}\le
\Cap E$ for all Borel sets $E.$ (Equivalently the  weak-type
analogue of (1.10)  holds.)
\newline (3)  $\omega$ satisfies
the    testing     condition    (1.11).
\newline  (4)  $K\omega\in \Cal Z_{q,K}$  i.e.  for   some
$\epsilon>0$  there  is  a  solution  $u$  of  the  equation
$u=Ku^q+\epsilon \, K\omega.$ \newline (5) There is a  constant
$C$ so that $K(K\omega)^q\le C \, K\omega.$ \endproclaim

The proof  of Theorem 1.4 involves a quantity which under
mild assumptions
is equivalent
to the capacity of a ball. The  two-sided estimates of
$\Cap B$ are obtained in Theorem 5.4
without any restrictions
on the underlying measure $\sigma$ for a wide class of $K$.
 (This is a generalization
of  D. Adams's formula for the weighted capacity of a ball
proved in  \cite{1} in the case of Riesz potentials
for $\sigma \in A_\infty$; an upper estimate for arbitrary $\sigma$
can be found in \cite{51}.)
This should
be compared to a similar estimate of  $||\chi_B||_{\Cal Z_{q, K}}$
 in Section 4 which clarifies the role of the
infinitesimal inequality.

In Sections 6 and 7 we give some applications of our ideas to concrete
problems of the type introduced earlier.  Section 6 is devoted to
convolution operators on $\bold R^n.$  In particular we translate our
results for the Riesz potential of order $\alpha$ i.e.
$I_{\alpha}=(-\Delta)^{-\alpha/2}.$  The results developed in Sections
3, 4 and 5 can be translated directly to this setting taking into account
the change between ordinary Euclidean distance and the quasi-distance
induced by the quasi-metric $\rho$  (see Theorem 6.2 below).  We also
apply our results to the Poisson kernel to derive an extension of the
recent result of Treil and Volberg \cite{50} mentioned above.

In Section 7 we return to the problem which motivated this research, i.e.
the equations (1.1) and (1.2). We  consider more general
differential operators and inhomogeneous boundary
conditions.  If $G$ denotes the Green's kernel for the Laplacian $\Delta$
on $\Omega$ then (1.1) can be transformed to the equation:
$$ u= G(v \, u^q) + Gw. \tag 1.16$$
If we let $d\sigma=v(x)dx$ then we can consider this equation as being in
the form of (1.5). It is easy to see that
 the Green's kernel $G$ fails to satisfy the
quasi-metric assumption (1.6) in general, even for the simplest domains
(e.g. the Euclidean ball or the half-space).  However $G$ does satisfy
the so-called $3G-$inequality (see \cite{12}) i.e.
$$ \frac{G(x,y) \, G(y,z)}
{G(x,z)}\le C(|x-y|^{2-n}+|y-z|^{2-n})\tag 1.17$$
In \cite{7} Theorem 3.6 it is claimed that for Lipschitz domains one
can replace the right-hand side by $G(x,y)+G(y,z)$ which would establish
(1.6) for the Green's kernel but there is an error in the proof.

However if the boundary $\partial \Omega$ is smooth enough, the problem
can be transformed to meet the condition (1.6).  If $\partial \Omega$ is
$C^{1,1}$ then the \Naim kernel (introduced in an equivalent form
by Linda \Naim \cite{38}
in the theory related to Martin's kernels)
is defined by
$$N(x,y)=\frac{G(x,y)}{\delta(x)\delta(y)}\tag 1.18$$ where
$\delta(x)=d(x,\partial \Omega)$ is the distance to the boundary.  We
show that $N$ does indeed satisfy (1.6), and this enables us to
transform
(1.16) to an equation for which our general theory is applicable.
We observe  that this inequality is sharper than
(1.17) and gives the
right estimates of the Green's kernel at the boundary.

 These
methods can be applied to more general second-order non-divergence
uniformly elliptic differential operators $L$ with regular (bounded
H\"older-continuous) coefficients in place of
the Laplacian.  Under these
assumptions, it follows from the well-known  estimates of the
Green's kernels (see \cite{55}, \cite{56},
 \cite{5}, \cite{25}, and the discussion in Section 7)
that the corresponding \Naim kernel given by (1.18) satisfies
condition (1.6).

By using this method we are able to give  very general results on the
solvability of the equations (1.1) and (1.2) as well as a
characterization of trace inequalities for Green's potentials.
 In Theorem 7.5 we show that
the equation $$ -\Delta u =\sigma \, u^q + \epsilon \,
\omega \tag 1.19$$ with
$u\ge 0$ and $u=0$ on $\partial \Omega$ has a solution for some
$\epsilon>0$ if and only if for some constant $C$ we have
$$G [\sigma \, (G\omega)^q] \le C \, G\omega, \tag 1.20$$
 where
$G\omega(x)=\int_\Omega G(x,y) \, d\omega(y)$. Here $\sigma$ and $\omega$
are arbitrary nonnegative measurable
functions (or positive measures)  on $\Omega$.
We also give equivalent formulations in geometric terms involving the
infinitesimal inequality and the
testing inequality quoted above.

Finally we return to the Adams-Pierre theorem (Theorem 1.1) and show that
our methods in this situation give a complete solution (up to a constant)
to the  problem of the existence of positive solutions.
This corresponds to the case $\sigma\equiv 1$, but similar results are
proved for arbitrary  $\sigma$.
 We can then
apply Theorem 1.4 and show that the solvability of (1.19) can be
characterized in terms of a capacitary condition.  We are thus able to
remove the assumption that $\omega$ is compactly supported in $\Omega$.

Let $$\Cap E=\inf \, \left \{\int_\Omega g^p \, \delta(x)^{1-p} \, dx:
 \ Gg(x)\ge \delta(x) \, \chi_E(x),\ g\ge 0 \right \}$$
for any set $E \subset \Omega$.

\proclaim{Theorem 1.5} Let $\omega \in \Cal M_+(\Omega)$ and
$\sigma\equiv 1$.
Then the   Dirichlet problem (1.19) has a solution
for some $\epsilon >0$ if and only if
there is a constant $C$ so that
$$ \int_E \delta(x) \, d\omega(x) \le C \, \Cap E\tag 1.21$$
for every compact
set $E$.

Moreover, (1.21) is equivalent to the pointwise condition
(1.20).\endproclaim

As was mentioned above, if in (1.20)
the constant $C \le q^{-1} p^{1-q}$, then for $0<\epsilon\le 1$
  (1.19) has  a solution $u$
so that $ G \omega \le u \le
p \,  G \omega$.

In the case when $\omega$ is compactly supported the capacitary
characterization of Theorem 1.5 reduces to the
Adams-Pierre theorem since the capacity defined above can be then shown
equivalent to the nonlinear Newtonian capacity associated with the
Sobolev space $W^{2,p}(\bold R^n)$ used in \cite{3}.

After this paper was finished we learned that
H. Brezis and X. Cabre \cite{9}
considered very recently  another special case of (1.19)
where  the inhomogeneous term $\omega$ is a positive uniformly
bounded function
on $\Omega$. \footnote{See also the Addendum at the end of this paper.}
In particular,
they proved in a different way the following
``nonsolvability'' result (for bounded domains $\Omega$ with
smooth boundary):  (1.19)
with bounded $\omega$  has no solutions
unless $G (\delta^q \, \sigma) \in L^\infty (\Omega)$.
 Note that
in this case $G \omega (x) \asymp \delta (x)$, and our characterization
  (1.20)
 boils down to  a sharper necessary and sufficient condition
$$ G (\delta^q \, \sigma) (x) \le C \, \delta (x), \quad x \in
 \Omega.\tag 1.22$$
By Theorem 7.5 (see Sec. 7), it follows that (1.22) also
 characterizes completely the solvability of
 (1.19) with  $\omega \in L^\infty (\Omega)$
 for more general  uniformly
elliptic second order differential operators in place of the Laplacian.

All the results in Section 7 apply to much more general situations as
explained therein.
For instance,  the only property of
the differential operator with Green's function $G$ which is
important for us  is the fact that the kernel
$$N(x, y) = s(x) \, G(x, y) \, s(y)^{-1}$$
has
the quasi-metric property for some  weight
function $s > 0$; then  the pointwise
condition (1.20) which is invariant
under this transformation of the kernel
characterizes
the solvability of the corresponding Dirichlet problem.
We conjecture that this
holds true for a wide class of
differential operators with bounded measurable
coefficients and  non-smooth domains  $\Omega$.
There are many other potential applications of these
ideas which we plan to explore in future work.

We would like to express our thanks to our colleague Zhongxin Zhao for
his helpful comments concerning the $3G$-inequalities.

\vskip10pt\heading {2. Superlinear problems and
related function
spaces}\endheading
\vskip10pt

In this section we will introduce certain Banach function
 spaces which
will play an important role in the later sections of the
paper.  We will
also give an alternative approach to some results of Baras
and Pierre \cite{6} on
the solvability of the ``superlinear" problem $u=Tu^q+f$
discussed in the
introduction.

Let  $X$  be  a  metric  space  and  suppose  $\sigma$  is a
$\sigma-$finite measure on  $X.$ We denote  by $L^0(\sigma)$
the space of all (equivalence classes) of real-valued  Borel
functions  on  $X$.    The  topology of $L^0(\sigma)$ is the
usual topology of convergence  in measure on sets  of finite
measure. We use $L^0_+$ to denote the positive cone
$\{f: f \in  L^0, \ f \ge 0\}.$

Let us  first state  two fundamental  results which  will be
used in the sequel:

\proclaim{Theorem  2.1}Let  $H$  be  a closed bounded convex
subset of $L^0_+.$  Then: \newline (1) (Nikishin  \cite{39})
There exists  a weight  function $w\in  L^0_+$ with  $w>0$
a.e.  such  that  $\sup_{f\in  H}\int   fw\,d\sigma<\infty.$
\newline (2) (Komlos \cite{28})  If $(f_n)$ is a  sequence in
$H$ then $(f_n)$ has a subsequence $(g_n)$ such that  $\frac
1n(g_1+\cdots+g_n)$ is a.e. convergent.  \endproclaim

Nikishin's  result  is  \cite{39}  Theorem  4,  or see Maurey
\cite{33} Th\'eor\`eme 13.   We observe that Komlos's  result
is usually  stated for  bounded sequences  in $L^1$  but our
statement follows from the  usual Komlos theorem in  view of
(1).  Also note that  the conclusion of Komlos's theorem  is
much stronger than we have stated here:  one can ensure that
every subsequence of $(g_n)$ is Cesaro convergent a.e.  to
some fixed $h.$ The version  we will use can be  established
by much more simple  means.  In fact  we do not need  Cesaro
means;  any   suitable  sequence   of  convex   combinations
suffices.

We say  that a  convex subset  $H$ of  $L^0_+$ is {\it
solid} if
$f\in H$ and $0\le  g\le f$ a.e. implies that  $g\in
H.$ We  will also  say that  $H$ is  {\it nondegenerate}
 if  there
exists $f\in H$  with $f>0$ a.e.   Our next  two results are
surely well-known to specialists but we know of no reference
where they are established in exactly this form.

\proclaim{Proposition 2.2}Let $H$  be a solid  convex subset
of $L^0_+.$  Then $H$  is closed  if and  only if whenever
$f_n\uparrow f$ a.e. with $f_n\in H$ for all $n,$ then $f\in
H.$\endproclaim

\demo{Proof} Suppose $h_n\in
H$ and $h_n\to f$ a.e.   Then by applying Egoroff's  theorem
we can find an increasing sequence of Borel subsets $E_m$ of
$X$  so  that  $h_n\to  f$  uniformly  on  each $E_m$, $f\ge
m^{-1}$ on $E_m$ and $\cup  E_m=\supp f.$ Now for fixed  $m$
we    can    find    $\epsilon_n\downarrow    0$   so   that
$(1-\epsilon_n)f\le h_n\le  f$ on  $E_m$.   It follows  that
each $(1-\epsilon_n)f\chi_{E_m}\in H$ and hence we can apply
the hypothesis twice to obtain $f\in H.$\qed\enddemo

If $H$ is a solid closed bounded convex subset of  $L^0_+$
then we can define a Banach function space $\Cal X=\Cal X_H$
associated  to  $H.$  Precisely  if  $f\in  L^0$  we  set $$
\|f\|_{\Cal X}=\inf \,
\{\lambda:  |f|\in \lambda C\}$$ and then
$\Cal X=\{f:\|f\|_{\Cal X}<\infty\}.$  It is easily  checked
that $\Cal X$ is a Banach lattice continuously embedded into
$L^0.$ Notice that the norm $\|\,\|_{\Cal X}$ has the  Fatou
property   i.e.      $0\le   f_n\uparrow  f$  implies  $0\le
\|f_n\|_X\uparrow \|f\|_{\Cal X}.$  We shall say  that $\Cal
X$  is  nondegenerate  if  it  contains  a strictly positive
function:    this  is  easily  seen  to  be  equivalent   to
nondegeneracy of $H.$

For  any  convex  subset  $H$  of  $L^0_+$  we  can define
$$H'=\{g\in  L^0_+:\int  fg\,d\sigma  \le  1,\ \forall f\in
H\}.$$
It is clear that $H'$ is closed solid and convex.   If
$H$ is bounded then by  Nikishin's Theorem 2.1 (1) $H'$  is
nondegenerate; if  $H$ is nondegenerate  then $H'$  is bounded.
When $H$  is nondegenerate  bounded  convex  and solid  then
$\Cal X_{H'}$ is  K\"othe dual space
$$\Cal X'=\{f:   \ \int
|f| \, |g| \, d\sigma<\infty \quad \forall \  g\in\Cal X\}$$
 of  $\Cal X,$
equipped with the dual  norm.  The following  observation is
important  for  later  applications  so  we  display it as a
lemma:

\proclaim{Lemma  2.3}  If  $H$  is  a  nondegenerate bounded
closed solid convex set  then $H=H''$ and $\Cal  X=\Cal X''$
with equality of norms.\endproclaim

\demo{Proof}Using  Nikishin's  theorem  again  there  exists
$u\in  H'$  with  $u>0$  a.e.     Then  $H$  is  closed   in
$L^1(u\,d\sigma)$ or equivalently  $u H$ is closed  in $L^1$.
Hence  $uH=\{f\ge  0:\int  fg\le  1  \ \forall \ 0 \le g\in
u^{-1} H' \cap      L^{\infty}\}$ and the lemma
follows.\qed\enddemo

Next we define a generalized notion of a positive
 operator.
Let $L^{0,\#}_+$ denote the space of (equivalence classes
of) measurable  functions $f:X\to  [0,\infty].$ We  define a
{\it positive   operator}   to   be   a   map
$T:L^{0,\#}_+\to
L^{0,\#}_+$ such  that:   \newline (P1)  $T(\alpha f+\beta
g)=\alpha Tf+\beta Tg$  if $\alpha,\beta\ge 0,$  and $f,g\in
L^{0,\#}_+.$ \newline  (P2) If  $0\le f_n\uparrow  f$ a.e.
then $Tf_n\uparrow Tf$ a.e.  \newline (P3) There exists  $f$
with $f>0$ a.e. such that $Tf<\infty$ a.e.

We say that $T$ is {\it strictly positive} if in addition
$T f = 0$ a.e. implies $f=0$ a.e.

In  all  our  applications  $T$  will  be  given by a kernel
function $K$.  That is, we will be given a nonnegative Borel
function $K$ on  $X\times X$ and  then $$ Tf(x)=Kf(x)  =\int
K(x,y)f(y)\,d\sigma(y).$$ Conditions (P1) and (P2) are  then
obvious and (P3) is a  condition on the kernel.

 The {\it domain}
of $T$  is defined  by $\Cal  D(T)=\{f:   T|f|<\infty \text{
a.e.}\}.$ It is clear that $T$ extends to a positive  linear
operator $T:\Cal D(T)\to L^0.$

Given $T$ we can construct a formal adjoint
 $T^*$ such  that
if  $f,g\in  L^{0,\#}_+$  then  $\int  (Tf)g\,d\sigma=\int
(T^*g)f \,d\sigma.$ This can be done using the Radon-Nikodym
theorem and it is not difficult to check that $T^*$ is  also
a positive  operator.   To verify  (P3) for  $T^*$ one needs
only observe that if $f$ is chosen as in (P3) for $T$ and if
$\int u (Tf)<\infty$ then $T^*u<\infty$ a.e.

We first discuss  weighted norm inequalities  for $T$.   Fix
$1<p<\infty$, and let $\frac1p+\frac1q=1.$

Let $H$ be the set of $f\in L^0_+$ such that there  exists
$0\le g\in L^p(\sigma)$ with  $\|g\|_p\le 1$ and $f\le  (Tg)
^p.$

\proclaim {Lemma 2.4} $H$ is a solid closed convex set.  $H$
is bounded  if and  only if  $L^p\subset \Cal  D(T).$ $H$ is
nondegenerate if  and only  if whenever  $f\in L^0_+$ with
$T^*f=0$ a.e.  then $f=0$  a.e.   (i.e.   $T^*$ is  strictly
positive).\endproclaim

\demo{Proof}Let us prove convexity.   If $f_1,f_2\in H$  and
$0\le t\le 1$ we first find $0\le g_j$ with $Tg_j^p\ge  f_j$
and    $\|g_j\|_p\le    1$    for    $j=1,2.$    Then    let
$g=(tg_1^p+(1-t)g_2^p)^{1/p}.$  Then   $\|g\|_p\le  1$   and
$Tg^p\ge  tf_1+(1-t)f_2$.    $H$  is  easily seen to be also
solid.

We  next  check  $H$  is  closed.    Suppose  $f_n\in H$ and
$f_n\uparrow f$  a.e.   Then there  exist $0\le  g_n\in L^p$
with $\|g_n\|_p\le 1$ and $(Tg_n)^p\ge f_n.$ Since $L^p$  is
reflexive when $1<p<\infty$ we  can find a weak  limit point
$h$     of     $(g_n)$     and     a     sequence      $g'_n
\in\co\{g_n,g_{n+1},\ldots\}$  such  that  $g'_n\to  g$   in
$L^p$-norm.  Then $Tg'_n \to Th$ in $L^0.$ Now
$$(Tg'_n)^p\ge
\inf \, \{(Tg_k)^p:\ k\ge n\}\ge f_n$$
 so that $(Th)^p\ge f.$  It
follows that $f\in H$ and by Proposition 2.2 the set $H$  is
closed.

Clearly if $0\le  f\in L^p$ and  $E=\{x:Tf(x)=\infty\}$ then
$\alpha\chi_E\in H$  for every  $\alpha>0.$ Hence  if $H$ is
bounded  $ L^p \subset \Cal  D(T).$  Conversely assume $ L^p
\subset \Cal  D(T).$  To show $H$ is bounded, let us suppose
$2^nf_n\in  H$.    The  above  argument  used to show $H$ is
closed shows the existence of  $0\le h\in L^p$ with $Th  \ge
\sum_{n=1}^{\infty}f_n$ so that $f_n\to 0$ a.e.

If $H$ is  nondegenerate then there  exists $0\le f\in  L^p$
with   $Tf>0$   a.e.      If   $T^*h=0$   a.e.   then  $\int
h(Tf)d\sigma=0$ so that $h=0$ a.e.

Conversely  let  us  show  $H$  is  nondegenerate, under the
hypothesis  that  $T^*f=0$  a.e.  implies  $f=0$  a.e.   Let
$(f_n)$ be a sequence in  dense in the positive quadrant  of
the      unit      ball       of      $L^p$.             Let
$g=(\sum_{n=1}^{\infty}2^{-n}f_n^p)^{1/p}.$              Let
$E=\{Tg=0\}.$ Then  $\chi_ETf_n=0$ a.e.  for all  $n$ whence
$T^*\chi_E=0$ so that $|E|_{\sigma}=0.$ \qed \enddemo

\proclaim{Definition}The   space   $\Cal   V=\Cal   V_{p,T}$
consists  of  all  $f\in  L^0$  such that there exists $h\in
L_{p,+}$ with $|f|\le (Th)^p$ a.e.  We define
$$ \|f\|_{\Cal V}=\inf \, \{\int h^pd\sigma:\
 0\le h,\  |f|\le (Th)^p\}.$$
 If $L^p\subset \Cal D(T)$  then $\Cal V$  is a Banach
function space.  $\Cal V$ is a nondegenerate Banach function
space if
in addition $T$ is strictly positive.  \endproclaim

\proclaim{Definition}The  space  $\Cal  W=\Cal  W_{p,T}$  of
$L^p$-weights  for  $T$  is  defined  to  be  the  space  of
functions $f\in L^0$ such that for some $\gamma$ we have
$$\int |f| \, |Tg|^p \,
 d\sigma  \le \gamma \int |g|^p \, d\sigma,$$
for all $g \in L^p$. We
define $\|f\|_{\Cal  W}$ to  be the  least constant $\gamma$
for which the preceding inequality holds.  \endproclaim

It is clear that if we have both $L^p\subset \Cal D(T)$  and
$T$ strictly positive then $\Cal  W$ is the K\"othe dual  of
$\Cal V$, and is thus  a Banach function space for  the norm
$\|f\|_{\Cal W}.$ Note that this implies $\Cal W'=\Cal V$ as
we observed above in Lemma 2.3.

We  will  also  introduce  the  $q$-convexification of $\Cal
W_{p,T}$ say $\Cal Y_{p,T}=\Cal W_{p,T}^{1/q}.$

\proclaim{Definition}$\Cal  Y_{p,T}$  is  the  space
of $w\in L^0$
such that $|w|^q\in\Cal W_{p,T}$ with the associated norm $$
\|w\|_{\Cal  Y}  =  \||w|^q\|_{\Cal  W}^{1/q}.$$  This  is a
Banach function  space whenever  $\Cal W_{p,T}$  is a Banach
function space.\endproclaim

Note that $\|w\|_{\Cal  Y}$ is the  least constant
$\gamma$
for which the inequality
$$
\int w^q \,  |Tg|^p \, d \sigma
\le \gamma^q \int |g|^p \, d\sigma$$
holds for all $g \in L^p$.

We next  consider the  nonlinear equation  $$ u=Tu^q  +f\tag
2.1$$  where  $T$  is,  as  before,  a positive operator and
$1<q<\infty.$  We  suppose  $f\ge  0$  and  seek  a positive
solution  $u\in  L^0.$  With  the  restriction that $T$ is a
kernel operator  but for  more general  convex functions  in
place of $u\to u^q$  this problem was previously  considered
by Baras and Pierre \cite{6}.

Let us start with some very elementary observations.
We  denote by $\Cal A: L^{0,\#}_+\to L^{0,\#}_+ $ the
nonlinear map $\Cal A f=Tf^q$ and rewrite (2.1) as
$$u = \Cal A u + f.$$
Note that $\Cal A$ is a convex operator:
$\Cal A (t \, f + (1-t) \, g) \le t \,
\Cal A f + (1-t) \,
\Cal A g$
for all $f, g \ge 0$ and $0<t<1$.

\proclaim{Proposition  2.5}Suppose  $f\in L^0_+.$ Define
$u_0=0$ and  then $u_n= \Cal A u_{n-1} +f$  for $n\ge 1.$
Then the following  are  equivalent: \newline  (1)
There exists
$v\in L^0_+$ with $v= \Cal A v+f.$\newline (2)
The sequence  $(u_n)$
is bounded in $L^0.$\newline (3) $\displaystyle{\sup_{n\ge 0}
 u_n}<\infty$ a.e. \newline (4)
There   exists   $ w\in   L^0_+$   with  $w\ge
\Cal A w+f.$
\endproclaim

\demo{Proof}Of  course  (1)  implies  (4).   (4) implies (3)
since  we  will  have  $0\le  u_n\le  w$  for  all  $n.$ (3)
trivially implies (2).  If  (2) holds then since $(u_n)$  is
increasing   we   have    $v=\displaystyle{\sup_{n\ge 0} u_n}
\in   L^0$    solving
(2.1).\qed\enddemo

\proclaim{Definition}We define  $S=S_{q,T}$ the  set of  all
$f\ge 0$ so that (2.1) has a solution.\endproclaim

Notice  that  the  iterative  procedure  of  Proposition 2.5
yields a {\it minimal solution} $u$ of (2.1) corresponding to
 each
$f\in S$ and that the map $f\to u$ (where $u$ is the minimal
solution) is monotone.

\proclaim{Proposition  2.6}The  set  $S$  is  a solid convex
set.\endproclaim

\demo{Proof}This follows from 2.5.  The only part  requiring
proof is convexity.   If $f,g\in  S$ and $0\le  t\le 1$ then
there  exist  $u,v$  with  $u=\Cal A u+f,\  v= \Cal A v+g.$
But then using the convexity of $\Cal A$ we have
$tu+(1-t)v\ge  \Cal A (tu+(1-t)v)+tf+(1-t)g$     so  that
$t f + (1-t) g \in S.$\qed\enddemo

\proclaim{Proposition 2.7}(1) Suppose $f\in S$.  Then
$\displaystyle{\sup_{n\ge 0}
\Cal A^n f<\infty}$  a.e.\newline (2)  Suppose $\co\{\Cal A^n
f:\ n\ge 0\}$ is  bounded in $L^0.$ Then  $p^{-1}q^{1-p}f\in
S.$ In particular,  if $\displaystyle{\sup_{n\ge 0}
\Cal A^n f<\infty} <\infty$ a.e.
 (e.g. if for some $n$ we have $\Cal A^nf\le \Cal A^{n-1}f$)
 then $p^{-1}q^{1-p}f\in S.$\newline
(3) If $\Cal A f \le q^{-1} p^{1-q}  \, f$
then (2.1)
has a solution $u$ such that
$f \le u \le p \, f.$\newline
(4) If $\Cal A^2 f\le q^{-q} p^{q (1-q)}  \,
\Cal A f$
then (2.1) has a solution $u$ such that
$$f + \Cal A f \le u \le f + p^q \, \Cal A f.$$
\endproclaim

\demo{Proof}(1)  is  immediate.    If  $u=\Cal A u +f$
then $\Cal A^n f\le \Cal A^n u\le u.$

For (2)  let $G$  be the  set $\co\{\Cal A^n f:\ n\ge  0\}$.
Define $u_0=0$ and then  $u_n=\Cal A u_{n-1}+p^{-1}q^{1-p}f$
for $n\ge  1.$ We  will show  by induction  that there  is a
sequence $v_n\in G$ such that $u_n\le q^{-p/q} v_n.$ This is
trivial for $n=0;$ now assume  it is proved for $n=k.$  Then
$\Cal   A u_k\le q^{-p} \, \Cal   A  v_k$   so  that
 $$u_{k+1} =  \Cal A u_{k}+p^{-1}q^{1-p} \, f  \in
((q-1)q^{-p}+q^{-p})G= q^{1-p}G=q^{-p/q}G,$$
i.e. $u_{n} \le q^{-p/q} v_{n}$ for all $n$. Hence $(u_n)$
is bounded in   $L^0$ and   we   can   apply    Proposition
2.5.

For (3) define $u_0=0$ and $u_n = \Cal A u_{n-1}+ f$.
Then for $n \ge 1$ by induction
$f \le u_n \le c_n \, f$ where $c_1=1$ and
$$c_{n+1} =  q^{-1} p^{1-q} \, c_n^q + 1.$$
Since $x_0=p$ is the only root of the equation $x=
q^{-1} p^{1-q} \, x^q + 1$ and $c_1=1$ it is easy to see
that $\lim_{n \to \infty} c_n = p$ and hence $v=\sup_{n}
u_n$ is a solution to (2.1) such that $f \le v \le p \, f$.

To prove (4)  we will need the inequality
$$\Cal   A (f + g) \le  [(\Cal A f)^{1/q} +
 (\Cal A g)^{1/q}]^q$$
for all $q \ge 1$ and $f, g \ge 0$. (For the reverse
inequality in case $0<q<1$
see [LT], p. 55). Note that if $T$ is a
kernel operator then this follows from Minkowski's
inequality. For an arbitrary $T$, using the
convexity of $\Cal   A$ we have
$$
\Cal   A (f + g) \le \lambda \, \Cal   A (\frac f \lambda)
+ \mu \, \Cal   A (\frac g \mu )
= \lambda^{1-q} \, \Cal   A f + \mu^{1-q} \, \Cal   A g
$$
for all  $\lambda, \mu > 0, \, \lambda + \mu =1$.
It then follows that we have:
$$ \Cal A(f+g) \le \inf\Sb\lambda+\mu=1 \\ \lambda, \, \mu>0\endSb
\{ \lambda^{1-q}\Cal
Af +\mu^{1-q}\Cal Af\}\quad \text{a.e.}$$

The desired inequality follows for a.e. $x$ by letting
$$\lambda =
\frac {(\Cal A f(x))^{1/q}} {(\Cal A f(x))^{1/q} +
 (\Cal A g(x))^{1/q}}, \quad \mu =
\frac {(\Cal A g(x))^{1/q}} {(\Cal A f(x))^{1/q} +
 (\Cal A g(x))^{1/q}},$$  when $\Cal Af(x),\Cal Ag(x)>0.$

Now for  $u_0= f$ and $u_n = \Cal A u_{n-1}+ f$
obviously $u_n \ge f + \Cal A f$ if $n \ge 1$.
To get the upper estimate we show
by induction that
$u_n \le f + c_n \, \Cal A f $
where
$$c_{n} = (1 +  q^{-1} p^{1-q} \, c_{n-1})^q.$$
This
is true for $n=1$ since $c_1=0$ and $c_2=1$. Assuming it is
true for $n=k$, we have
$$\aligned u_{k+1} = \Cal A u_{k}+ f & \le \Cal A (f + c_k \,
\Cal A f)^q + f \\
& \le [(\Cal A f)^{1/q} + c_k \, (\Cal A^2 f)^{1/q}]^q+ f
\\ & \le (1 + q^{-1} p^{1-q} \, c_k)^q \, \Cal A f + f.
\endaligned$$
Hence $ u_{k+1} \le f + c_{k+1} \, \Cal A f$ where
$c_{k+1}= (1 + q^{-1} p^{1-q} \, c_k)^q$. It remains to
note that $x_0 = p^q$ is obviously the only root of the
 equation $x= (1 + q^{-1} p^{1-q} x)^q$ and
$\lim_{n \to \infty} c_n = p^q$. Thus $v=\sup_{n}
u_n$ is a solution to (2.1) with the desired pointwise
 estimates.\qed\enddemo

\demo{Remark} It is not difficult to see that all the
constants in Proposition 2.7 are sharp.\enddemo

\proclaim{Proposition 2.8}If the set $S$ is bounded then  it
is also closed.\endproclaim

\demo{Proof}We need only check that if $0\le f_n\uparrow  f$
a.e. and  each $f_n\in  S$ then  $f\in S.$  Let $u_n$ be the
minimal solution corresponding to $f_n.$ Then  $ \Cal A
u_n \le u$
so that by Proposition  2.6, $p^{-1}q^{1-p}u_n\in S$.   Thus
$(u_n)$ is an increasing sequence which is bounded in $L^0.$
Let $u=\sup_n u_n.$ Then $u= \Cal A u  +f.$\qed\enddemo

\proclaim{Definition}  We  define  the  {\it solution space}
$\Cal Z_{q,T}$ to be the space of all $f\in L^0$ so that for
some  $\epsilon>0$ we have $\epsilon \, |f|\in  S.$ We
define
$$\|f\|_{\Cal Z}=\inf \, \{\alpha>0:\alpha^{-1}|f|\in S\}.$$
\endproclaim

We note that $S$ is nondegenerate if there exists $u>0$ a.e.
with $Ku^q\le u <\infty$ a.e.
Let us assume that $S$ is both bounded and
nondegenerate.  Then $\Cal Z$ is a Banach function space for
the norm $\|\,\|_Z.$

\proclaim{Proposition  2.9}Suppose  $f\in  L^0_+$.    Then
$f\in \Cal Z$ if and only if there exists $u\ge f$ and $C>0$
so that $u \in  L^0_+$ and  $ \Cal A u \le C u.$\endproclaim

\demo{Proof}If $f\in  \Cal Z$  then $f/\|f\|_{\Cal  Z}\in S$
and so  there exists  $v\ge f/\|f\|_{\Cal  Z}$ with
$\Cal A v \le
v.$  But  then  let  $u=\|f\|_{\Cal  Z}v.$  We have
$ \Cal A u \le
\|f\|_{\Cal Z}^{q-1}u.$ Conversely if $ \Cal A u \le Cu$
 then  let
$v=C^{-1/(q-1)}$ so  that $\Cal A v  \le  v.$ Then
by Proposition
2.6,  we  have  $p^{-1}q^{1-p}v\in  U$  so  that   $u\in\Cal
Z$.\qed\enddemo

\demo{Remark}If  we  define
$$|f|_{\Cal   Z}=\inf \, \{\alpha>0:
T|f|^q\le \alpha^{q-1} \, |f| \}$$
  then $|f|_{\Cal Z}$  is
an equivalent norm on $\Cal Z.$\enddemo

We next state  properties of $\Cal  Z$ which will  be useful
later.

\proclaim{Theorem 2.10} (1) For any  $f\in L^0_+$  we have
$f\in  \Cal  Z$  if  and  only  if  $\Cal A f \in\Cal  Z$ and
$$\|f\|_{\Cal Z}^q    \le \| \Cal A f \|_{\Cal    Z}   \le
pq^{p-1}\|f\|_{\Cal Z}^q.$$
(2) Suppose $\Cal X$ is a Banach
function space which contains $\Cal Z$.  Then for any  $f\in
\Cal   X_+$   we   have   $f\in\Cal   Z$   if  and  only  if
$\displaystyle{\sup_{n\ge 0}} \, \|\Cal  A^nf\|_{\Cal   X}
^{1/q^n}<\infty$  and   $$
\limsup_{n\to\infty} \, \|\Cal   A^nf\|_{\Cal   X}^{1/q^n}
 \le
\|f\|_{\Cal  Z}\le   pq^{p-1} \limsup_{n\to\infty} \, \|\Cal
A^n f\|_{\Cal X}^{1/q^n}.$$
(3) If $f\in L^0_+$, then $f \in Z$   if and only if
 $\displaystyle{\limsup_{n\to\infty} \,
 (\Cal  A^n f)^{1/q^n}} \in
L^{\infty}$ and   then   $$ \|\limsup_{n\to\infty} \, (\Cal
A^n f)^{1/q^n}\|_{\infty} \le \|f\|_{\Cal Z} \le
pq^{p-1} \|\limsup_{n\to\infty}(\Cal
A^n f)^{1/q^n}\|_{\infty}.$$
(4)  If $\Cal X$  is as
in   (2)   then
$$\limsup_{n\to\infty}\|\Cal   A^n f\|_{\Cal X}^{1/q^n} =
\|\limsup_{n\to\infty} \, (\Cal A^n f)^{1/q^n}\|_{\infty}.$$

\endproclaim

\demo{Proof}(1)  Assume  first  that  $f\in  \Cal Z$; we may
assume  $\|f\|_{\Cal  Z}=1.$  Then  there  exists  $u$  with
$u= \Cal A u +f.$ Then  $u^q \ge
(\Cal A u)^q +f^q$  and so  $\Cal A u  \ge
\Cal A^2 u   + \Cal A f$   so   that
 $\|\Cal A f  \|_{\Cal  Z}\le  1.$
Conversely if  $\| \Cal A f  \|_{\Cal Z}=1$  then we
conclude that
there exists $u$ with $u= \Cal A u + \Cal A f.$ But then
 $\Cal  A^nf\le
u$   for   $n\ge   1$,   so   that   by   Proposition   2.7,
$p^{-1}q^{1-p}f\in S$ or  $\|f\|_{\Cal Z}\le pq^{p-1}.$  The
result follows by homogeneity.

(2) Note that  if $f\in \Cal  Z$ then there  exists $u\in Z$
with  $u\ge  f/\|f\|_{\Cal  Z}$  with  $\Cal Au\le u.$ Hence
$\Cal  A^nf  \le  \|f\|_{\Cal  Z}^{q^n}u$.   It follows that
$$\limsup_{n\to\infty}\|  \Cal  A^nf\|_{\Cal  X}^{1/q^n} \le
\|f\|_{\Cal   Z}.$$   Conversely   suppose  $
\limsup_{n\to\infty} \|\Cal
A^n f\|_{\Cal X}^{1/q^n}=a<\infty.$ Then if $b>a$ we have
that $\{\Cal
A^n(b^{-1}f)\}_{n \ge 0}$ is bounded in $\Cal X.$
Applying Proposition 2.7(2)
this  yields  that  $p^{-1}q^{1-p}b^{-1}f\in  S$ i.e.  $f\in
\Cal Z$ and $\|f\|_{\Cal Z}\le bpq^{p-1}.$

(3) Assume $f \in \Cal Z.$  Then there exists $u \in
\Cal Z_+$
so that $ \Cal  A^n(f/\|f\|_{\Cal Z})\le u.$ Hence  $$ (\Cal
A^n f)^{1/q^n}  \le  \|f\|_{\Cal Z} u^{1/q^n}$$  so  that
 $$ \limsup_{n\to\infty}  \,
  \|(\Cal   A^n f)^{1/q^n}\|_{\infty}\le \|f\|_{\Cal Z}.$$
Conversely        if        $f\in        L^0_+$      and
$\|\limsup_{n\to\infty} \,
(\Cal A^n f)^{1/q^n}\|_{\infty}=a$
then   for   any   $b>a$   we   have  $\limsup_{n} \,
\Cal A^n (b^{-1}f)<\infty$ a.e.
 It follows from  Proposition 2.7
that $\|f\|_{\Cal Z}\le pq^{p-1}b$ and the result follows.

(4) The  convexity of  the map  $\Cal A$  is easily  seen to
imply that  both
$$\|\limsup_{n \to \infty}
 (\Cal A^n f)^{1/q^n}\|_{\infty} \quad
\text{and} \quad  \limsup_{n \to \infty}
\|\Cal A^n f\|_{\Cal X}^{1/q^n}$$
 define equivalent
norms on $\Cal Z$.  Both norms satisfy the identity  $\|\Cal
Af\|=\|f\|^q.$ This clearly makes them identical, since  one
obtains $\|\Cal  A^nf\|=\|f\|^{q^n}$ and  the two  norms are
equivalent.

\enddemo

\proclaim{Proposition 2.11}$\Cal Z'$ is an  order-continuous
Banach function  space so  that $\Cal  Z'$ is  separable and
$\Cal Z=\Cal  Z''$ can  be identified  as the  dual space of
$\Cal Z'.$\endproclaim

\demo{Proof}We show that $\Cal  Z$ is $q$-convex i.e.  there
exists a constant $C$ so that if $f_1,\ldots,f_n\in\Cal Z_+$
then      $\|(\sum_{k=1}^nf_k^q)^{1/q}\|_{\Cal Z}\le
C(\sum_{k=1}^n\|f_k\|_{\Cal Z}^q)^{1/q}.$ Indeed we note  as
in the preceding  Proposition that $\Cal A f_k \le
 \|f_k\|_{\Cal
Z}^{q-1}f_k.$ Let $g=(\sum_{k=1}^nf_k^q)^{1/q}.$ Then
$$\Cal A g
\le      \sum_{k=1}^n\|f_k\|_{\Cal      Z}^{q-1}f_k      \le
(\sum_{k=1}^n\|f_k\|_{\Cal   Z}^q)^{1/p}g$$
by  H\"older's
inequality.    Hence  again  arguing  as  in  the  preceding
Proposition if $h=g/(\sum_{k=1}^n\|f_k\|_{\Cal  Z}^q)^{1/q}$
then  $\Cal A h \le  h$  so  that  by  Proposition  2.5,  we
 have
$\|h\|_{\Cal  Z}\le  pq^{p-1}.$  Thus  we obtain the desired
inequality with constant $C=pq^{p-1}.$ (Note that we have in
effect  proved  that  the   norm  $f\to  |f|_{\Cal  Z}$   is
$q$-convex with constant one.)
Now (\cite{31}) this implies  that $\Cal Z'$ is  $p$-concave
and in particular order-continuous.\qed\enddemo

We  will  now  prove  a  result  which is except for certain
technical assumptions the same  as the main result  of Baras
and  Pierre  \cite{6},  but  our  approach  is   completely
different.

\proclaim{Theorem 2.12}If $g\in \Cal Z'$ then $$ \|g\|_{\Cal
Z'}=pq^{p-1} \, \inf
\left\{\int \frac{h^p}{(T^*h)^{p-1}}d\sigma:   \
h\in\Cal Z',\ h\ge |g|\right\}.$$ \endproclaim

\demo{Remark}It follows that if $f\ge 0$ a.e. then $f\in  S$
if  and  only  if  we  have $\int fh\,d\sigma\le 1$ whenever
$h\in  \Cal  Z'_+$  and  $\int (h^p/(T^*h)^{p-1})d\sigma \le
p^{-1}q^{1-p}.$ This is similar to the formulation in
 \cite{6}.  Notice that we do not assume that $T$
is given by  a kernel function  $K$, which is  necessary for
the  Baras-Pierre  approach.    Our  arguments  are   purely
functional analytic, and depend only on duality.  Of  course
we  are  restricting  our  attention  to  the case of
functions  of  the  type  $x\to  x^q$ while Baras and Pierre
consider more general convex functions.\enddemo

\demo{Proof}  Let  us  define  $V\subset\Cal  Z$ by $V=\{f:\
\exists u\in  \Cal Z_+:   u\ge Tu^q+f\}.$  Note that  $V\cap
\Cal Z_+=S,$ and that $V$ is convex.

We first show that $V$  is weak$^*$-closed.  To do  this, by
the  Banach-Dieudonn\'e  theorem  it  suffices  to show that
$V\cap  \alpha  B_{\Cal  Z}$  is  weak$^*$-closed  for   all
$\alpha>0.$  Since  $\Cal  Z'$  is  separable it suffices to
consider  a  sequence  $f_n\in  V\cap  \alpha  B_{\Cal   Z}$
converging to some $f\in \Cal Z$ for the  weak$^*$-topology.
Thus there exist $u_n\in \Cal Z_+$ with $u_n-Tu_n^q\ge f_n.$
Now since $\alpha^{-1}|f_n| \in S$ there exists  $v_n \in
\Cal Z_+$     with     $v_n - Tv_n^q=\alpha^{-1}|f_n|.$    Let
$(1+\alpha)w_n=u_n+\alpha    v_n.$    Then    by   convexity
$(1+\alpha)Tw_n^q\le     Tu_n^q+\alpha     Tv_n^q.$    Hence
$$(1+\alpha)(w_n-Tw_n^q)\ge  f_n  +|f_n|\ge  0.$$
It follows
that  $\|v_n\|_{\Cal  Z},\|w_n\|_{\Cal  Z}\le  pq^{p-1}$ and
hence $\|u_n\|_{\Cal  Z} \le  (1+2\alpha)p^{-1}q^{1-p}.$ The
sequence $(u_n)$ is thus also bounded in $\Cal Z.$

By Komlos's theorem we can now pass to a common subsequence,
which we  still denote  by $(u_n)$  and $(f_n)$  so that the
sequences         $z_n=\frac1n(u_1+\cdots+u_n)$          and
$g_n=\frac1n(f_1+\cdots+f_n)$  are  a.e.  convergent.     If
$\varphi\in \Cal  Z'$ and  $\varphi>0$ a.e.  then $(g_n)$ is
weakly convergent to $f$ in $L^1(\varphi)$ and hence it must
be also a.e. convergent to $f.$

If $z=\lim_n z_n$ a.e., then $Tz^q\le \liminf_n Tz_n^q$ so
 that
$z-Tz^q\ge  \limsup_n  z_n-Tz_n^q.$  However $z_n-Tz_n^q\ge
\frac 1 n \sum_{k=1}^n u_k-Tu_k^q \ge g_n.$ Hence
$z-Tz^q\ge f$
and  $f\in  V.$  This  completes  the  proof  that  $V$   is
weak$^*$-closed.

Now  let  $V^0=\{h\in  \Cal  Z':    \int  h f\,d\sigma\le 1\
\forall f\in V\}.$  We show that  $h\in V^0$ if  and only if
$h\in\Cal Z'_+$ and $$ \int  \frac{h^p}{(T^*h)^{p-1}}d\sigma
\le  pq^{p-1}.$$  First  note   that  it  is  obvious   that
$V^0\subset \Cal Z'_+.$

Suppose $h\in V^{0}.$ Select a sequence $u_n\in\Cal Z_+$  so
that $u_n\uparrow q^{1-p}(h/Th)^{p-1}$ a.e.  (Here $0/0$  is
interpreted as $0$.)   Then
$$ \int (u_n-Tu_n^q)h\,  d\sigma \le 1.$$
Noting that $u_n h \in L^1$ we can rewrite this as
$$\int (u_n h- u_n^q T^*h) \, d\sigma \le 1.$$
Pointwise    we    note    that   $u_n h-u_n^q T^*h
\uparrow
p^{-1}q^{1-p}h^p/(T^*h)^{p-1}.$    So by the Monotone
Convergence Theorem $$ \int  \frac{h^p}{(T^*h)^{p-1}}d\sigma
\le pq^{p-1}.$$

Conversely    suppose     $h\ge    0$     and    $$     \int
\frac{h^p}{(T^*h)^{p-1}}d\sigma  \le  pq^{p-1}.$$  Then   if
$f\le u-Tu^q$ where $u\in \Cal Z_+$ we have:
$$ \align \int
f h \, d\sigma  &\le  \int (uh  -  hTu^q) \,  d\sigma
\\ &= \int (uh-u^qT^*h) \, d\sigma \\
&\le \int  p^{-1}q^{1-p}
\frac{h^p}{(T^*h)^{p-1}} \, d\sigma\le 1. \endalign $$

Now by the theorem of bipolars (or the Hahn-Banach theorem)
we  have  $V^{00}=V.$  We  also have $(V^{0}-\Cal Z'_+)^{0}=
V\cap\Cal Z_+=S.$  From this  we obtain  by the  Hahn-Banach
theorem that $S^0$ is the closure in $\Cal Z'$ of the convex
set $V^{0}-\Cal Z'_+.$ In particular if $g\in S'=S^0\cap\Cal
Z'_+$ and $\epsilon>0$ then  there exists $h \in  V^{0}$ and
$g'\in \Cal Z'_+$ so that $\|g'\|_{\Cal Z'}\le \epsilon$ and
$h+g'\ge g.$

Fix $g=g_0\in S'$ and $\epsilon>0.$ Then by induction we can
construct  sequences  $(g_n)$  and  $(h_n)$  so that $h_n\in
\epsilon^{n-1}V^{0},$ $\|g_n\|_{\Cal Z'}\le \epsilon^n $ and
$g\le h_1+\cdots+h_n+g_n.$  Clearly $g\le
\sum_{n=1}^{\infty}h_n=h$ say.

Now
$$  \align\int  \frac{h^p}{(T^*h)^{p-1}}  d\sigma  &\le
(\sum_{n=1}^{\infty}(\int
\frac{h_n^p}{(T^*h)^{p-1}} \, d\sigma)^{1/p})^p
\\ &\le pq^{p-1}
(\sum_{n=1}^{\infty}\epsilon^{(n-1)/p})^p\\  &\le
pq^{p-1}(1-\epsilon^{1/p})^{-p}.\endalign $$

Since $\epsilon>0$ is arbitrary this  completes  the
proof.\qed\enddemo

\vskip10pt\heading{3. Quasi-metric kernels and infinitesimal
inequalities} \endheading \vskip10pt

We will now specialize the positive operator $T$  considered
in the previous section.  Let $K$ be a positive Borel kernel
function $K:X\times X\to  (0,\infty]$ (note that  $K(x,y)>0$
for all $x,y$ and that $K(x,y)=\infty$ is allowed).  We will
say that $K$ satisfies the {\it quasi-metric inequality}  if
$K$ is symmetric and there is a constant $\kappa\ge 1$  such
that for  all $x,y$  we have  $$ \frac{1}{K(x,y)}\le  \kappa
\left(\frac{1}{K(x,z)}+\frac{1}{K(z,y)}\right).\tag 3.1 $$

Under  these  conditions  it  is  natural  to  introduce the
quasi-metric $\rho(x,y)=(K(x,y))^{-1}.$ Note however that we
do not assume that  $K(x,x)=\infty$ and so $\rho(x,x)>0$  is
possible.   We can  then also  introduce the  ball of radius
$r>0$ i.e.  $$  B_r(x)=\{y:\rho(x,y)\le r\}$$ but note  that
this set can be empty.

One large class  of examples is  created by taking  $d$ as a
metric  on  $X$   and  $K(x,y)=d(x,y)^{-\alpha}$  for   some
$\alpha>0$; this  kernel defines  an operator  of fractional
integral type.  We will refer to a Borel set $B\subset X$ as
{\it bounded} if $\sup_{x,y\in B}\rho(x,y)<\infty.$

We suppose as in the previous section that there is given  a
$\sigma-$finite  Borel  measure  $\sigma$  on $X.$ Let $\Cal
M_+(X)$ be the space  of all positive $\sigma-$finite  Borel
measures on  $X$.   For each  $\omega\in\Cal M_+(X)$  we can
define       $K\omega\in       L^{0,\#}_+ (X,\sigma)$   by
$$K\omega(x)=\int_X    K(x,y)\,d\omega(y).$$    For    $f\in
L^0_+(X,\sigma)$     we     define     $$Kf(x)     =\int_X
K(x,y)f(y)\,d\sigma(y)          =K\omega(x)$$          where
$d\omega=f\,d\sigma.$  (Thus  we  identify  $L^0_+$  as  a
subset of $\Cal M_+(X).$) Sometimes we will write
$K f = K^\sigma f$ to emphasize that $K$ is defined on
$L^0_+(X,\sigma)$.

It is natural to require that  $K$
defines a positive operator as described in Section 2;  this
requires only the existence of a strictly positive  function
$w$ so that $Kw<\infty$ a.e.; however, this assumption  does
not  affect  the  results  of  the paper except to eliminate
triviality.

For  $1\le s \le  \infty$   we  say   $f\in  L^s_{loc}$  if
$f\chi_B\in L^s$  for every  bounded Borel  set $B.$  We say
that  $\omega\in\Cal  M_+(X)$  is  {\it  locally  finite} if
$|B|_{\omega}<\infty$ for every bounded Borel set $B.$

Our  first  proposition  gives  an  alternative  formula for
$K\omega.$

\proclaim{Proposition 3.1} Let $\omega\in \Cal M_+(X)$. Then
$$K\omega(x) = \int_0^\infty \frac { |B_r  (x)|_\omega}{r^2}
\,  {  dr},  \qquad  x  \in  X.  \tag  3.2  $$  \endproclaim
\demo{Proof}  For  a  fixed  $x  \in  X$,  we  can   rewrite
$K\omega(x)$,  using  the  distribution  function  of  $K(x,
\cdot)$, as $$K\omega (x) = \int_X K(x, y) \, d \omega (y) =
\int_0^\infty |\{y:  \, K(x, y) > t\}|_\omega \, d t.$$ Then
the substitution $r=1/t$ gives $$K\omega (x) = \int_0^\infty
\frac {  |B_r (x)|_\omega}  {r^2} \,  { dr},$$  which proves
Proposition 3.1. \qed\enddemo

\proclaim{Proposition 3.2}If $K$ satisfies the  quasi-metric
assumption then $K\omega<\infty$ a.e. implies that  $\omega$
is locally finite.\endproclaim

\demo{Proof}If    $B$     is    a     bounded    set     and
$|B|_{\omega}=\infty$ then we  can apply the  representation
3.2 to deduce  that $K\omega(x)=\infty$ everywhere.   Indeed
for any $x$  there is a  large enough $r$  so that $B\subset
B_r(x).$\qed\enddemo

One of our main tools is the following decomposition of  $K$
into its ``lower'' and ``upper'' parts.  For any $a>0$,  let
$$L_a (x, y) = \min \,  [K(x, y), a^{-1}], \quad (x, y)  \in
X\times X  . \tag  3.3$$ If  $K$ satisfies  the quasi-metric
assumption (3.1), then obviously, $L_a$ also satisfies (3.1)
with the same  constant $\kappa$.   We now split  the kernel
into the lower part $L_a$ and the upper part $U_a=K-L_a.$

\proclaim{Proposition 3.3}  Suppose $\omega\in\Cal  M_+(X).$
Then
$$L_a  \omega  (x)   =  \int_a^\infty  \frac  {   |B_r
(x)|_\omega} {r^2} \, { dr}, \qquad x \in X, \tag 3.4  $$
and
$$U_a  \omega (x)  = \int_0^a  \frac {|B_r  (x)|_\omega}
{r^2} \, { dr}, \qquad x \in X. \tag 3.5 $$
\endproclaim

\demo{Proof}  By  Proposition  3.1,  we  have
$$L_a  \omega(x) =
\int_0^\infty \frac { |\widetilde {B_r} (x)|_\omega} {r^2} {
dr},$$ where $\widetilde{B_r} (x)$ is a ball associated with
the kernel  $L_a$.   Obviously, $\widetilde{B_r}  (x) =  B_r
(x)$ if $r> a$, and $\widetilde{B_r} (x) = \emptyset$ if  $r
\le a$.   This  yields (3.4),  and hence  $$U_a\omega (x)  =
K\omega (x) - \int_a^\infty \frac { |B_r (x)|_\omega}  {r^2}
\, {  dr} =  \int_0^a \frac  { |B_r  (x)|_\omega} {r^2} \, {
dr}.$$ The proof of Proposition 3.2 is complete.
\qed\enddemo

As we will see in the next Proposition, the truncated kernel
$L_a$ or lower  part of the  operator has certain  stability
properties   for   kernels   satisfying   the    quasi-metric
inequality:  in particular $L_a\omega$ obeys a  Harnack-type
inequality on  a ball  $B_a(x)$ and  $(L_a\omega)^{-1}$ is a
quasi-concave function of $a>0.$

\proclaim{Proposition 3.4}  Suppose that  $K$ satisfies  the
quasi-metric condition (3.1) with a constant $\kappa \ge 1$.
Suppose $\omega\in\Cal M_+(X)$

(1) For all $a, b > 0$, $$L_a\omega(x) \le \max \, \left (1,
\frac {b} {a}\right) \, L_b\omega(x).  \tag 3.6 $$

(2)  For  all  balls  $B_a  (x)$,  $$ \dfrac {1} {2 \kappa}
\sup_{B_a (x)} L_a \omega\le L_a \omega (x) \le 2 \kappa  \,
\inf_{B_a (x)} L_a \omega . \tag 3.7$$ \endproclaim

\demo{Remark}It follows from (ii) that if $K$ is an operator
in the  sense defined  in Section  2 so  that for some $w>0$
a.e.  we  have  $Kw<\infty$  a.e.  then  we  must  have that
$\sigma$ is locally finite.\enddemo

\demo{Proof} If  $a \ge  b$, then  obviously $L_a\omega  \le
L_b\omega$.   Suppose $a<  b$, so  that $\delta  = b/a  >1$.
Then the substitution  $t = \delta  \, r$ gives  $$L_a\omega
(x) = \int_{a}^\infty \frac  { |B_r (x)|_\omega} {r^2}  \, {
dr}  =  \delta  \,  \int_{b}^\infty  \frac  {  |B_{t/\delta}
(x)|_\omega}   {t^2}   \,   {   dt}$$   $$\le   \delta    \,
\int_{b}^\infty \frac { |B_{t} (x)|_\omega} {t^2} \, { dt} =
\frac {b} {a} \, L_a\omega (x).$$ This proves statement (1).

To prove  statement (2),  notice that  if $y  \in B_a (x)$,
then by (3.1) $B_r(y) \subset B_{2 \kappa r} (x)$ for all $r
\ge a$.  Hence $$L_a \omega (y) = \int_a^\infty \frac { |B_r
(y)|_\omega} {r^2} \, { dr} \le \int_a^\infty \frac {  |B_{2
\kappa r}  (x)|_\omega} {r^2}  \, {  dr}$$ $$=  2 \kappa  \,
\int_{2 \kappa a}^\infty  \frac { |B_{r}  (x)|_\omega} {r^2}
\, {  dr} =  2 \kappa  \, K_{2  \kappa a}  \omega (x)  \le 2
\kappa \, L_a\omega (x).$$ This proves the lower estimate in
(3.7).  The  upper estimate is  proved in a  similar manner.
The proof of Proposition 3.3 is complete.  \qed\enddemo

Now  let  us  fix  $1<q<\infty$  and  consider the nonlinear
equation $$u=Ku^q +f\tag 3.8$$ for  $u\ge 0, f\ge 0.$ As  in
the previous section let $S=S_{q,K}$  be the set of $f$  for
which (3.8) has a solution $u\in L^0_+.$

\proclaim{Proposition 3.5}If $K$ satisfies the  quasi-metric
assumption  then  either  $S=\{0\}$  or  $S$  is bounded and
nondegenerate.  Furthermore, if  $S\neq \{0\}$ then for  any
bounded  Borel  set  $B$  there  exists $\epsilon>0$ so that
$\epsilon\chi_B\in S.$\endproclaim

\demo{Proof}If  $S\neq  \{0\}$  then  there exists $u\neq 0$
such that $u\ge Ku^q.$ But $Ku^q>0$ a.e. so that $u>0$  a.e.
Now by Proposition 2.7 we have $pq^{p-1}u\in S$ since  $\Cal
A^nu\le u$  for all  $n\in \Bbb  N.$ This  shows that $S$ is
nondegenerate.

To check boundedness fix $x\in X$ and suppose that for  some
$a>0$ we have $|B_a(x)|_{\sigma}>0.$ Suppose $0\neq f\in  S$
and so  $f\le u$  where $Ku^q\le  u.$ Then  $u>0$ a.e.,  and
$u\in     L^q_{loc}.$     Further     $$    L_au^q(x)    \ge
a^{-1}\int_{B_a(x)}u(x)^q \, d\sigma.$$ By Proposition 3.4
this
yields    that    $$u\ge    L_au^q    \ge   \frac{1}{2\kappa
a}\left(\int_{B_a(x)} u^q \,
d\sigma\right)\chi_{B_a(x)}.$$    We
deduce  immediately  that  $|B_{a}(x)|_{\sigma}<\infty$  (so
that $\sigma$ must be  locally finite if $S\neq  \{0\}).$ We
also  note  that   $\epsilon \, \chi_{B_a(x)}\le  u$  for
 some
$\epsilon>0$  which  justifies  the  last  statement  in the
Proposition.

We continue the proof that $S$ is bounded;  we
have    $$    \int_{B_a(x)} u^q \, d\sigma \le 2\kappa a
|B_a(x)|_{\sigma}^{-1}   \int_{B_a(x)}u \, d\sigma.$$
 Applying
H\"older's   inequality    $$\int_{B_a(x)}u \, d\sigma \le
|B_a(x)|_{\sigma}^{1/p}\left(\int_{B_a(x)} u^q \,
d \sigma\right)^{1/q}.$$
Combining we  get an  estimate:   $$ \int_{B_a(x)}u^qd\sigma
\le  (2\kappa  a)^p|B_a(x)|_{\sigma}^{-p/q}.$$  Since   this
holds for  any such  ball, it  follows that  the set  $S$ is
bounded in $L^0.$\qed\enddemo

Thus  $\Cal Z_{q,K}$  either  reduces  to  $\{0\}$  or is a
nondegenerate Banach function space with the norm induced by
$S_{q,K}.$ We  also have  $\chi_B\in\Cal Z$  for any  metric
ball   $B=B_a(x).$   We    now   prove   an    estimate   on
$\|\chi_B\|_{\Cal Z}$; this is somewhat more complicated and
requires some preliminary work.

\proclaim{Lemma  3.6}Suppose   $x\in  X$   and  $a>0.$   Let
$B=B_a(x)$ and then $B_j=B_{a2^{-j}}(x)$ for $j\in \Bbb  N.$
Let $c_j=(4\kappa  a)^{-1}2^j|B_j|_{\sigma}.$ Then  $K\chi_B
\ge         \phi_B$         where         $\phi_B=\sum_{j\ge
0}c_j\chi_{B_j}.$\endproclaim

\demo{Proof}Let   $b=4\kappa   a.$   We   have   $$   \align
K\chi_B(y)&=\int_0^{\infty}\frac{|B                     \cap
B_r(y)|_\sigma}{r^2}dr\\                                &\ge
\sum_{j=0}^{\infty}\int_{b2^{-(j+1)}}^{b2^{-j}}\frac     {|B
\cap     B_r(y)|_{     \sigma}}{r^2}dr\\     &\ge    \frac1b
\sum_{j=0}^{\infty}2^j|B \cap B_{b2^{-(j+1)}}(y)|_{\sigma}\\
&\ge             \frac1b              \sum_{j=0}^{\infty}2^j
|B_j|_{\sigma}\chi_{B_j}(y) \endalign $$

since    if    $y\in    B_j$    then    $B_j\subset    B\cap
B_{b2^{-(j+1)}}.$\qed\enddemo

\proclaim{Lemma 3.7}Suppose $0<s<\infty$.  Then $$ K\phi_B^s
\ge \frac1{(s+1)}\phi_B^{s+1}.$$ \endproclaim

\demo{Proof}    Let     $\phi=\phi_B.$    We     also    let
$\beta_j=\sum_{i=0}^jc_i$ (and $\beta_{-1}=0$).  We will use
the   inequality   that  $\beta_j^{s+1}-\beta_{j-1}^{s+1}\le
(s+1)c_j\beta_j^s.$

We     start     with     the     observation     that    $$
\phi^s=\sum_{j=0}^{\infty}\alpha_j\chi_{B_j}$$         where
$\alpha_j=\beta_j^s-\beta_{j-1}^s$.
 Note that $K\chi_{B_j}\ge
\sum_{i=j}^{\infty}c_i\chi_{B_i}$ by Lemma 3.6.  It  follows
that         $$         \align         K\phi^s          &\ge
\sum_{j=0}^{\infty}\alpha_j\sum_{i=j}^{\infty}c_i\chi_{B_i}
\\     &=     \sum_{i=0}^{\infty}c_i\sum_{j=0}^i    \alpha_j
\chi_{B_i}\\        &=        \sum_{i=0}^{\infty}        c_i
\beta_i^s\chi_{B_i}\\  &\ge   \frac1{s+1}\sum_{i=0}^{\infty}
(\beta_i^{s+1}-\beta_{i-1}^{s+1})\chi_{B_i}\\             &=
\frac1{s+1}\phi^{s+1}.    \endalign  $$  This  completes the
proof.\qed\enddemo

Now  we  introduce  the   quantity  $$  M(x,a)  =   \int_0^a
\frac{|B_r(x)|_{\sigma}}{r^2}dr=  U_a(\bold  1)(x).$$   Here
$\bold 1=\chi_X.$ Thus $U_a(\bold 1)=U_a\sigma.$

\proclaim{Theorem  3.8}If  $\Cal  Z\neq\{0\}$ and $B=B_a(x)$
then
$$C(q) \, \left(\frac{M(x,a)}{4\kappa}\right)^{p/q}\le
\|\chi_B\|_{\Cal Z}.$$ \endproclaim

\demo{Proof}This is trivial if $|B|_{\sigma}=0.$  Otherwise,
define $\phi = \phi_B$ as in Lemmas 3.6 and 3.7.
 Let $\Cal  Af=Kf^q$
as  in  Section  2.  Then  $\Cal  A\chi_B\ge \phi$ and $\Cal
A\phi^s\ge  (s q +1)^{-1}\phi^{s q +1}.$  It
follows by induction
that
$$\align \Cal A^n \chi_B &\ge \,
\prod_{j=1}^{n-1} (1+q+q^2+\cdots+q^{j})^{-q^{n-j-1}} \, \,
\phi^{1+q+\cdots+q^{n-1}} \\ &= C(n,q) \, \phi^{(q^n-1)/(q-1)}.
\endalign$$
Then
$$\limsup_{n\to\infty} \, (\Cal A^n  \chi_B )^{1/q^n} \ge
\lim_{n\to\infty} C(n,q)^{1/q^n} \, \phi^{-1/(q-1)}.$$
Clearly,
$$\align \lim_{n\to\infty} C (n,q)^{1/q^n}  &=
\prod_{j=1}^{\infty} (1+q+q^2+\cdots+q^{j})^{-q^{-j-1}}\\
& = \prod_{j=1}^{\infty}
q^{-j q^{-j-1}}
\,  \prod_{j=1}^{\infty} (1+q^{-1}+q^{-2} +\cdots+
q^{-j})^{-q^{-j-1}}\\&\ge \prod_{j=1}^{\infty}
q^{-j q^{-j-1}} \,
\prod_{j=1}^{\infty} (1-q^{-1})^{q^{-j-1}}\\& =
q^{-(q-1)^2} \, (1-q^{-1})^{1/q(q-1)}.\endalign  $$
Thus
$$ \limsup_{n\to\infty} \, (\Cal A^n  \chi_B )^{1/q^n}
\ge  q^{-(q-1)^2} \, (1-q^{-1})^{1/q(q-1)}    \,
 \phi^{1/(q-1)}.$$
Appealing to Theorem 2.10 we obtain
$$\|\chi_B\|_{\Cal   Z}\ge   q^{-(q-1)^2} \,
(1-q^{-1})^{1/q(q-1)}
 \,  \|\phi\|_{\infty}^{p/q}.$$
It remains to notice that by definition of $\phi$
$$\align  \|\phi\|_{\infty}&=\sum_{j=0}^{\infty}c_j\\
&=\frac{1}{4\kappa  a}\sum_{j\ge  0}2^j|B_j|_{\sigma}\\
& \ge
\frac1{4\kappa}\int_0^a\frac{|B_r(x)|_{\sigma}}{r^2}dr.
\endalign $$
This completes the proof.\qed\enddemo

\proclaim{Corollary  3.9}A  necessary  condition  for  $\Cal
Z\neq\{0\}$    is    that    for    any    $a>0$   we   have
$M^*(x,a)=\sup_{y\in                        B_a(x)}M(y,a)\in
L^{\infty}_{loc}.$\endproclaim

\demo{Proof}Note that if  $y\in B_a(x)$ then  $B_a(y)\subset
B_{2\kappa    a}(x)$     so    that
$$M(y,a)^{p/q}    \le C(q) \,
(4\kappa)^{p/q}\|\chi_{B_{2\kappa  a}(x)}\|_{\Cal
Z}.$$
\qed\enddemo

The following Theorem introduces  an inequality we name  the
{\it  infinitesimal   inequality};  this   is  a   necessary
condition on a measure $\omega$ so that $K\omega\in\Cal Z.$

\proclaim{Theorem   3.10}(The   infinitesimal   inequality.)
Assume that $K$ satisfies the quasi-metric condition.   Then
there is a constant $C=C(\kappa)$ so that if  $\omega\in\Cal
M_+(X)$ satisfies  $K\omega\in\Cal Z$  then for  every $x\in
X,$                                                       $$
\sup_{a>0}\left\{\int_0^a\frac{|B_t(x)|_{\sigma}}
{t^2}dt\right\}^{p/q}
\left\{\int_a^{\infty}\frac{|B_t(x)|_{\omega}}{t^2}dt\right\}
\le C\|K\omega\|_{\Cal Z}.  \tag 3.9 $$ \endproclaim

\demo{Remarks} (a) The  conclusion  is   that  (3.9)  holds
  {\it
everywhere}:  however,  when working with  the infinitesimal
inequality later it  will only be  necessary to assume  that
(3.9) holds $\sigma-$a.e. for each $a>0.$ \newline
(b) There is an appealing alternative form
 of  this
inequality (in  the almost  everywhere form),  namely
$$\sup_{a>0} \|(U_a\sigma)^{1/q}(L_a\omega)^{1/p}\|_{\infty}
 \le C'\|K\omega\|_{\Cal Z}^{1/p}.\tag 3.10 $$ \enddemo

\demo{Proof}We  note  that  for   all  $a>0, \,
x, y \in  X$  if
  $B=B_a(x),$ then
$$ \chi_B (y) \, L_a\omega (x) \le 2\kappa \, L_a\omega (y)$$
by Proposition 3.4.  Hence
$$ L_a\omega(x) \,  \|\chi_B\|_{\Cal Z} \le  2\kappa \,
\|K\omega\|_{\Cal  Z},$$
and  by Theorem 3.8
 $$L_a\omega(x) \, M(x,a)^{p/q} \le  C\|K\omega\|_{\Cal Z}$$
 where $C=C(\kappa).$ This yields (3.9).\qed\enddemo

\proclaim{Corollary 3.11}There  is a  constant $C=C(\kappa)$
so that  if $f\in  \Cal Z_+,$  then for  every $x\in  X,$ if
$d\omega=f^qd\sigma$
$$\sup_{a>0} \,
\left\{\int_0^a\frac{|B_t(x)|_{\sigma}}{t^2}dt\right\}^{1/q}
\left\{\int_a^{\infty}\frac{|B_t(x)|_{\omega}}
{t^2}dt\right\}^{1/p}
\le C\|f\|^{q-1}_{\Cal Z}.  \tag 3.11 $$ \endproclaim

\demo{Proof}Here $K\omega=Kf^q\in\Cal Z$ so that by  Theorem
3.10
$$\sup_{a>0}\left\{\int_0^a\frac{|B_t(x)|_{\sigma}}
{t^2}dt\right\}^{1/q}
\left\{\int_a^{\infty}\frac{|B_t(x)|_{\omega}}
{t^2}dt\right\}^{1/p}
\le  C \|Kf^q\|^{1/p}_{\Cal  Z}.    $$
However   $\|Kf^q\|\le
pq^{p-1}\|f\|_{\Cal Z}^q$ (Theorem  2.10) and so  the result
follows.\qed\enddemo

\vskip10pt\heading{4. Weighted  norm   inequalities  and
nonlinear integral equations}\endheading \vskip10pt

In  this   section  we  prove some of the main
results of the paper (Theorems 4.8 and 4.9).
We develop  connections   between  the
solvability of (2.1)  and weighted norm  inequalities, which
clarifies the  role  of  the  infinitesimal
inequalities  and  testing  conditions  of  Sawyer type.
Simple criteria of Koosis type (see \cite{18}) are given for
$\Cal W_{p,K}$ and $\Cal Z_{q,K}$ to be nondegenerate
Banach function spaces. Sharp estimates of
 $||\chi_B||_{\Cal Z_{q,K}}$
for a ball $B=B_r(x)$ are obtained.
We
retain the  assumption that  $K$ satisfies  the quasi-metric
condition (3.1).

We recall from Section 2 that $f\in L^0$ is an  $L^p$-weight
for $K$ if there is a constant $\gamma>0$ such that
$$  \int
|f| \, |Kg|^p \, d\sigma \le C \int |g|^p \,
d\sigma \tag 4.1$$ for all
$g\in L^p(\sigma).$  The space  of $L^p-$weights  is denoted
$\Cal  W_{p,K}$.    Since  $K(x,y)>0$  for  all $x,y$ if $K$
defines a positive operator  it is strictly positive  and we
have that  $\Cal W_{p,K}$  under the  norm $f\to \|f\|_{\Cal
W}=\inf \, \{C:\  (4.1)  \text{  holds}\}$  is  a
nondegenerate
Banach   function   space   on   $(X, \sigma)$  provided  the
assumption  that  $L^p\subset  \Cal  D(K)$  holds  (and this
condition is  clearly also  necessary).   Even without  this
assumption it is clear that $\{f:\|f\|_{\Cal W}\le 1\}$ is a
bounded subset of  $L^0.$ We recall  also that $\Cal  Y=\Cal
Y_{p,K}=\{f:|f|^q\in\Cal W_{p,K}\}$; if $\Cal W$ is a Banach
function space then so is $\Cal Y$ with the associated  norm
$f\to\||f|^q\|_{\Cal W}^{1/q}.$

We first prove a simple criterion for $\Cal W_{p,K}$ to be a
Banach  function  space  on  $(X,\sigma).$  To  do  this  we
introduce the quantity
$$ N(x,a) = \int L_a(x,y)^q \, d\sigma(y)
=\int_a^{\infty}\frac{|B_t(x)|_{\sigma}}{t^{q+1}}dt.$$
It is
immediate  that  if  $N(x,a)<\infty$  for  some  $x,a$  then
$\sigma$  is  locally  finite  and  that $N(y,b)<\infty$ for
every  $y\in  X,$  and  $b>0.$  This  follows  from applying
Proposition 3.4 to the modified kernel $K(x,y)^q$.

\proclaim{Proposition 4.1}(1) If $L^p\subset \Cal D(K)$ then
$\sigma$   is   locally  finite   and   for   every   $a>0,$
$N(x,a)<\infty.$ \newline (2) If  for every $a>0$ and  $x\in
X$  we  have  $N(x,a)<\infty$  and  $M^*(x,a)<\infty$   then
$L^p\subset  \Cal  D(K)$  (and  hence there exists $w\in\Cal
W_{p,K}$ with $w>0$ a.e.).  Furthermore we have the estimate
$$\|\chi_{B_a(x)}\|_{\Cal W}  \le
C(|B_a(x)|_{\sigma}N(x,a)^{p/q}+M^*(x,2\kappa a)^p)$$ where
$C$ depends only on $\kappa.$ \endproclaim

\demo{Proof}(1)    Pick     a    ball     $B=B_a(x)$    with
$|B|_{\sigma}>0$.  Then we  can apply Proposition 3.4  again
to get  that
$$  Kf\ge (2\kappa)^{-1}L_af(x) \,
\chi_B.   $$
If $L^p\subset \Cal D(K)$ this implies that $L_af(x)<\infty$
for   all   $f\in   L^p$   so   that  $N(x,a)<\infty.$  This
contradiction implies $N(x,a)<\infty$  for all $x\in  X$ and
$a>0.$

(2) It will suffice to prove that for every ball  $B=B_a(x)$
we have $\chi_B\in\Cal W.$ Note first that if $y \in B$
then by Proposition 3.4 and H\"older's inequality
$$  L_a f (y) \le
2\kappa \, L_a f(x) \, \le    2\kappa \,
N(x,a)^{1/q} \, \|f\|_p.$$
 Hence  $$  \int_B (L_a f)^p \, d \sigma \le
(2\kappa)^p \, N(x,a)^{p/q} \, |B|_{\sigma} \,
\int f^p \, d\sigma.$$

Now if $y\in B$  then $U_a f(y) \le M(y,a)\|f\|_{\infty}
 \le
M^*(x,a)\|f\|_{\infty}.$ Hence
$$  \|\chi_B U_a f\|_{\infty} \le
M^*(x,a) \, \|f\|_{\infty}.$$

Also
$$\int_B U_a f(y) \, d\sigma(y)\le  \int_{B_{2\kappa a}(x)}
\left(\int_{B_{  a}(y)} K(y,z) \, d\sigma(y)\right) f(z)
\, d\sigma(y)
\le M^*(x,2\kappa a) \, \|f\|_1.$$

By the  Riesz interpolation  theorem
$$  \int_B (U_a f)^p \, d \sigma \le
M^*(x,2\kappa  a)^p \, \int  f^p \, d\sigma.$$
The  Proposition now
follows by combining the estimates for $L_a f$ and $U_a f$.
  \qed\enddemo

Proposition 4.1(2) gives  a situation in  which $\Cal W$  is
nondegenerate  which  is   sufficient  to  cover   our  main
interests in this paper.  However it is possible to prove
more  general  results  which  can  be  regarded as
extending
previous Koosis-type  theorems due  to
Rubio de Francia \cite{42} and  Sawyer  \cite{43}
for the Riesz potential.

\proclaim{Proposition   4.2}Suppose   $K$   satisfies    the
quasi-metric  condition  and  that  $\mu$  is  a  nontrivial
$\sigma-$finite Borel  measure on  $X$ such  that for  every
$x\in  X$  and  $a>0$  we  have $$ M_{\mu}^*(x,a)=\sup_{y\in
B_a(x)}\int_0^a\frac{|B_t(y)|_{\mu}}{t^2}dt<\infty.$$    Let
$\sigma$ be a $\mu$-continuous $\sigma-$finite Borel measure
on  $X$.    Then  for  $(X,\sigma)$  we  have $\Cal W_{p,K}$
nondegenerate     if     and     only     if    $$    N(x,a)
=\int_a^{\infty}\frac{|B_t(x)|_{\sigma}}{t^2}dt    <\infty$$
for some (and hence all) $x\in X,\ a>0.$\endproclaim

\demo{Proof}The  necessity  of  the  condition  is proved in
Proposition 4.1.  For the converse, we remark first that  we
can replace $\mu$ by its $\sigma-$continuous part and  hence
assume $d\mu=w\,d\sigma$ for  some strictly positive  weight
function $w.$ We must  show that if $0\le f \in
\L^p (\sigma)$
then $Kf<\infty$ a.e.  Indeed in this case let $B=B_a(x)$ be
any ball.  Then  $$ \sup_{y\in B}L_af(y)\le 2\kappa  L_af(x)
\le  2\kappa  N(x,a)^{1/q}  \|f\|_p.$$
 Also
$$\align \int_B
w \, U_af\,d\sigma  &= \int_{y\in  B}\int_{B_{2\kappa a}(x)}
U_a(x,y) \, w(y) \, d \sigma(y)\, d\sigma(x)\\
&=   \int_{B_{2\kappa
a}(x)}  M(y,a)  f(y)\,d\sigma(y)\\  &\le  M^*(x,2\kappa a)
|B_{2\kappa   a}|_{\sigma}^{1/q}\|f\|_p.\endalign   $$  This
shows that $\chi_BKf<\infty$ $\sigma-$a.e.\qed\enddemo

The  following  result  was
proved for the special case when $K$ is the Riesz  potential
in \cite{42} and  \cite{43}.

\proclaim{Proposition 4.3} Suppose   $K$   satisfies    the
quasi-metric   condition    and   that    $\sigma$   is    a
$\sigma-$finite Borel  measure on  $X$ such  that for  every
$x\in  X$  and  $a>0$  we  have $M^*(x,a)<\infty.$ Let $v\in
L^0_+$.   In order  that there  exist $w\in  L^0_+$ with
$w>0$ whenever $v>0$ and
$$ \int (Kf)^pw\,d\sigma \le  C\int
f^pv\,d\sigma\tag 4.2$$
for all $f\in L^p(v\,d\sigma)$ it is
necessary      and      sufficient      that
$$     \int
L_a(x,y)^q \, v(y)^{1-q} \, d\sigma(y)<\infty$$
for some (and hence
all) $x\in X$ and $a>0.$\endproclaim

\demo{Proof}By  replacing  $X$  by  $X_0=\{x: v(x)>0\}$ and
restricting $\sigma$  we can  suppose $v>0$  everywhere. let
$d\nu=v^{1-q}d\sigma$ and let $K^{\nu}f=K(v^{-q/p}f).$  Then
if we put $f=\varphi v^{-q/p}$ then (4.2) becomes equivalent
to
$$ \int (K^{\nu} \varphi)^p \, (w/v) \, d\nu
\le C \int f^p \, d\nu$$
and so the result reduces to Proposition 4.2.  The necessary
and   sufficient   condition   required   is  that
$$\int
L_a(x,y)^qd\nu   <\infty$$
for    some   $x\in X$  and
$a>0.$\qed\enddemo

We  may  also  study weighted  norm  inequalities   for
measures
$\omega\in\Cal M_+(X).$ We  define $\widetilde
{\Cal W}_{p,K}$
to be the cone of $\omega$ such that for some constant
$C>0$ we have $$ \int |Kg|^p \, d\omega \le C \int  |g|^p
\, d\sigma \tag 4.3$$  for all  $g\in L^p.$  Of course
 if $\omega$  is
$\sigma-$continuous  then  we  can write $d\omega=f \,
d\sigma$
and $f$ satisfies (4.1).

The following characterization
of (4.3) is due to Sawyer and
Wheeden \cite {43},  \cite {44} (see also substantial improvements in
\cite {45} and \cite {54}) under  the  hypothesis that $X$ is a
homogeneous space in the sense of Coifman-Weiss \cite{11}
equipped
with a quasi-metric $d (x,y)$
and doubling measure
$\mu$ such that $|B_{2 r} (x)|_\mu \le C \, |B_{r} (x)|_\mu$
(where $B_r (x) = \{y: \,
d (x, y) < r\}$): Suppose $K(x, y)$ is a kernel
satisfying the inequalities
$$
K(x, y) \le C_1 \, K(x', y) \quad \hbox {if} \quad d (x', \, y)
 \le C_2 \, d (x, \, y), \tag {K1}$$
and
$$K(x, y) \le C_1 \, K(x, y') \quad \hbox {if} \quad d (x, \, y')
 \le C_2 \, d (x, \, y) \tag {K2}$$
for some $C_1 > 1$ and  $C_2 > 1$. Suppose that all annuli
 with respect to $d$ are nonempty. Then
 (4.3) holds  if and  only if  both of  the
following  {\it  testing   conditions}  hold,
$$\int_X   (K\sigma_B)^p \,  d \omega  \le C  \, |B|_\sigma
\tag 4.4$$
and
$$\int_X (K^* \omega_B)^q  \, d  \sigma \le  C \, |B|_\omega,
\tag  4.5$$  for  all  balls  $B$ associated with
 $d$;  here  $d\sigma_B=\chi_B\, d\sigma$, $K^*$ is  a formal
adjoint,
and as usual $1/p+1/q=1$, $1<p<\infty$. Under certain mild assumptions
it can be shown that
(4.4) and (4.5) are equivalent to similar inequalities with
integration
over $B$ in place of $X$ on the left-hand sides.
 (See  \cite{46} and \cite{54}.)

Observe that if $K(x,y)= 1/d(x, y)$ then
$(K1)$ and $(K2)$ hold automatically with $C_1$ = $C_2 = 1$ and
thus the testing inequalities (4.4), (4.5) characterize
$\widetilde {\Cal W}_{p,K}$.
As is shown in \cite{53}, under some
 assumptions on $X$ every kernel $K$
 which satisfies $(K1)$ and $(K2)$
 is controlled
by its symmetric ``dyadic'' analogue $K^d$ which satisfies
the quasi-metric inequality and is pointwise
smaller
than $K$; moreover, weighted norm inequalities
for integral operators generated by $K$ and $K^d$ hold
simultaneously so in the setting of homogeneous spaces
our quasi-metric assumptions and $(K1) - (K2)$ are in a
sense equivalent.

\demo{Remarks} (a) We emphasize that these  deep two weight
results  with  difficult and rather
technical proofs are not used in the present
paper. For our purposes it suffices to use simpler
pointwise
 characterizations of weighted norm
inequalities discussed in this section.
However, the  second  testing  condition (4.5), with $B$
in place of $X$ on the left-hand side, plays
an important role in our approach to the
 solvability problem for the equation  $u
= Ku^q + f$.
\newline
(b) It is easy to give a nonsymmetric version
of our  solvability results
 in the framework outlined above; we do not consider
it here to avoid unnecessary complications. For the
applications we have in mind it is more convenient to
restrict ourselves to  the case
of possibly nonsymmetric kernels $K$ such that $K(x, y)
\asymp s_1(x) \, K_0(x, y) \, s_2 (x)$ where $K_0(x, y)$
is a quasi-metric kernel and $s_1, s_2$ are arbitrary
positive weight functions. These generalizations
are obtained in Sec. 7 together with  applications to
nonlinear Dirichlet problems.\enddemo

We start with the following proposition which  is basic to
our results
and explains the significance of the quasi-metric assumption
for weighted norm inequalities.

\proclaim{Proposition  4.4}  Let  $K$  be  a  kernel   which
satisfies the  quasi-metric assumption.   Then  for any $g\in
L_{0_+},$  and  $1\le  s<\infty   ,$  $$(Kg)^s  \le  C   K(g
(Kg)^{s-1}),\tag  4.6  $$  where  $C  =  s(2 \kappa)^{s-1}$.
\endproclaim

\demo{Remarks} (a) Proposition  4.4 shows that
integral operators  with quasi-metric kernels
resemble Hardy's operator
$K  g (x) = \int_0^x g (t) \, dt$ which obviously has
the property $(Kg)^s  = s  \,  K(g
(Kg)^{s-1}).$\newline
(b) A different proof valid for
kernels satisfying assumptions (K1) and (K2)
is given in \cite{53}. \enddemo

\demo{Proof}Let $d  \nu =  g \,  d \sigma$.   By Proposition
3.1,  $$Kg  (x)  =  K\nu  (x)  =  \int_0^\infty \frac { |B_r
(x)|_\nu} {r^2} \, {  dr}.$$
Clearly, we
have $$  \align (K  g (x))^s  &= \left  (\int_0^\infty \frac
{|B_r (x)|_\nu}  {r^2} \,  { d  r} \right  )^s \\  & =  s \,
\int_0^\infty  \left  (\int_r^\infty  \frac  {|B_t (x)|_\nu}
{t^2} \,  { d  t} \right  )^{s-1} \,  \frac {|B_r  (x)|_\nu}
{r^2} \, { d r} \\ &= s \, \int_0^\infty \, {\int_{B_r  (x)}
\left (\int_r^\infty \frac {|B_t  (x)|_\nu} {t^2} \, {  d t}
\right  )^{s-1}  \,  d  \nu  (y)}  \,  \frac  { d r} {r^2} .
\endalign $$ We estimate the inside integral by  Proposition
3.3.   For $y  \in B_r  (x)$, we  have $$\int_r^\infty \frac
{|B_t (x)|_\nu} {t^2} \, { d t} = L^r g (x) \le 2 \kappa  \,
L^r g (y) \le  2 \kappa \, K  g (y),$$ where $L^r  g$ is the
``lower part'' of $Kg$.  Then $$\align (Kg (x))^s &\le s  \,
(2 \kappa)^{s-1} \, \int_0^\infty \, \frac { \int_{B_r  (x)}
(K^\sigma g)^{s-1} \, d \nu (y) } {r^2} \, { d r} \\ &= s \,
(2   \kappa)^{s-1}   \,   K   (g   \,   (K^\sigma  g)^{s-1})
(x),\endalign $$  which completes  the proof  of Proposition
4.4.\qed\enddemo

\proclaim{Proposition   4.5}Suppose   $K$   satisfies    the
quasi-metric  condition  and  that  $\omega\in\Cal  M_+(X).$
Suppose $K\omega\in Y_{p,K}$ or equivalently $(K\omega)^q\in
\Cal   W_{p,K}$   is   an   $L^p-$weight   for   $K.$   Then
$\omega\in\widetilde{\Cal W}_{p,K}$      (i.e. (4.3)
holds).\endproclaim

\demo{Proof} Suppose $0\le  g\in L^p(\sigma).$ Then,  if $C$
is   the   constant   of   Proposition  4.4,  $$  \align\int
(Kg)^pd\omega  &\le  C\int  K(g(Kg)^{p-1})d\omega\\ &= C\int
g(Kg)^{p-1}K\omega    d\sigma    \\    &\le   C   \left(\int
g^pd\sigma\right)^{1/p}\left(\int
(Kg)^p(K\omega)^qd\sigma\right)^{1/q}\\                 &\le
C\|(K\omega)^q\|_{\Cal      W}^{1/q}\|g\|_p^p      \\     &=
C\|K\omega\|_{\Cal Y}\|g\|_p^p.\endalign  $$ This  completes
the proof. \qed\enddemo

\proclaim{Theorem 4.6}Suppose  $f\in L^0_+$  and $f\in
\Cal Z_{q,K}$.    Then  $f\in\Cal  Y_{p,K}$  (so that $f^q$
is an
$L^p-$weight   for   $K$)   and   there   is   a    constant
$C=C(q,\kappa)$ such that  $\|f\|_{\Cal Y}\le C  \|f\|_{\Cal
Z},$ or equivalently  $\|f^q\|_{\Cal W}
\le C  \|f\|^q_{\Cal Z}.$  \endproclaim

\demo{Proof}Suppose first  that $Kf^q\le  f.$ If  $f$ is not
zero  then  $f>0$  a.e.;  furthermore under this hypothesis,
$\sigma$ is locally  finite i.e.   $|B|_{\sigma}<\infty$ for
all  bounded  sets  $B.$  Let  $0\le  g\in  L^{\infty}$ with
$\|g\|_p=1$ be  supported in  a bounded  set $B.$  Then note
that on $B$  $$ Kf(x) \ge  \alpha \int_B f\,d\sigma$$  where
$\alpha=\inf \,
\{K(x,y):\ x,y\in B\}.$ Hence $f\ge c>0$ on  $B$
and so $\chi_Bf\le c^{1-q}f^q$ and $K(\chi_Bf)\le c^{1-q}f.$
Now  for  any  $1\le  s<\infty$  if  $C_s=s(2\kappa)^{s-1},$
$$\align     \int     (Kg)^sf^q\,d\sigma     &\le    C_s\int
K(g(Kg)^{s-1})f^q\,d\sigma        \\        &=       C_s\int
\chi_Bg(Kg)^{s-1}Kf^q\,d\sigma\\ &\le  C_s\|g\|_{\infty}\int
(Kg)^{s-1}\chi_B\,f\,d\sigma.\\                         &\le
C_sc^{1-q}\|g\|_{\infty}     \int     (Kg)^{s-1}f^q\,d\sigma
\endalign $$

Since $|B|_{\sigma}<\infty$ this  implies by induction  that
$Kg\in L_s$ for every integer $s$ and hence all $s.$ We  now
take $s=p.$  $$\align \int  (Kg)^pf^q\,d\sigma &\le  C_p\int
g(Kg)^{p-1}Kf^q\,d\sigma\\     &\le     C_p\left(\int    g^p
d\sigma\right)^{1/p}       \left(\int       (Kg)^p       f^q
d\sigma\right)^{1/q}.\endalign $$  Since the  left-hand side
is   finite and  $\|g\|_p=1$   we   can   now   cancel  and
 obtain  $$  \int
(Kg)^pf^qd\sigma \le C_p^p $$ so that $f^q\in \Cal  W_{p,K}$
and $\|f\|_Y \le C=C(p,\kappa).$

Now if $f\in S=S_{q,K}$ there exists $u\ge f$ with
$Ku^q\le u$
and  hence  $\|f\|_{\Cal  Y}\le  C.$  The  general  case now
follows by homogeneity.\qed\enddemo

\demo{Remark}We  conclude  that  if  $\Cal  Z\neq\{0\}$ then
$\Cal   W\neq\{0\}$   and   $K$   satisfies   the  condition
$L^p\subset \Cal D(K).$ It follows that $W$ and $Y$ are then
both Banach function spaces.\enddemo

It is now possible to characterize the solution space  $\Cal
Z$  by  means  of  weighted  norm inequalities, although the
result is rather technical:

\proclaim{Theorem 4.7}Suppose $f\in L^0_+.$ Let $f_0 = f$,
and define $(f_n)$
inductively by $f_n=Kf_{n-1}^q.$ In order that $f\in \Cal Z$
it is necessary and  sufficient that each
$f_n^q\in\Cal  W_{p,K}$
and that  if $C_n$  is the  least constant  so that  $$ \int
(Kg)^p \, f_n^q \, d\sigma\le   C_n\int   g^p \, d\sigma$$
for all $0 \le g \in L^p$, then
$\displaystyle{\sup_{n \ge 0} \,
C_n^{1/q^n} <\infty}.$\endproclaim

\demo{Proof}It is a direct consequence of Theorem 2.10  that
$f\in\Cal    Z$     if    and     only    if    $\sup_n \,
\|\Cal
A^nf\|_Y^{1/q^n}<\infty.$ (This  proof tacitly  assumes that
$\Cal  Z\neq  \{0\}$;  however  this  case  can  be   argued
similarly since the set of $h$ such that $\|h\|_{\Cal  W}\le
1$ is bounded in $L^0$).\qed\enddemo

\proclaim{Theorem 4.8}Let $\omega\in \Cal M_+(X).$  Consider
the   following   conditions   on   $\omega:$  \newline  (1)
$\omega\in\widetilde{\Cal W}_{p,K}$ i.e.
there is a  constant $C$
so that $$ \int  (Kg)^pd\omega \le C\int g^pd\sigma\tag  4.7
$$  for  all  $0\le  g\in L^p(\sigma).$
\newline (2) $\omega$
satisfies the second testing inequality (4.5) i.e. there  is
a  constant  $C$  so  that  every  ball $B$ we have
$$\int_B (K \omega_B )^q \, d \sigma \le  C \, |B|_\omega.
 \tag  4.8$$
 (3) $\omega$ satisfies the infinitesimal inequality (3.9)
a.e., i.e. for some constant $C$ and $\sigma$-a.e. $x\in X,$
$$
\sup_{a>0}\left\{\int_0^a\frac{|B_t(x)|_{\sigma}}
{t^2}dt\right\}^{p/q}
\left\{\int_a^{\infty}\frac{|B_t(x)|_{\omega}}{t^2}dt\right\}
\le C.  \tag 4.9 $$
Then the  following are  equivalent:
\newline  (i)  $K\omega\in  \Cal Z_{q, K}$.
\newline (ii) For
some     constant     $C$     we    have    $K(K\omega)^q\le
C \, K\omega < \infty$ a.e.
\newline  (iii)   $\omega$  satisfies   (1)  and
(3).
\newline    (iv)    $\omega$    satisfies    (2)     and
(3).\endproclaim

\demo{Remarks} (a) It follows from Proposition 2.7 (3)
that if
$K(K\omega)^q\le q^{-1} p^{1-q}
\, K\omega$
then  $ K\omega \in S_{q, K}$ i.e. the equation
$u = K u^q +   K\omega$
has a solution $u$ such that
$   K\omega \le u \le p \, K\omega.$\newline
(b) Let $\Cal F$ be an  ideal  space of measurable
functions (e.g. $\Cal F= L^r, \, 0<r \le \infty$).
Then by Theorem 4.8 and the preceding remark,
there  exists a solution $u \in\Cal F$
of the equation $u = K u^q +  \epsilon \, f$
for some $\epsilon > 0$
if and only if $ K\omega \in \Cal F \cap \Cal Z_{q, K}$.
\enddemo

\demo{Proof of Theorem 4.8} That  (ii)   implies  (i)   is
 immediate   from
Proposition 2.7.  For (i) implies (iii), we note first  that
Theorem 4.7  implies that  $K\omega\in \Cal  Y_{p,K}$ and so
Proposition 4.5  gives $\omega\in  \Cal W_{p,K}$  i.e.   (1)
holds.  Theorem 3.10 implies (3) holds.  It is trivial  that
(iii) implies  (iv).   Hence it  remains only  to prove that
(iv) implies (ii).

We let  $C$ denote  the constant  in both  (4.8) and  (4.9).
Note that if $r>0,$
$$ (K\omega)^q \le  2^{q-1} \, ((L_r\omega)^q
+ (U_r\omega)^q).$$
Let $d\nu=(K\omega)^q d\sigma$;  for each
$r>0$       let       $d\mu_r=(U^r\omega)^qd\sigma$      and
$d\lambda_r=(L_r\omega)^qd\sigma.$  Suppose  $x\in  X.$  Let
$B_r=B_r(x).$ Since $K(K \omega)^q = K \nu$,  we
have to prove that
$$K \nu (x) =\int_0^{\infty}\frac{|B_r|_{\nu}}{r^2}dr \le C' \,
 K \omega  (x),$$
where $C'$ depends on $C$, $q$, and $\kappa$.

It is easy to see that
$U_r \omega = U_r \omega_{B_{2\kappa r}}$ on $B_r$. Using
this together with (4.8) we have
$$ |B_{r}|_{\mu_r} =  \int_{B_{r}}(U_r\omega)^q \,
d\sigma=
\int_{B_{   r}}(U_r\omega_{B_{2\kappa r}})^q \, d\sigma \le
C|B_{2\kappa r}|_{\omega}.   $$
  Hence
   $$\align\int_0^{\infty}\frac{|B_r|_{\mu_r}}{r^2}dr  &\le
C\int_0^{\infty}\frac{|B_{2\kappa  r}|_{\omega}}{r^2}dr   \\
&= 2C \, \kappa \,  K\omega(x).\tag 4.10\endalign$$

On  the  other  hand
$$ |B_{r}|_{\lambda_r} =  \int_{B_r}(L_r\omega)^q
 \, d\sigma\le
(2\kappa)^q (L_r\omega(x))^q |B_r|_{\sigma}$$
by Proposition
3.4.  Thus
$$ \int_0^{\infty}\frac{|B_r|_{\lambda_r}}{r^2}dr
\le   (2\kappa)^q
\int_0^{\infty}\left(\int_r^{\infty}\frac{|B_t|_{\omega}}
{t^2}dt\right)^q
\frac{|B_r|_{\sigma}}{r^2}dr.$$  Now  we  use integration by
parts   to   replace   the   right-hand   side   $R$  by  $$
R=q(2\kappa)^q\int_0^{\infty}\left(\int_0^t\frac
{|B_{\tau}|_{\sigma}                               }{\tau^2}
d\tau\right)\left(\int_t^{\infty}\frac
{|B_\tau|_{\omega}}{\tau^2}d\tau\right)^{q-1}
\frac{|B_t|_{\omega}}{t^2}dt.$$

At this point the  infinitesimal inequality (4.9) allows  to
estimate
$$  R \le C q (2\kappa)^q
\int_0^\infty \frac{|B_t|_{\omega}}{t^2}dt  =  Cq(2\kappa)^q
K\omega(x).\tag 4.11$$

Combining these two estimates (4.10) and (4.11) and the fact
that   $\nu\le   2^{q-1}(\lambda_r+\mu_r)$   gives  that  $$
\int_0^{\infty}\frac{|B_r|_{\nu}}{r^2}dr \le  C' \,
K\omega(x)$$
for a suitable constant $C'=C \, C''(q,\kappa).$
Thus (iv) implies
(ii) as claimed.\qed\enddemo

\proclaim{Theorem 4.9}Suppose $K$ satisfies the quasi-metric
condition.    Suppose  $f\in  L^0_+.$  Then  the following
conditions are equivalent:  \newline (1) $f\in \Cal Z_{q,K}$
 i.e.  there  exists  $\epsilon>0$  such  that  there  is  a
solution of  $u=Ku^q+\epsilon f.$\newline  (2) There  exists
$C>0$ such that $K(Kf^q)^q\le C \, Kf^q < \infty$ a.e.
\newline (3)  $f\in
\Cal Y_{p,K}$   and   the   measure   $\omega$   given  by
$d\omega=f^qd\sigma$ satisfies the infinitesimal  inequality
(4.9).\newline   (4)   The   measure   $\omega$  defined  by
$d\omega=f^qd\sigma$  satisfies  both  the testing condition
(4.8) and the infinitesimal condition (4.9).  \endproclaim

\demo{Remarks} (a) By Proposition 2.7 (4) it follows that
if
$K(K f^q)^q\le q^{-q} p^{q(1-q)}
\, K f^q$
then  $ f \in S_{q, K}$ i.e. the equation
$u = K u^q +  f$
has a solution $u$ such that
$$  f+  K f^q \le u \le f + p^q \, K f^q.$$
(b) Theorem 4.9 and the preceding remark yield
the following criterion of the existence of solutions
to the equation $u = K u^q +  \epsilon \, f$
 belonging to some ideal  space of measurable
functions $\Cal F$ (e.g. $\Cal F= L^r, \, 0<r \le \infty$):
  There exists a solution
$u \in \Cal F$ (for some $\epsilon \ge 0$)
if and only if $f, K f^q  \in \Cal Z_{q,K} \cap \Cal F$.
\enddemo

\demo{Proof of Theorem 4.9} This is simply a restatement
of Theorem 4.8 once
one makes  the observation  that $f\in\Cal Z_{q,K}$ if  and
only  if  $Kf^q\in \Cal Z_{q,K}$  (either Proposition 2.7 or
Theorem 2.10 ).\qed\enddemo

\proclaim{Theorem 4.10} Suppose    $K$    satisfies    the
quasi-metric condition and that $x\in X$ and $a>0.$ Then the
following  conditions  are  equivalent:   \newline (1)
$\Cal Z_{q,K}\neq \{0\}.$ \newline
 (2) We have $N(x,a)<\infty$ and
$\sup_{r\ge a}r^{-q/p}M^*(x,r)<\infty.$ \endproclaim

\demo{Proof}Assume  (1).    Then  $\Cal  W\supset  \Cal   Z$
contains  a   strictly  positive   function  and   hence  by
Proposition  4.1,  we  have  $N(x,a)<\infty.$  To  prove the
second condition we  can replace $a$  by any $b\ge  a$ where
$|B_b(x)|_{\sigma}>0.$    Let    us    therefore     suppose
$|B_a(x)|_{\sigma}>0.$  Then  $\chi_{B_a(x)}\in  \Cal Z.$ It
follows  that  we  have  the  infinitesimal  inequality  for
$\chi_{B_a(x)}$. Thus
$$ \sup_{y\in
X}\sup_{r>0}M(y,r)^{p/q}\int_r^{\infty}\frac{|B_a(x)\cap
B_t(y)|_{\sigma}}{t^2}dt  \le  C.$$
Assume  $r>a$ and $y\in
B_r(x).$ Then  for $t\ge  2\kappa r$  we have $B_t(y)\supset
B_a(x).$ Hence  $$ M(y,r)^{p/q}  |B_a(x)|_\sigma
 \le  2C \, \kappa \, r.$$
This implies the second part of (2).

Conversely  assume  (2).    Then  again  we  can assume that
$|B_a(x)|_{\sigma}>0.$  By  Proposition  4.1,  we  have that
$\chi_{B_a(x)} \in\Cal W.$  Now we verify  the infinitesimal
inequality
$$ \sup_{r> 0, y \in X} M(y,r)^{p/q} \,
\int_r^{\infty}\frac{|B_a(x)\cap
B_t(y)|_{\sigma}}{t^2}dt < \infty.$$

  Let $C=\sup_{r\ge a}r^{-q/p}M^*(x,a)$.
Suppose first $y \in B_{2 \kappa a} (x)$.
Then if $r \le a$,
 $$\aligned \int_r^{\infty}
\frac{|B_a(x)\cap B_t(y)|_{\sigma}}{t^2}dt
&\le  \int_0^{a}\frac{|B_t(y)|_{\sigma}}{t^2}dt +
\int_a^{\infty}\frac{|B_a(x)|_{\sigma}}{t^2}dt\\
&\le M(y,a)+\frac{|B_a(x)|_{\sigma}}{a}\le C(\kappa) \,
M^*(x,2 \kappa a).\endaligned$$
Hence   $$   M(y,r)^{p/q}\int_r^{\infty}   \frac{|B_a(x)\cap
B_t(y)|_{\sigma}}{t^2}dt \le C(\kappa, p) \,
M^*(x,2 \kappa a)^{p}.$$
If  $r>a$
$$   M(y,r)^{p/q} \, \int_r^{\infty}\frac{|B_a(x)\cap
B_t(y)|_{\sigma}}{t^2}dt   \le
M(y,r)^{p/q}|B_a(x)|_{\sigma}r^{-1}  \le
C^{p/q}|B_a(x)|_{\sigma}. $$

Now if $y \notin B_{2 \kappa a} (x)$ we set
$b= \rho (x, y) /(2 \kappa)$.
It follows from the quasi-metric inequality that
$B_a(x) \cap B_t(y) = \emptyset$ if $t < b$. Hence
$$\int_r^{\infty}\frac{|B_a(x)\cap
B_t(y)|_{\sigma}}{t^2}dt \le \int_{b}^{\infty}
\frac{|B_a(x)|_{\sigma}}{t^2}dt =  b^{-1}
\,  |B_a(x)|_{\sigma}.$$
Since  $y \in B_{4 \kappa b} (x)$, we have
$$ \aligned M(y,r)^{p/q} \, \int_r^{\infty}\frac{|B_a(x)\cap
B_t(y)|_{\sigma}}{t^2}dt  &\le b^{-1} \,
M^*(x, 4 \kappa b)^{p/q}
 \,  |B_a(x)|_{\sigma}
\\ &\le 4 \, \kappa \, C^{p/q} \,
 |B_a(x)|_{\sigma}.\endaligned$$

Thus the  infinitesimal
inequality holds and so $\chi_{B_a(x)}\in \Cal Z.$ But  this
means   that   $\Cal Z$   contains   a  strictly  positive
function.\qed\enddemo

There is an implicit two-sided estimate of
$\|\chi_{B_a(x)}\|_{\Cal Z}$
in the preceding theorem:

\proclaim{Corollary 4.11} Suppose $K$    satisfies    the
quasi-metric condition and that $x\in X$ and $a>0.$
Let
$$\Phi (x, a) =
|B_a(x)|_{\sigma}^{1/q} \, N(x,a)^{p/q^2} +
 M^*(x,a)^{p/q}$$
$$ + |B_a (x)|_\sigma^{1/q} \, \sup_{r\ge a} \, r^{-1/q} \,
M^*(x,r)^{p/q^2}.$$
Then
$$C_1 \, \Phi (x, \frac a {2 \kappa}) \le
\|\chi_{B_a(x)}\|_{\Cal  Z} \le C_2 \, \Phi (x,
2 \kappa a),$$
where $C_1$ and $C_2$ depend only on $\kappa$ and $p$.
\endproclaim

\demo{Remarks} (a) It is  not difficult  to see  that if
$\sigma$ is  doubling then
   $$  \|\chi_{B_a(x)}\|_{\Cal  Z}   \asymp
|B_a(x)|_{\sigma}^{1/q} \, N(x,a)^{p/q^2} +
 M^*(x,a)^{p/q}$$
$$ + |B_a (x)|_\sigma^{1/q} \, \sup_{r\ge a} \, r^{-1/q} \,
M^*(x,r)^{p/q^2},$$
where the constants of equivalence depend only on $\kappa$
and $p$.\newline
(b) It follows from Corollary 4.11 that
$ \| \bold 1\|_{\Cal  Z}
\asymp \|K \bold 1 \|^{p/q}_{L^\infty}$. \enddemo

\vskip10pt\heading{5.  Capacitary inequalities and  criteria
of solvability} \endheading \vskip10pt

In the previous section we saw that if $f\ge 0$ is such that
there is a positive solution of the equation  $$u=Ku^q+f\tag
5.1$$ then  (if $K$  is a  quasi-metric kernel)  there is  a
corresponding    weighted    norm    inequality    $$   \int
(Kg)^p \, f^q \, d\sigma \le C \int g^p \,
d\sigma,  \tag 5.2$$
for all $g \ge 0$. This can
be  rewritten   in  terms   of  function   spaces  as
$\Cal Z_{q,K}\subset \Cal Y_{p,K}.$

It  is  natural  to  ask  for  a  converse  result, i.e. for
conditions on $K$ so that the weighted norm inequality (5.2)
implies  that  there  exists  $\epsilon>0$  such that
  the equation
$u=Ku^q+\epsilon f$  has a  positive solution.   This  would
imply  that  $\Cal  Y_{p,K}= \Cal Z_{q,K}.$  The aim of this
section is establish such conditions.

We  will  use  the   notion  of  capacity  associated   with
$(K,p,\sigma).$ We define
$$\Cap  E=
\text{Cap}_{p,K}(E)  =  \inf \, \left
\{\int  g^pd\sigma:    \   Kg\ge
\chi_E,\  g\ge  0 \right\}$$
 where  we  stress  that  we   require
$Kg(x)\ge  1$  for  every  (not  almost every) $x\in E.$ The
theory of these  capacities and the  corresponding potential
theory (which is  usually called nonlinear  potential theory
in the nonclassical case $p \not= 2$) is due to B.  Fuglede,
N.  G.  Meyers,  V.  G.  Maz'ya  and  V.  P.  Havin, Yu.  G.
Reshetnyak; a  weighted theory,  for $d  \sigma =  v \, dx$,
where $v \in A_\infty$, was developed by D. Adams  \cite{1}
(see also \cite{2},
\cite{24}, \cite{35}, \cite{57},  and
the literature cited there).

If $\omega\in\Cal M_+(X)$ then the {\it capacity  condition}
$$ |E|_{\omega}\le C\,\Cap E \tag 5.3$$ for every Borel  set
$E$  is  easily  seen  to  be  equivalent  to  the weak-type
weighted     norm      inequality,
$$ |\{Kg(x) \ge \lambda\}|_{\omega} \le
 C\lambda^{-p}\left(\int
g^pd\sigma\right)^{1/p} \tag  5.4$$ for  all $g\ge  0.$
In general (5.3) does not  imply
the corresponding strong-type inequality
$$   \int
(Kg)^p \, f^q \, d\omega \le C \int g^p \, d\sigma,
\tag 5.5$$
 but  in certain  cases this  is true.   In particular
(5.3) implies  (5.5)
 for  Riesz potentials  or more  general
convolution  operators  with  radial  decreasing  kernels on
$\bold R^n$ if $\sigma\in A_{\infty};$ see \cite{2},  \cite {21}
and \cite {34}.  In this section we will find some  new
classes  of  measures  $\sigma$  for  which this implication
holds (see also \cite{52}).

We will also introduce another concept of capacity given  by
$$  \ca   E=  \inf \,
\{\int   g^pd\sigma:     \  Kg\ge  \chi_E, \
\sigma-\text{a.e.},\ g\ge 0\}.$$ It is clear that $\ca E\le
\Cap  E$;  in  fact  $\ca  E=\inf \,
\{\Cap  F:\  F\subset  E,\
|E\setminus  F|_{\sigma}=0\}.$  We   also  note  that   $\ca
E=\|\chi_E\|_{W'}$ by applying the results of Section 2.  If
$\omega$ is absolutely  continuous with respect  to $\sigma$
then (5.3) is equivalent to $$ |E|_{\omega} \le C\,\ca E$$
for every Borel set $B.$

We first  observe that  the capacity  condition implies  the
second testing inequality.

\proclaim{Proposition 5.1}If  $\omega$  satisfies  the
capacity
condition  (5.3)   then  $\omega$   satisfies  the   testing
condition (4.8).\endproclaim

\demo{Proof}By the  remarks above,  and duality  (\cite{31})
$K:L^{q,1}(\omega)\to   L^q(\omega)$   is   bounded    where
$L^{q,1}$ is the  Lorentz space of  all Borel functions  $f$
such that $$\int_0^{\infty} t^{1/q-1}f^*(t)dt<\infty$$ (here
$f^*$   is   the   decreasing   rearrangement  of  $f$). One
immediately   obtains   (4.5)   by   applying   $K$   to   a
characteristic function.\qed\enddemo

\proclaim{Proposition 5.2} Suppose $E$ is a Borel subset  of
$X$.    Assume  that  $L^p\subset  \Cal  D(K)$ so that $\Cal
W_{p,K}$ is a Banach function  space on $X.$ Then given  any
Borel  set  $E$  with  $\ca  E<\infty$ and any $\epsilon>0$
there exists $w\in L^1(E,\sigma)$ with $\int_Ew\,d\sigma=\ca
E$  and  so  that   we  have  $$  \int   (Kf)^pw\,d\sigma\le
(1+\epsilon)\int f^p\,d\sigma$$ for all $f\in  L^p(\sigma).$
\endproclaim

\demo{Proof}This follows from our earlier identification  of
$\Cal W'=\Cal V$  (see Section 2).   From the  definition of
$\ca E$  we have  $\|\chi_E\|_{\Cal V}=\ca  E$ and  so there
exists $w\in \Cal W$ with $\|w\|_{\Cal W}\le 1+\epsilon$ and
$\int w\chi_E=\ca E$.  Replacing $w$ by $w\chi_E$ gives  the
result.\qed\enddemo

Our next result is closely related to the concept of
an  equilibrium  measure  (see,  for  example   \cite{2}
in the case of Riesz potentials and $d \sigma = d x$).

\proclaim{Proposition  5.3} Let  us  assume  that $K$ is  a
quasi-metric  kernel  such  that  $\rho=K^{-1}$ is
 continuous.
Suppose also that $\sigma$ is locally finite.  Then for  any
compact  set  $E$  with  $\Cap  E<\infty$  there  is a Borel
measure    $\omega$    supported    on    $E$    such   that
$|E|_{\omega}=\Cap E$  and we  have the  inequality $$  \int
(Kf)^pd\omega  \le   \int  f^pd\sigma$$   for  every   $f\in
L^p(\sigma).$\endproclaim

\demo{Proof} Fix  any $x_0\in  X.$ For  all $r>0$  and $s>0$
it is easily seen
that if $f\in L^{p}(B_r(x_0))$ then $L_sf$ is continuous  on
$E.$ Since $L_s(x,y)$ is continuous and bounded this follows
immediately from the Dominated Convergence Theorem.

Now consider  the convex  set $F$  of all  functions in  the
space  of  continuous  real-valued  functions  on $E,$ $\Cal
C(E)$ of the form $\sum_{i=1}^n (L_{s_i}f_i)^p-\chi_E$ where
$0<s_i<a,$ $f_i\in B_{r_i}(x_0)$ for some $a<r_i<\infty$ and
$\sum\int f_i^pd\sigma \le \Cap E.$ Let $P=\{f\in \Cal C(E):
\  f(x)>0\  \forall  x\in  E\}.$  We claim that $F\cap P=0.$
Indeed  if  not  there  exist  $s_i,f_i$  as  above with and
$\epsilon>0$  so  that  $$  \sum_{i=1}^n  (L_{s_i}f_i)^p \ge
(1+\epsilon)\chi_E.$$ But then $$ K(\sum_{i=1}^nf_i^p)^{1/p}
\ge             (\sum_{i=1}^n(Kf_i)^p)^{1/p}             \ge
(1+\epsilon)^{1/p}\chi_E$$ which contradicts the  definition
of capacity.

Now since the cone $P$  is open the Hahn-Banach theorem  and
Riesz  representation  theorem  combine  to  give  a measure
$\omega$ such that $\int f\,d\sigma>0$ for all $f\in P$  and
$\int f\,d\sigma\le 0$ for all $f\in F.$ If we normalize  so
that $|E|_{\omega}=\Cap E$ then we have $$ \int (L_sf)^p \le
1$$  whenever  $s>0$,  $f$  has  bounded  support  and $\int
f^pd\sigma\le 1.$ This implies the result.\qed\enddemo

Now let us introduce for $x\in X$ and $a>0$ the quantity  $$
N(x,a) =\int_a^{\infty}\frac{|B_t(x)|_{\sigma}}{t^{1+q}}dt =
\int_X L_a(x,y)^qd\sigma(y).$$

\proclaim{Definition}We shall say that  a locally finite
Borel measure  $\mu$ is
{\it stable with respect to $K$} if there exists
 a constant $C>0$
and  $\delta>0$  so  that   and  for  every  $x\in   X$  and
$0<r<R<\infty$  we have
$$ |B_{2 R} (x)|_{\mu} \le C \, |B_R(x)|_{\mu} \quad
\text{and} \quad |B_r(x)|_{\mu} \le
C \, \left(\frac{r}{R}\right)^{1 + \delta}|B_R(x)|_{\mu}.$$
\endproclaim

We remark  that of  course if  $K$ is  a Riesz  potential on
$\bold   R^n$   (i.e.     $K(x,y)=c\|x-y\|^{-\alpha}$  where
$0<\alpha<n$) then Lebesgue measure $\lambda$ is stable  for
$K.$ More generally  if $X$ is  a space of  homogeneous type
with  quasi-metric  $d$  and  a  doubling measure $\mu$ then
$\mu$   is   stable   for   $K(x,y)=d(x,y)^{-\alpha}$  where
$0<\alpha<1$ (\cite{20}).

\proclaim{Theorem   5.4}    Suppose   $K$    satisfies   the
quasi-metric condition.  \newline(1) There is a constant $C$
so  that  for  any  for  any  ball  $B_a(x)$ we have $$ \Cap
B_a(x)\le C N(x,a)^{-p/q}.\tag  5.6$$  \newline  (2) Suppose
that there  exists a  measure $\mu\in\Cal  M_+(X)$ which  is
stable for $K.$ Then there exists a constant $C$ so that for
every  ball  $B=B_a(x)$  we  have $$ C^{-1}N(x,a)^{-p/q} \le
\Cap B_a(x)  \le CN(x,a)^{-p/q}.   \tag 5.7$$  (3) Suppose
that there exists a $\sigma$-continuous measure $\mu$  which
is stable for $K$.  Then there exists a constant $C$ so that
for  every  ball  $B=B_a(x)$  we have
$$ C^{-1}N(x,a)^{-p/q}
\le\ca B_a(x)\le \Cap B_a(x)  \le C \,
N(x,a)^{-p/q}. \tag 5.8
$$ \endproclaim

\demo{Remarks}In  particular  we  obtain  equivalence  $\Cap
B_a(x)\sim N(x,a)^{-p/q}$ for the case $\bold R^n$ and any
 measure
$\sigma;$ we obtain the further equivalence $\ca  B_a(x)\sim
N(x,a)^{-p/q}$ for any $\sigma$ of the form $d\sigma=w(x)dx$
where $w(x)>0$ a.e.  For measures in class $A_{\infty}$ this
result was  previously shown  by Adams  \cite{1}; the upper
estimate (5.6) was proved  for arbitrary measures by
Turesson
\cite{51} in the case of Riesz potentials.
  For  the  case  $d\sigma=dx$  and more general
radial convolution operators, a similar result was obtained
by Aikawa  (see \cite{4}).  \enddemo

\demo{Proof}(1)  Let   $g(y)=L_a(x,y)^{q-1}.$  Then   $$\int
g(y)^pd\sigma(y)  =\int  L_a(x,y)^q \,
d\sigma(y)  =N(x,a).$$
 We
also have
$$ Kg(x)  \ge L_a g(x)=\int  L_a(x,y)^q \,
d\sigma(y) = N(x,a).$$
Applying   Proposition   3.4    we   obtain   $$    Kg   \le
(2\kappa)^{-1}N(x,a)\chi_{B_a(x)}.$$  Hence  $\Cap B_a(x)\le
(2\kappa)^p \, N(x,a)^{-p/q}.$

(2)  Suppose  $g\in  L^p(\sigma)$  and  $Kg\ge \chi_B$ where
$B=B_a(x).$     Then     $$     L_ag(x)    \le    \left(\int
L_a(x,y)^qd\sigma(y)\right)^{1/q}\left                 (\int
g(y)^pd\sigma(y)\right)^{1/p}.$$    Then    $$    L_ag(x)\le
N(x,a)^{1/q}\|g\|_p$$ and thus if $y\in B_a(x),$ $$  L_ag(y)
\le     2\kappa     N(x,a)^{1/q}\|g\|_p.$$     Hence      $$
\int_{B}L_ag\,d\mu     \le      2\kappa|B|_{\mu}N(x,a)^{1/q}
\|g\|_p.\tag 5.9$$

On  the  other   hand
$$\align  \int_{B}U_a g \, d\mu   &=\int
U_a \mu_B \, g \,d\sigma\\ &\le \int_{B^*} U_a \mu  \,
 g \, d\sigma\\
&\le  \left(\int_{B^*} (U_a\mu)^q \, d\sigma \right)^{1/q}
\|g\|_p,\endalign $$
where $B^* = B_{2 \kappa a} (x)$.

Now  we estimate
$$ U_a\mu(y) =\int_0^a\frac{|B_t(y)|_{\mu}}{t^2}dt  \le
C_1 \, |B_a(y)|_{\mu}$$
for  a  suitable  constant  $C_1$ since
$\mu$  is   stable  for   $K.$  But   if  $y\in   B^*$  then
$B_a(y)\subset     B_{\kappa(2\kappa+1)a}(x)$     so    that
$|B_a(y)|_{\mu} \le  C_2 |B|_{\mu}$  again by  the stability
condition.   Hence
$$     \int_B U_a g\,d\mu    \le
C_3 \, |B|_{\mu} \, |B^*|_{\sigma}^{1/q} \, \|g\|_p.$$
Observe    that
$|B^*|_{\sigma}^{1/q}    \le    C_4    N(x,a)^{1/q}$   where
$C_4=C_4(\kappa).$ Combining with (5.9) we obtain
$$  \int_B Kg \,d\mu  \le  C_5  |B|_{\mu} \, N(x,a)^{1/q}
\, \|g\|_p $$
which  gives  the   lower  estimate  since   $Kg\ge  \chi_B$
(everywhere).

(3)  The  proof  is  the  same  as  (2) except for the final
observation note that if $Kg\ge\chi_B$   $\sigma$-a.e. then
$Kg\ge\chi_B$ $\mu$-a.e.\qed\enddemo

Before moving  to our  main result  of this  section, let us
prove a preliminary lemma.

\proclaim{Lemma  5.5}Suppose  $K$  is  a quasi-metric kernel
with the property that there  is a constant $C$ so  that for
every $x\in X$ and $a>0$ we have $M(x,a)\le Ca^{q-1}N(x,a).$
Then there exists a constant $C'$ so that for every $x\in X$
and           $a>0$           we           have           $$
M(x,a)^{p/q}\int_a^{\infty}\frac{N(x,t)^{-p/q}}{t^2}dt   \le
C'.\tag 5.10 $$ \endproclaim

\demo{Proof}First we  observe if  $M(x,a)\le Ca^{q-1}N(x,a)$
then       $$       \int_0^a\frac{|B_t(x)|}{t^2}dt       \le
Ca^{q-1}\int_a^{\infty} \frac{|B_t(x)|}{t^{1+q}} dt  \le
C\int_a^{\infty}\frac{|B_t(x)|_{\sigma}}{t^2}dt,$$  so  that
for every $x\in X$ we have $\lim_{a\to\infty}M(x,a)=\infty.$
For fixed $x\in X$ and  $a>0$ let us define $a^*>a$  so that
$$M(x,a^*)=(1 +\frac1{2C})M(x,a).$$

We claim  now that  $M(x,a) \le  2Ca^{q-1}N(x,a^*).$ To  see
this note  that $$  \align M(x,a)  &\le Ca^{q-1}N(x,a)\\  &=
Ca^{q-1}\int_a^{a^*}\frac{|B_t(x)|_{\sigma}}{t^{1+q}}dt    +
Ca^{q-1}N(x,a^*)\\      &\le      C(M(x,a^*)-M(x,a))       +
Ca^{q-1}N(x,a)\\       &=       \frac12       M(x,a)       +
Ca^{q-1}N(x,a).\endalign $$ The estimate then follows.

Now for fixed  $a$ define a  sequence $(a_j)_{j=0}^{\infty}$
inductively by $a_0=a$ and then $a_j=a_{j-1}^*.$ We have  $$
\align
M(x,a)^{p/q}\int_a^{\infty}\frac{N(x,t)^{-p/q}}{t^2}dt    &=
M(x,a)^{p/q}
\sum_{j=0}^{\infty}\int_{a_j}^{a_{j+1}}
\frac{N(x,t)^{-p/q}}{t^2}dt\\
&                                                        \le
M(x,a)^{p/q}\sum_{j=0}^{\infty}N(x,a_{j+1})^{-p/q}a_j^{-1}\\
&\le                 (2C)^{p/q}                 M(x,a)^{p/q}
\sum_{j=0}^{\infty}M(x,a_j)^{-p/q}\\      &=      (2C)^{p/q}
\sum_{j=0}^{\infty}(1+\frac{1}{2C})^{-pj/q}\\  &=C'\endalign
$$  say.     This   completes  the   proof  of   the  lemma.
\qed\enddemo

We can now state our main theorem of the section:

\proclaim{Theorem  5.6}Let  $K$  be  a  quasi-metric kernel.
Assume that $N(x,a)<\infty$ for  all $x\in X$ and  $a>0$ and
that there is a constant $C$ so that for every $x\in X$  and
$a>0,$   $$   \int_0^a\frac{|B_t(x)|_{\sigma}}{t^2}dt    \le
Ca^{q-1}\int_a^{\infty}\frac{|B_t(x)|_{\sigma}}{t^{1+q}}dt
\tag
5.11$$  (i.e.     $M(x,a)\le   C a^{q-1} N(x,a)$).
 Suppose
$\omega\in\Cal M_+(X).$  Then $\Cal Z_{q,K}\neq \{0\}$  and
the  following  statements  are  equivalent:    \newline (1)
$\omega\in \widetilde{\Cal W}_{p,K}$ i.e.  $\omega$
satisfies  the
weighted   norm   inequality   $$   \int  |Kg|^p \,
d\omega  \le
C\int|g|^p \,
d\sigma$$ for all $g\in L^p(\sigma).$\newline  (2)
$\omega$ satisfies  the capacity  condition $|E|_{\omega}\le
\Cap E$ for all Borel sets $B.$ (Equivalently the  weak-type
inequality (5.4)  holds.)   \newline (3)  $\omega$ satisfies
the    testing     condition
$$  \int
(K \omega_B)^q \,
d\sigma\le C|B|_{\omega}$$  for all balls $B.$
\newline  (4)  $K\omega\in \Cal Z_{q,K}$  i.e.  for   some
$\epsilon>0$  there  is  a  solution  $u$  of  the  equation
$u=Ku^q+\epsilon K\omega.$ \newline (5) There is a  constant
$C$ so that $K(K\omega)^q\le CK\omega.$ \endproclaim

\demo{Proof}Let us  first prove  the nondegeneracy  of
$\Cal Z_{q,K}.$ Fix  any $x\in  X$ and  $a>0.$ If
$r>a$ and $y\in
B_r(x)$ then $M(y,r)\le Cr^{q-1}N(y,r) \le  C'r^{q-1}N(x,r)$
where  $C,C'$  do  not  depend  on  $r.$ Hence $M^*(x,r) \le
C'r^{q/p}N(x,a)$ and we can apply Theorem 4.10.

The equivalence of (4) and (5) is proved in Theorem 4.7.  We
have seen that  (1) implies (2)  and (2) implies  (3) (Lemma
5.1).  We also have that (5) implies (1) by Theorem 4.8.  It
remains to show that (3)  implies (5).  For this  by Theorem
4.8  we  need  only  establish  the infinitesimal inequality
(4.9).  To this end, if $B$ is any ball, note that if  $0\le
f\in L^p(\sigma)$  and $Kf\ge  \chi_B$ then  $$ |B|_{\omega}
\le  \int  \chi_BKf  \,d\omega=\int (K\omega_B)f\,d\sigma.$$
Hence $$ |B|_{\omega} \le C^{1/q}|B|_{\omega}^{1/q}\|f\|_p$$
by  the   testing  condition.     Thus   $$  |B|_{\omega}\le
C_1\|f\|_p^p$$ for a suitable $C_1.$ However we then  deduce
that   $$    |B|_{\omega}\le   C_1\Cap    B.$$   Hence    $$
\align\int_a^{\infty}\frac{|B_t(x)|_{\omega}}{t^2}dt    &\le
C_1\int_a^{\infty}    \frac{\Cap    B_t(x)}{t^2}dt\\    &\le
C_2\int_a^{\infty}\frac{N(x,t)^{-p/q}}{t^2}dt\\   &\le   C_3
M(x,a)^{-p/q}\endalign $$ by Theorem 5.4 and Lemma 5.5.  But
this implies that
$$\left(
\int_0^a\frac{|B_t(x)|_{\sigma}}{t^2}dt \right)^{p/q}\left (
\int_a^{\infty}\frac{|B_t(x)|_{\omega}}{t^2}dt\right)    \le
C_4$$ i.e. the infinitesimal inequality holds.\qed\enddemo

Let us remark at this point that (5.11) can be  equivalently
formulated         as         $$         M(x,a)          \le
C|B_a(x)|_{\sigma}^{1/p}N(x,a)^{1/q}.\tag 5.12$$ That (5.12)
implies (5.11) follows from the fact that $|B_a(x)|_{\sigma}
\le    q \, a^q \, N(x,a).$    Now    assume    (5.11),
  and   let
$\theta=aM(x,a)/\alpha|B_a(x)|_{\sigma},$              where
$\alpha=\max(2,4C/q).$ Then $$ M(x,\theta^{-1}a) \ge  M(x,a)
- (\theta-1)a^{-1}|B_a(x)|_{\sigma}\ge \frac12M(x,a).$$  Now
$$          \align          M(x,\theta^{-1}a)           &\le
C\theta^{q-1}a^{q-1}N(x,\theta^{-1}a)\\ &\le  C\theta^{-q+1}
a^{q-1}N(x,a) + \frac{C}q\theta a^{-1}|B_a(x)|_{\sigma})  \\
&\le C\alpha^{q-1}M(x,a)^{1-q}|B_a(x)|_{\sigma}^{q-1}N(x,a)+
\frac{C}{q\alpha}               M(x,a)\\                &\le
C\alpha^{q-1}M(x,a)^{1-q}N(x,a) +\frac14 M(x,a).   \endalign
$$           Thus           $$           M(x,a)          \le
4C\alpha^{q-1}M(x,a)^{1-q}|B_a(x)|^{q-1}N(x,a)$$  which  can
be reorganized as (5.12).

\proclaim{Theorem 5.7}Under the  same hypotheses as  in
Theorem
5.6, we have  $\Cal Y_{p,K}=\Cal Z_{q,K}$  or, equivalently,
the following conditions  on $f\in L^0_+$  are equivalent:
\newline (1) There exists $\epsilon>0$ so that the  equation
$u=Ku^q+\epsilon f$ has  a solution $u\in  L^0_+.$\newline
(2) There is a constant $C$ so that
$$ \int (Kg)^pf^qd\sigma
\le      C\int      g^pd\sigma$$
   for      all     $g\in
L^p(\sigma).$\endproclaim

\demo{Proof}Just observe  that $f\in\Cal  Z$ if  and only if
$Kf^q\in\Cal Z$ and apply the equivalence of (1) and (4)  in
Theorem 5.6.\qed\enddemo

We now turn to the problem of converses.  We will show  that
under   mild   conditions   (5.11)   is  necessary  for  the
equivalences of (1) and (4) of Theorem 5.6 or of (1) and (2)
in Theorem 5.7.

\proclaim{Theorem 5.8}(1) Let  $K$ be a  quasi-metric kernel
and suppose there is a stable measure $\mu$ for $K$ which is
$\sigma-$continuous.  Then if $\Cal Y_{p,K}=\Cal Z_{q,K}\neq
\{0\}$ (i.e.   (1)  and (2)  of Theorem  5.7 are equivalent)
then  (5.11)  holds  i.e.  for  some  $C$ we have $M(x,a)\le
Ca^{q-1}N(x,a)$ for  all $x\in  X$ and  $a>0.$ \newline  (2)
Suppose $K$  is a  continuous quasi-metric  kernel such that
each ball  $B_a(x)$ is  compact.   Suppose that  there is  a
stable  measure  $\mu$  for  $K.$  Then if for every measure
$\omega\in \widetilde {\Cal W}_{p,K}$  we  have
 $K\omega\in
\Cal Z_{q,K}$ then (5.11) holds.  \endproclaim

\demo{Proof}(1) The hypotheses imply that $\Cal Y_{p,K}=
\Cal Z_{q,K}$ are Banach  function spaces with
equivalent norms.
Hence for some  constant $C$ and  any Borel set  $E$ we have
$$\|\chi_E\|_{\Cal Z}\le C\|\chi_E\|_{\Cal
Y}=\|\chi_E\|_{\Cal W}^{1/q}  \le
|E|^{1/q}_{\sigma}\|\chi_E\|_{\Cal
W'}^{-1/q}=|E|_{\sigma}^{1/q}(\ca E)^{-1/q}.$$  If we  apply
this to a ball $B_a(x)$ using  Theorem 3.8 we obtain
$$  M(x,a)
\le  C_1 \, |B_a(x)|_{\sigma}^{1/p}  (\ca B_a(x))^{-1/p}$$
and
hence    by    Theorem    5.4    $$    M(x,a) \le C_2 \,
|B_a(x)|_{\sigma}^{1/p} \, N(x,a)^{1/q} \le
C_3 \, a^{q-1} \,  N(x,a).$$

(2) The hypotheses can easily be seen to imply the
existence
of a constant $C$ so  that if $\omega$ is any  Borel measure
such that  $$ \int  (Kf)^pd\omega \le  \int f^pd\sigma  \tag
5.12$$ for all $f\in L^p$ then the infinitesimal  inequality
(3.9)    holds    with    constant    $C$ i.e.
$$\sup_{a>0}\left\{\int_0^a\frac{|B_t(x)|_{\sigma}}
{t^2}dt\right\}^{p/q}
\left\{\int_a^{\infty}\frac{|B_t(x)|_{\omega}}{t^2}dt\right\}
\le C.$$
Now by Proposition  5.3, we can find  a
measure    $\omega$     supported    on     $B_a(x)$    with
$$|B_a(x)|_{\omega}=\Cap B_a(x)$$
  and so  that (5.12)  holds.
Then we  have $$  M(x,a)^{p/q} \,
\Cap B_a(x)  \le C_1.$$  Again
appealing to Theorem 5.4 gives the result.\qed\enddemo

It may happen that the equivalence of (1) and (4) of Theorem
5.6 or (1) and (2)  of Theorem 5.7 hold even  when condition
(5.11) fails, however.   Of course  this can only  happen if
there is no stable measure for $K.$

\proclaim{Theorem 5.9}Suppose $K$  is a quasi-metric  kernel
and that $\sigma$ satisfies the conditions (5.13) and (5.14)
for some constant $C$:
$$ M(x,2a) \le C \, M(x,a) \tag 5.13  $$
and $$ M(y,a)  \le C \, M(x,a)\tag 5.14$$
 whenever $x,y \in  X$
and $\rho(x,y)\le a.$ Suppose also that $N(x,a)<\infty$  for
some $x\in  X$ and  $a>0.$ Suppose  $\omega\in \Cal M_+(X)$.
Then $\Cal Z_{q,K}\neq\{0\},$  and the following  conditions
are equivalent:  \newline
(1) $\omega\in \widetilde{\Cal W}_{p,K}$
i.e.   $\omega$ satisfies  the weighted  norm inequality  $$
\int  |Kg|^pd\omega  \le  C\int|g|^pd\sigma$$  for all $g\in
L^p(\sigma).$\newline  (2)  $\omega$  satisfies the capacity
condition $|E|_{\omega}\le \Cap E$  for all Borel sets  $B.$
(Equivalently   the   weak-type   inequality  (5.4)  holds.)
\newline
(3) $\omega$ satisfies the testing condition  (4.8)
i.e.
 $$ \int (K\omega_B)^q \, d\sigma \le C \, |B|_{\omega}$$
for  all
balls $B.$ \newline
(4)  $K\omega\in \Cal Z_{q,K}$ i.e.  for
some $\epsilon>0$ there  is a solution  $u$ of the  equation
$u=Ku^q+\epsilon K\omega.$
\newline
(5) $K \omega \in L^0_+ (\sigma)$ and there is a  constant
$C$ so that $K(K\omega)^q\le C \, K\omega.$ \endproclaim

\demo{Remark}We  remark  that  (5.13)  is  equivalent to the
requirement      that      we      have      an     estimate
$|B_{2a}(x)|_{\sigma}\le   C \, a \, M(x,a)$   since
$M(x,2a)  \le
(2a)^{-1}|B_{2a}(x)|_{\sigma}+M(x,a)\le    2M(x,4a).$     In
particular it is sufficient that we have a doubling estimate
$|B_{2a}(x)|_{\sigma}\le    C|B_a(x)|_{\sigma}$    for   any
constant $C.$ \enddemo

\demo{Proof}We first establish $\Cal Z_{q,K}\neq\{0\}.$  Fix
$a>0$ and $x\in X$  so that $N(x,a)<\infty.$ Then  for $r>a$
we  have  an  estimate  $|B_r(x)|_{\sigma}\le  Cr^q.$ Now
$$M(x,r) =M(x,a)+\int_a^r\frac{|B_t(x)|}{t^2}dt \le
M(x,a)+C'r^{q-1}.$$
Hence by (5.14)
$$M^*(x,r) \le M(x,a) + C' \, r^{q-1}.$$
Since $1-q = -q/p$, we have $\sup_{r>a}
 r^{-q/p} \,  M^*(x,r) < \infty$,  and Theorem 4.10 applies.

To complete the  proof we argue  as in Theorem  5.6.  It  is
only  necessary  to  show   that  if  (3)  holds   then  the
infinitesimal inequality (4.9) holds.  We note as in Theorem
5.6   that   (3)   implies   the   capacity   condition
$$|B|_{\omega}\le C \, \Cap B$$
for any ball $B.$

Let us  estimate $\Cap  B_a(x).$ In  fact $K\chi_{B_{2\kappa
a}(x)}(y)\ge M(y,a)$ if $y\in  B_a(x)$ so that we  obtain an
estimate  $$  \Cap  B_a(x)  \le  C_1 \,  M(x,a)^{-p} \,
|B_{2\kappa a}(x)|_{\sigma}.$$       Thus  we have
$$\int_a^{\infty}\frac{|B_t(x)|_{\omega}}{t^2}dt  \le
C_2 \, \int_a^{\infty} M(x,t)^{-p} \, \frac{|B_{2\kappa
t}(x)|_{\sigma}}{t^2} dt.$$

Now  we  obviously  have  an  estimate  $M(x,2\kappa   t)\le
C_3 \, M(x,t)$
so that upon substituting $\tau=2\kappa t$ we can
obtain
$$ \int_a^{\infty}\frac{|B_t(x)|_{\omega}}{t^2}dt \le
C_4 \, \int_{2\kappa a}^{\infty} M(x,\tau)^{-p}
\frac{|B_{\tau}(x)|_{\sigma}}{t^2}dt=
C_4 \, (p-1)^{-1} \, M(x,2\kappa a)^{1-p}.$$
 Using  the   estimate
(5.14)  again  this  gives  us  the  infinitesimal condition
(4.9).\qed\enddemo

\vskip10pt\heading{6. Trace inequalities, Carleson measure
theorems,
and nonlinear convolution equations}
 \endheading \vskip10pt

In this section we give some applications
 of  the results of the previous section  to trace
inequalities for Riesz potentials as well as  more general
convolution
operators on $\bold  R^n$, and the solvability problem for
 the corresponding nonlinear convolution
equations. We also obtain a weighted version of the
Carleson measure
theorem for Poisson integrals on $\bold  R^{n+1}_+$ which
generalizes a result of Treil and Volberg \cite{50}
where the case $p=2$ was considered.
(The proof in \cite{50} makes use of a test for boundedness of
 quadratic forms; in the classical unweighted
case this idea is due to S.A. Vinogradov \cite{40}.) Our
approach is closer to the well-known proof of Hardy's inequality
(see e.g. \cite {34}) and works for all  $1<p<\infty$.

\proclaim{Theorem  6.1} Let  $k:\bold  R^n\to\bold  R$  be a
positive  radial  function  of  the form $k(x)=h(\|x\|)$ for
$x\neq    \{0\}$    where    $h$    is    decreasing     and
$$\inf_{r>0}\frac{h(2r)}{h(r)} > 0.$$
Suppose $k\in L^1+L^q.$ Let
$\lambda$ be  Lebesgue measure  on $\bold  R^n$ and  suppose
$\omega$  is  any  locally finite  Borel  measure. Then the
following are equivalent:
\newline (1) There exists  $C$  so  that  we  have
$$ \int (k*f)^p \, d\omega \le C \, \int  f^p \, d\lambda$$
for all $f \ge 0$.\newline
(2) $k*\omega \in L^0_+$ and
 $$k*(k*\omega)^q\le  C'\, k*\omega$$
for some
constant $C'$.\newline
 (3) $k*\omega \in\Cal Z$, or equivalently the equation
$u = k*u^q + \epsilon \, (k*\omega)$ has a solution
for all sufficiently small $\epsilon> 0$.\newline
 (4) The testing inequality
$$\int_B (k*\omega_B)^q \, d \lambda \le C \, |B|_\omega$$
holds for all Euclidean balls $B$ in  $\bold  R^n$.\newline
(5) There exists  $C$  so  that
$$|E|_\omega \le C \Cap_{p, k} (E)$$
for all compact sets $E$.
\endproclaim

\demo{Remarks} (a) Characterizations   of  trace
inequalities in terms of capacities are due to
V. Maz'ya, D. Adams, and B. Dahlberg in the case of Riesz
potentials (see \cite{3}, \cite{35}). For radial kernels similar
 characterizations in terms of capacity inequalities
or testing inequalities with totally different and more difficult
proofs are due to K. Hansson \cite{21} and E. Sawyer and R.
Kerman
\cite{27}.

(b) A characterization of solvability for the equation in (3)
with an arbitrary inhomogeneous term $f \ge 0$ in place of
$k*\omega$ is given by $k*(k*f^q)^q\le  C \, k*f^q.$

(c) We remark that we make use of Theorem 5.9 and not 5.6
which would only apply to special cases of such convolution
operators.

\enddemo

\demo{Proof}If  $K(x,y)=k(x-y)$   then  $K$   satisfies  the
quasi-metric assumption.  The assumptions on $h$ ensure that
$M(x,a)$   and   $N(x,a)$    are   everywhere   finite   for
$\sigma=\lambda.$ In this case $M(x,a)$ is constant for each
$a.$  Note  that  in  this  case  $B_a(0)=\{x:   h(\|x\|)\ge
a^{-1}\}$ and so we have an estimate that  $B_{2a}(0)\subset
CB_a(0)$    which    implies    $|B_{2a}(0)|_{\lambda}   \le
C^n|B_a(0)|_{\lambda}.$  By  the  remarks  following Theorem
5.9,  this  means  we  can  apply  this theorem to yield the
result.  Note that in (5) one can use compact sets in  place
of Borel sets due to the known capacitability results (see
\cite{3}, p. 28).
  \qed\enddemo

Next we characterize trace inequalities for Riesz
potentials  of order $\alpha$,
$I_\alpha = (- \Delta)^{- \alpha/2}$,
on $\bold R^n$, and the solvability
problem for the integral equation
$$u = I_\alpha(u^q \, d \sigma) +
f,  \quad 0 \le u < \infty \quad d \sigma\text{-a.e.}
\tag 6.1$$
Note that in contrast
to Theorem 6.1 now $\sigma$ is not necessarily Lebesgue
measure.

 We will make use of the class of
$A_\infty^\beta$-weights
introduced in \cite{45}, which contains both Muckenhoupt
$A_\infty$-weights (in case $\beta = n$) and (reverse)
doubling
weights $RD_\beta$ such that
$$|B_r|_\sigma \le C \, \left (\frac {r} {R} \right)^\beta
 \, |B_R|_\sigma$$
for all concentric balls $B_r$ and $B_R$ with
$0<r < R< \infty$.
We set  $\Cal Z_{q, \alpha} = \Cal Z_{q, I_\alpha}$ and
$S_{q, \alpha} = S_{q, I_\alpha}$. We also set
$$I_\alpha^\sigma f (x) = I_\alpha (f d \sigma) (x) =
C(n, \alpha) \, \int \frac {f (y) \, d \sigma}
{|x-y |^{n -\alpha}},$$
where
$C(n, \alpha)=
 \pi^{- n /2} \, 2^{-\alpha} \, \Gamma( n /2- \alpha /2) \,
\Gamma (\alpha/ 2)^{-1}$.

\proclaim{Theorem 6.2} Let $1< q < \infty$ and
$0 < \alpha <n$.
Let $\sigma$ and $\omega$ be locally finite measures on
$\bold R^n$, and let $f = I_\alpha \omega \in L^0_+ (\sigma)$.
 Then the following statements are true. \newline
(1) $f \in \Cal Z_{q, \alpha}$ if and only if the inequality
$$ I_\alpha^\sigma (I_\alpha \omega )^q \le C \,
 I_\alpha \omega \quad d \sigma\text{-a.e.} \tag 6.2$$
holds.  Moreover, if (6.2) holds with  $C= p^{1-q} q^{-1}$
then $f \in S_{q, \alpha}$, and
(6.1) has a solution $u$ such that $I_\alpha \omega
\le u \le p \, I_\alpha \omega$.
\newline
(2)  $f \in \Cal Z_{q, \alpha}$  if and only if both the
trace inequality
$$||I^\sigma_\alpha h||_{L^p (\omega)} \le C \,
||h||_{L^p (\sigma)}, \quad h \in L^p (\sigma),
\tag 6.3$$
and the infinitesimal inequality
$$\sup_{x \in \bold R^n, \, r > 0}
  \left \{ \int_0^r \frac {|B_t (x)|_\sigma}
{t^{n- \alpha +1}}
 \, d t \right \}^{1/q}
 \, \left \{ \int_r^\infty
 \frac {|B_t (x)|_{\omega}}
{t^{n- \alpha +1}} \, d t \right \}^{1/p}
< \infty, \tag 6.4$$
hold, where $B_r (x)$ is a Euclidean ball of radius $r$
centered at $x$.
\newline (3) $f \in \Cal Z_{q, \alpha}$    if and only if both
the infinitesimal inequality (6.4) and
the testing inequality
$$\int_B (I_\alpha^\omega  \chi_B)^q \, d \sigma
\le C \, |B|_\omega,
\tag 6.5$$
hold, where $C$ is independent of $B= B_r (x)$.
\newline (4) If $\sigma$ satisfies the estimate
$$\int_0^r \frac {|B_t (x)|_\sigma} {t^{n- \alpha +1}}
 \, d t  \le C \, r^{(n- \alpha) (q-1)} \,  \int_r^\infty
\frac {|B_t (x)|_\sigma} {t^{(n- \alpha)q+1}}  \, d t,
\tag 6.6$$
then  (6.3)$\Leftrightarrow$(6.5)$\Leftrightarrow$(6.2).
Moreover, (6.6) is necessary in order that
(6.3)$\Leftrightarrow$(6.2).
\newline
(5) If $\sigma, \, \omega \in A_\infty^\beta$ with
$\beta > n - \alpha$,
then (6.2) is equivalent to the infinitesimal
inequality (6.4).
\newline
(6) If $\sigma, \, \omega \in RD_\beta$ with
$\beta > n - \alpha$, then (6.2) is equivalent to the
following condition
of Muckenhoupt type,
$$\sup_{x \in \bold R^n, \, r > 0}
 \frac {|B_r (x)|_\sigma^{1/q} \,
 |B_r (x)|_{\omega}^{1/p}}
 {r^{n- \alpha}}
 < \infty. \tag 6.7$$
\endproclaim

\demo{Remarks}
(a) For $\sigma \in RD_\beta$ with $\beta > n - \alpha$
similar results were proved earlier in \cite{52}. Note that in this
case
(6.6) holds and hence (6.4) follows from (6.5) which simplifies the
proofs.
\newline
(b) There are analogous criteria of solvability of (6.1)
for an arbitrary $f \in L^0_+$ in place of $I_\alpha \omega$
 (see Theorem 4.9). In particular,  $f \in
\Cal Z_{q, \alpha}$ if
and only if $$I_\alpha^\sigma (I_\alpha^\sigma f^q)^q
\le C \, I_\alpha^\sigma f^q. \tag 6.8$$
Moreover, if (6.8) holds with $C= p^{q(1-q)} q^{-q}$, then
 $f \in S_{q, \alpha}$
 i.e. (6.1) has a solution $u$ such that
$$f + I^\sigma_\alpha f^q \le u \le  f +
p^q \, I^\sigma_\alpha f^q.$$
\enddemo

\demo{Proof} The proofs of statements (1)$-$(4) of
the Theorem  follow easily from Theorems 4.8, 5.6, and 5.8.
One only  need to notice that the balls $B_r (x)$ associated
 with the Riesz metric  $\rho(x,y) = C(n, \alpha)^{-1} \,
|x-y|^{n-\alpha}$ correspond to Euclidean balls with the
same center and radius $C(n, \alpha)^{1/(n - \alpha)} \,
r^{1/(n - \alpha)}$.

Under
the assumptions of statement (5)  of
the theorem the weighted norm
inequality (6.3)
is equivalent to (6.7) by Theorem 2 of \cite{45}. To prove (5),
 note
that obviously (6.4)$\Rightarrow$(6.7), and apply (2).

To prove
statement (6), observe that since $\sigma \in RD_\beta$
then clearly (6.6) holds. Hence by (4)
(6.2)$\Leftrightarrow$(6.5).  It is also well known that
in this case (6.7) is necessary and sufficient in order
that
the weighted norm inequality (6.3) hold
(see \cite{41}, \cite{45}). Thus (6.5)$\Leftrightarrow$(6.7) which
proves (6).
\qed\enddemo

We now give an application to weighted Carleson measure
inequalities for Poisson integrals.
Let $\omega$ and $\sigma$ be locally finite Borel measures
on $\bold R^{n+1}_+=\bold R^{n}\times\bold R^n_+$
 and $\bold R^{n}$
respectively. Let $1<p< \infty$.
 We consider the inequality
$$||P^\sigma f||_{L^p (\omega)} \le C \,
||f||_{L^p (\sigma)} \tag 6.9$$
for  the Poisson integral
$$P^\sigma f (x, t) = P [f \, d \sigma] (x,t) =
\int_{\bold R^{n}} P_t (x - y) \, f(y) \, d \sigma (y).$$
Here $(x,t) \in \bold R^{n+1}_+$ and
$P_t (y) = C_n \, t/(||y||^2 + t^2)^{(n+1)/2}$ is
 the Poisson kernel on the upper half-space.
Similarly, for $\omega \in M_+({\bold R^{n+1}_+})$ we set
$$P^\omega g (x, t)=
 \int_{{\bold R^{n+1}_+}}
P_{t+ \tau} (x-y) \, g(y, \tau) \, d \omega (y, \tau).$$

If $\sigma$ is Lebesgue measure then the
Carleson measure
 theorem \cite{10} says that (6.9) holds if and only if
$\omega$ is a Carleson measure, i.e.
$|{\widehat B}|_\omega \le C \, |B|$; here
${\widehat B}$ is the cylinder with height $|B|^{1/n}$
whose base  is a ball $B$ in
$\bold R^{n}$. It is easy to see that $\omega$ is a
Carleson measure
if and only if
$$ P^\omega \bold 1 (x, t) = \int_{\bold R^{n+1}_+}
P_{t+ \tau} (x-y) \,
d \omega (y, \tau) \le C < \infty   \tag 6.10 $$
where $C$ is independent of $(x, t) \in \bold R^{n+1}_+$.
(See \cite{19}, \cite{40}.) A simple proof of the nontrivial
implication
(6.10)$\Rightarrow$(6.9) (in case $\sigma$ is Lebesgue
measure)
is immediate from the following
 theorem.

\proclaim{Theorem  6.3}
Let $1<p< \infty$
and let $\omega \in
M_+ ({\bold R^{n+1}_+})$,  $\sigma \in
M_+ (\bold R^{n})$. Then (6.9) holds if
$P^\sigma \bold 1 \not\equiv + \infty$ and
$$P^\omega (P^\sigma \bold 1)^p (x,t) \le C \, P^\sigma
\bold 1 (x,t). \tag 6.11$$
\endproclaim

The case $p=2$ of Theorem 6.3 was established
in \cite{50}. Note that the pointwise condition (6.11) is used
in a dual form so that if
$\sigma$ is Lebesgue  measure we have $P^\sigma \bold 1  = \bold 1$;
 then (6.11) coincides with (6.10).

\demo{Proof} Consider $d \omega_1 (x,t) = t \, d \omega (x,t)$
 and
$d \sigma_1= \chi_{\bold R^n} \, d \sigma$ as measures
 on $X =\overline {\bold R^{n+1}_+}$. Let $K$ be
a quasi-metric kernel on $X \times X$ defined by
$$K(\bar x, \bar y) = [||x-y||^2 + (t+ \tau)^2]^{-(n+1)/2}$$
where $\bar x = (x,t)$ and  $\bar y=(y, \tau)$.

It is easy to see that (6.11) is equivalent to
$$K^{\omega_1} (K^{\sigma_1} \bold 1)^p  \le C \, K^{\sigma_1} \bold 1.
\tag 6.12$$
Then by Theorem 4.8 applied to $L^q (\omega_1)$  in
place of $L^p (\sigma)$ it follows
$$||K^{\omega_1} g||_{L^q (\sigma_1)} \le C \,
||g||_{L^q (\omega_1)}$$
for all $g \in L^q (\omega_1)$,
which is obviously equivalent to the  inequality dual to (6.9).
 Note that here we  have only used an easy part of  Theorem 4.8
 which is essentially contained in the elementary estimates
of Propositions 4.4 and  4.5.
\qed
\enddemo

\vskip10pt\heading{7. Existence of positive solutions for
 superlinear
Dirichlet problems}\endheading\vskip10pt

In this section, we obtain some estimates for the Green
kernels and
\Naim  kernels  related to the so-called
$3G$-inequalities. This makes it possible, as an application of
 the results of
Sections 2-5, to characterize the problem of the existence
of positive solutions
for the superlinear Dirichlet problem
 $$\left \{ \aligned
- & \Delta u = v (x) \,  u^q +  w (x), \quad  u \ge 0 \quad
  \text {on} \quad
 \Omega,\\
 & u =  \phi \quad \text {on} \quad \partial \Omega,
\endaligned \right. \tag 7.1$$
  on a regular domain $\Omega \subset\bold R^n$ in the
 ``superlinear
case'' $q > 1$; here we assume that $v, w \in L^1_{loc}
  (\Omega)$ and
$\phi \in L^1_{loc}  (\partial \Omega)$
are arbitrary nonnegative functions.
We denote by  $G=G_{\Delta, \Omega}$ the Green function of
 the Laplacian $\Delta$ on $\Omega$, and by $G u$ the
 Green potential
$$G u (x) = \int_\Omega G (x, y) \, u (y) \, d y.$$
By $P \phi$ we denote the Poisson integral (harmonic
extension) of $\phi$.
The solvability of (7.1) is understood in the sense  (see \cite{6},
\cite{29})
that $u \in L^q_{loc} (\Omega)$ satisfies the corresponding
nonlinear integral equation
$$u = G \, (v \, u^q) +  G w + P \phi \quad
\text{a.e. on  $\Omega$}. \tag 7.2$$

More generally,
we consider the Dirichlet problem
$$\left \{ \aligned
- &L u = \sigma \, u^q + \omega, \quad  u \ge 0 \quad
  \text {on} \quad
 \Omega, \\  & u = \phi \quad \text {on} \quad
\partial \Omega, \endaligned \right.\tag 7.1{$'$}$$
where  $\sigma$, $\omega$ are locally finite measures on
$\Omega$ (possibly
singular with respect to Lebesgue measure), and
$L$ is a uniformly  elliptic differential operator of
second order,
$$L u = \sum_{i,j =1}^n a_{ij} (x) \,
 D_{ij} u + \sum_{i=1}^n b_i (x)
 \, D_i u
 + c(x) \, u, \tag 7.3$$
with the assumptions on the coefficients and $\Omega$
 specified below.  (See
Proposition 7.2 and Lemma 7.1.) Let $G=G_{L, \Omega}$ be
the Green function of $L$.
We say that $u \in L^q_{loc} (d \sigma, \Omega)$
 is a solution to $(7.1')$
if $u$ satisfies the integral equation
$$u = G^\sigma \,  u^q  + G^\omega \bold 1 + P \phi \quad
d \sigma\text{-a.e. on  $\Omega$}. \tag 7.2{$'$}$$
Here $G^\nu u (x)$ is  the Green potential
 with respect to a locally finite measure $u \, d \nu$
defined by
$$G^\nu u (x) = \int_\Omega G (x, y) \, u (y) \, d \nu (y),$$
and $P \phi = P_{L, \Omega} \phi$ is the corresponding Poisson
integral ($L$-harmonic extension of $\phi$).

 Our assumptions on $L$ and $\Omega$
will be stated in terms of certain inequalities for
the Green function $G=G_{L, \Omega}$.
As was mentioned in the Introduction, $G$ does
not satisfy
our basic quasi-metric assumptions even in the case of the
 Laplacian
on a ball or half-space. In the so-called $3 G$-inequality
(see \cite{12},
\cite{13}),
$$\frac{G (x, y) \, G(y,z)} {G(x,z)} \le \kappa
 \, (|x-y|^{2-n} +
|y-z|^{2-n}),  \tag 7.4$$
unfortunately, one cannot replace the right-hand side by a
smaller term $G(x,y) + G(y,z)$.

However, we can reduce the solvability problem
for  $(7.2')$ to the problems studied above by using the
 so-called \Naim kernel.
For $x \in \Omega$, we denote  by $\delta (x)$ the
distance from
$x$ to the boundary $\partial \Omega$. If $\partial \Omega$
is smooth
enough,
we set
$$N(x, y) =  \frac {G(x,y)} {\delta (x) \,
 \delta (y)}.\tag 7.5$$
Then we will show that $N$ satisfies
the inequality
$$\frac{N (x, y) \, N(y,z)}  {N(x,z)} \le
\kappa \, [N(x,y) +
N(y,z)],  \tag 7.6$$
which is obviously equivalent to our quasi-metric assumption.
It is worthwhile to note that (7.6) is always stronger than the
original $3 G$-inequality  (7.4), and takes into account the behavior
of the Green kernel at the boundary in a proper way. (See \cite{49}
for
another  refinement of (7.4) in a different direction,
which is, however, again  not sharp at the boundary.)

We observe that  this
approach is applicable at least to operators $L$ with
H\"older-continuous
coefficients on bounded $C^{1,1}$ domains $\Omega$.
For more general domains, as was proposed in \cite{38},
one may replace the distance
 to the boundary $\delta (x)$ in the definition (7.5) by
$s(x) = G(x, x_0)$,
where $x_0$ is a fixed pole in $\Omega$. We conjecture that in this
setting (7.6) holds, with obvious modifications, for
bounded Lipschitz domains $\Omega$ and operators  $L$ with
bounded measurable  coefficients.

The proof of (7.6) under the assumptions stated above
is based on the well known two-sided estimates
$$G_{L, \Omega}(x, y)
 \asymp
|x-y|^{2-n} \,
\min \, \left [1,  \frac {\delta (x) \, \delta (y)} {|x - y|^2}
 \right]\tag 7.7$$
for the Green function $G_{L, \Omega}$.
The upper estimate in  (7.7) is due to K.-O. Widman \cite{55}
for $C^{1, \alpha}$ (or more general Dini type) domains, and
the lower
one was established by
Z. Zhao \cite{56} for $C^{1, 1}$ domains and $L= \Delta$.
For second order uniformly elliptic operators $L$ with
H\"older-continuous
coefficients, it was proved by Hueber and Sieveking \cite{25}
that $G_{L, \Omega} \asymp G_{\Delta, \Omega}$.
(See also \cite{5}, \cite{49}, and the references given there.)

\proclaim{Lemma 7.1} Let $\Omega$ be a bounded $C^{1,1}$ domain in
$\bold R^n$, $ n \ge 3$. Let $L$ be a uniformly elliptic second
order differential operator (7.3) with bounded H\"older-continuous
coefficients, and $c \le 0$. Then $\rho (x, y) = N(x, y)^{-1}$
defines a quasi-metric on $\Omega$, and thus (7.6) holds.

\endproclaim
 \demo{Remark} In Lemma 7.1, one can replace the restrictions
on $c$  and $b_i$ by some milder assumptions which guarantee
that $G_{L, \Omega} \asymp G_{\Delta, \Omega}$
(see \cite{5}, \cite{12}, \cite{13}).
Analogues of Lemma 7.1 also hold in the cases
$n=1, 2$ which require
usual modifications. \enddemo

\demo{Proof}  We first
 show that  (7.7) is equivalent to the estimate
$$G_{L, \Omega}(x, y) \asymp \frac {\delta (x) \, \delta (y)}
{|x-y|^{n-2} \, [|x-y|^2 + \delta (x)^2 + \delta (y)^2]}.
\tag {7.7$'$}$$
We have
$$\min \, \left [1,  \frac {\delta (x) \, \delta (y)}
{|x - y|^2} \right] = \frac {1} {\max \,  \left [1,  \frac {|x - y|^2}
{\delta (x) \, \delta (y)}\right]}$$
$$\asymp   \frac {\delta (x) \, \delta (y)} {|x-y|^2 +
  \delta (x) \, \delta (y)}.$$
It remains to notice that
 from the obvious inequality $| \delta (x) - \delta (y)| \le
 |x-y|$ it follows
$$|x-y|^2 +
 \,  \delta (x) \, \delta (y) \asymp |x-y|^2 + \delta (x)^2 +
\delta (y)^2.$$
 Hence (7.7$'$) holds. Now for the \Naim kernel
defined by (7.5) we have
$$N(x, y)   \asymp \frac {1}
{|x-y|^{n-2} \, [|x-y|^2 + \delta (x)^2  + \delta (y)^2 ]}.$$
We set
$$d (x, y) =  |x-y|^{n-2} \,
[|x-y|^2 + \delta (x)^2 + \delta (y)^2],$$
where $d(x, y) \asymp N(x, y)^{-1}$. Then  to prove (7.6)
 it suffices to
show that $d (x, y)$
satisfies the quasi-metric inequality
$$d (x, y) \le C \, [d (x, z) + d (y, z)]. \tag 7.8$$
We notice
$$d (x, y) = |x-y|^{n} + \delta (x)^2  \, |x-y|^{n-2}
 + \delta (y)^2 \, |x-y|^{n-2}.$$
Estimating the first term on the right  from above,
we obviously have
$$ |x-y|^{n} \le 2^{n-1} \,  [|x-z|^{n} + |y-z|^{n}]
\le  2^{n-1} \, [d (x, z) + d (y, z)].$$
It remains
to estimate $\delta (x)^2  \, |x-y|^{n-2}$, since a similar
bound for $\delta (y)^2  \, |x-y|^{n-2}$ follows by
interchanging the roles of $x$ and $y$.

To prove the inequality
$$\delta (x)^2  \, |x-y|^{n-2} \le C \,
[d (x, z) + d (y, z)],$$
we consider two cases, (i) $\text { } \delta (x)
\le \delta (y)$
 and (ii) $\text { } \delta (y) \le \delta (x)$.

In case (i), clearly,
$$ \delta (x)^2  \, |x-y|^{n-2} \le 2^{n-3} \, \delta (x)^2
 \,  [|x-z|^{n-2} +  |y-z|^{n-2}] $$
$$\le  2^{n-3} \, \delta (x)^2 \,  |x-z|^{n-2}
+   2^{n-3}
\, \delta (y)^2 \,  |y-z|^{n-2} \le  2^{n-3} \, [d (x, z) +
d (y, z)],$$
which implies (7.8).

In case (ii), use the inequality $\delta (x) \le \delta (y) +
|x-y|$. Then we get
$$\delta (x)^2  \, |x-y|^{n-2} \le 2 \, |x-y|^{n} +
 2 \, \delta (y)^2 \, |x-y|^{n-2}.$$
As above, for the first term  on the right-hand side of
the preceding
inequality we have
 $ |x-y|^{n} \le  2^{n-1} \,  [d (x, z) + d (y, z)].$
Estimating the second term, we have
$$\delta (y)^2 \,
|x-y|^{n-2} \le   2^{n-3} \, [ \delta (y)^2 \, |x-z|^{n-2}
+ \delta (y)^2 \,
 |y-z|^{n-2}]$$
$$ \le  2^{n-3} \, [ \delta (x)^2 \, |x-z|^{n-2} + \delta (y)^2 \,
 |y-z|^{n-2}] \le  2^{n-3}
 \, [d (x, z) + d (y, z)].$$
This proves (7.8) in case (ii). Thus, (7.8) holds, which
implies (7.6).
The proof of Lemma 7.1 is complete.\qed\enddemo

Next, we show that the main results of
Sections 3-5 hold true in  a more general setting where
 the kernel
$G(x, y)$ is not necessarily symmetric and may fail to satisfy
the
quasi-metric inequality. We assume
that there exist positive measurable
 functions $s_1$ and $s_2$
such that
$$K (x,y) \asymp s_1(x) \, G(x, y) \, s_2 (y),
\quad x, y \in X, \tag 7.9$$
where $K (x,y)$ is symmetric and
satisfies the quasi-metric inequality. It follows
that this is the case
for the Green function $G(x, y)=G_{L, \Omega}$, under
the assumptions
of Lemma 7.1, if $s_1(x) = s_2(x) = \delta^{-1} (x)$;
 then $K$ coincides with
the \Naim kernel (7.5). Recall that $f \in \Cal Z_{q, G}$ if
the nonlinear integral equation
$$u = G^\sigma u^q + \epsilon \,
f \quad d \sigma\text{-a.e. on $X$}  \tag 7.10$$
has a solution for some $\epsilon > 0$.

\proclaim{Proposition 7.2} Suppose that $G$ is a kernel on
$X \times X$ such that (7.9) holds,
where $s_1$ and  $s_2$ are
 positive measurable function on $X$,
and $K(x,y)$ is a symmetric quasi-metric kernel.
Then $f \in \Cal Z_{q, G}$
 if and only if
$$G^\sigma (G^\sigma f^q)^q \le C \, G^\sigma f^q < \infty \quad
 d \sigma\text{-a.e.}  \tag 7.11$$
\endproclaim

\demo{Proof} We rewrite (7.10) in the equivalent form
$$\tilde u = K^{\tilde \sigma} \tilde u^q + \tilde f,
\tag 7.12$$
where
$$\tilde u = s_1 \, u, \quad \tilde f = s_1 \, f, \quad d \tilde
\sigma = s_2^{-1} \, s_1^{-q} \, d \sigma. \tag 7.13$$
By Theorem 4.1, this equation has a solution $\tilde u$ if and
only if
$$K^{\tilde \sigma} (K^{\tilde \sigma} \tilde f^q)^q \le C \,
K^{\tilde \sigma} \tilde f^q  < \infty \quad
 d \sigma\text{-a.e.}  \tag 7.14$$
Multiplying both sides of the preceding inequality by $s_1^{-1}$
and
using (7.9) and (7.13), we see that (7.14) is equivalent to (7.11).
 \qed\enddemo

\demo{Remark} One of the advantages of using
pointwise characterizations (7.11) is that they are invariant
under the transformation of the kernels given by (7.9).
It is easy to see that all other results of
Sections 3-5 have complete analogues for kernels $G$
such  that
$K(x, y) \asymp s_1(x) \, G(x, y) \, s_2 (y)$ satisfies the
quasi-metric
inequality.

 However, all geometric conditions should be modified
because of the new quasi-metric. Moreover,
 (7.9) leads to a change of weights in the
corresponding weighted norm inequalities.
Using (7.9) and
(7.13) in the same manner as in the proof of Proposition 7.2,
it is easily seen
that the weighted norm inequality
$$||G^\sigma h||_{L^p (d \omega)} \le C \, ||h||_{L^p (d \sigma)},
 \quad h \in L^p (d \sigma) \tag 7.15$$
for $G$ is equivalent to a similar inequality
$$||K^{\sigma_1} g||_{L^p (d \omega_1)}
\le C \, ||g||_{L^p (d \sigma_1)},
 \quad g \in L^p (d \sigma_1), \tag 7.15{$'$}$$
for $K$, where $d \sigma_1 = s_2^{1-q} \, d \sigma$ and
$d \omega_1=  s_1^{-p} \, d \omega$.
\enddemo

Using these results with $s_1 = s_2 = \delta^{-1}$,
together with the  testing
characterizations
 of weighted norm inequalities \cite{45},
we obtain the following characterization of
 two weight  inequalities for Green's potentials.

\proclaim{Theorem 7.3} Let $\Omega$ and $L$ satisfy the
assumptions
of Lemma 7.1. For $x \in \Omega$
and $a>0$, denote by $B=B_a(x)$
a ``ball'' associated with the \Naim  kernel:
$$B_a(x) = \{y \in \Omega: \, G(x, y) \ge a^{-1} \, \delta(x) \,
\delta (y) \}. \tag 7.16$$
 Then the two
weight inequality (7.15) holds for the Green potential $G^\sigma$
if and only if, for all balls $B=B_a(x)$,  both
$$\int_\Omega (G^\sigma \delta^{q-1} \chi_B)^p \, d \omega \le C
\, \int_B \delta^q \, d \sigma \tag 7.17$$
and
$$\int_\Omega (G^\omega \delta^{p-1} \chi_B)^q \, d \sigma \le C
\, \int_B \delta^p \, d \omega  \tag 7.18$$
hold.

\endproclaim

The same argument as in Proposition 7.2, together with
Theorem 4.9,
yields the following characterization of the solvability
 problem for
the nonlinear integral equation (7.10). (Note that
Theorem 7.3 is
not used in this proof.)

\proclaim{Theorem 7.4} Let $\Omega$ be a bounded $C^{1,1}$
domain in
$\bold R^n$, $ n \ge 3$, $1 < q < \infty$, and let $L$ be a uniformly
elliptic second
order differential operator as in Lemma 7.1.
 Let $f \in L^q_{loc} (d \sigma, \Omega)$
be a nonnegative  function on $\Omega$. Then the
following statements
are equivalent.
\newline (1) Equation (7.10) has a solution for
some $\epsilon > 0$.
\newline (2) There exists a constant $C >0$ such that
$$G^\sigma (G^\sigma f^q)^q \le C \, G^\sigma f^q < \infty \quad
d \sigma\text{-a.e.}\tag 7.19$$
\newline (3) Both the weighted inequality
$$\int_\Omega (G^\sigma h)^p \, f^q \, \delta^{1-p} \, d \sigma
\le C \, \int_\Omega h^p \, \delta^{1-p} \, d \sigma, \tag 7.20$$
and the  infinitesimal inequality
$$\left \{ \int_0^a \frac {\int_{B_t (x)} \delta^{1+q} \,
d \sigma} {t^2} \, d t \right \}^{1/q}
 \, \left \{ \int_a^\infty
 \frac {\int_{B_t (x)} f^q \, \delta \, d \sigma} {t^2} \,
d t \right \}^{1/p}
\le C, \tag 7.21$$
hold.
\newline (4) Both (7.21) and the testing inequality
$$\int_B (G^\sigma f^q \,  \chi_B)^q \, \delta \,
d \sigma \le C
\, \int_B f^q \, \delta \, d \sigma \tag 7.22$$
hold.
\endproclaim

Now we are in a position to
characterize the solvability of the nonlinear Dirichlet problem
(7.2$'$). The latter is
related to the integral equation (7.10) with
$f = f_1 + f_2$, where $f_1= G^\omega \bold 1$
and $f_2 = P \phi$.
We first  consider the homogeneous problem
where $\phi = 0$:
$$\left \{ \aligned
- &L u = \sigma \, u^q + \epsilon \, \omega, \quad
 u \ge 0 \quad
  \text {on} \quad
 \Omega, \\  & u = 0 \quad \text {on} \quad
\partial \Omega, \endaligned \right. \tag 7.23$$
which is equivalent to the integral equation
$$u = G^\sigma u^q + G^\omega \bold 1. \tag 7.24$$

 Since in this case the
inhomogeneous
term of (7.24) is a Green potential,
the characterizations of Theorem 7.4 may be simplified by using
Theorem 4.8 in place of Theorem 4.9. This yields

\proclaim{Theorem 7.5} Under the assumptions of Theorem 7.4 the
following statements
are equivalent.
\newline (1) The Dirichlet problem
(7.23)
has a solution for
some $\epsilon > 0$.
\newline (2) $G \omega < \infty$ and
there exists a constant $C >0$ such that
$$G^\sigma (G^\omega \bold 1)^q \le C \, G^\omega \bold 1. $$
\newline (3) Both the weighted inequality
$$\int_\Omega (G^\sigma h)^p \,  \delta^{1-p} \, d \omega
\le C \, \int_\Omega h^p \, \delta^{1-p} \, d \sigma, $$
and the  infinitesimal inequality
$$\left \{ \int_0^a \frac {\int_{B_t (x)} \delta^{1+q} \,
d \sigma} {t^2} \, d t \right \}^{1/q}
 \, \left \{ \int_a^\infty
 \frac {\int_{B_t (x)} \delta \, d \omega} {t^2} \,
d t \right \}^{1/p}
\le C$$
hold.
\newline (4) Both the infinitesimal
and the testing inequality
$$\int_B (G^\omega \chi_B)^q \, \delta \,
d \sigma \le C
\, \int_B  \delta \, d \omega$$
hold.
\endproclaim
\demo{Remark}  In the one-dimensional case, we consider
the problem
$$\left \{\aligned  - &u''(x)= \sigma \, u(x)^q  + \omega,
\quad 0 < x < 1, \\ & u(0) = u (1) = 0. \endaligned
\right.$$
 Here $\Omega = (0, 1)$ and $G (x, y) =
\min [x (1-y), y (1-x)]$.

The corresponding \Naim kernel $N(x, y) =
[x y (1-x) (1-y)]^{-1} \, G (x, y)$
satisfies the quasi-triangle inequality with $\kappa = 1$,
so that
$\rho(x, y) = 1/N(x,y) = [1- \min(x,y)] \, \max (x,y)$,
and thus Theorem 4.8 is applicable.  Note
that in the easy case where
 $\omega \in L^1 (0, 1)$ (or, more generally,
 if $\omega$ is a finite measure),
we have $G \omega (x) \asymp x (1-x)$. Then
(7.19) boils down to the requirement that the  integral
$\int_0^1 [x(1-x)]^q \, d \sigma$ is finite
and  small enough.
This is
clearly necessary for the existence of solutions in case
$\omega$ is a finite measure.

However, $\omega$ need not be finite on $(0, 1)$:
 the only natural restriction
 on $\omega$, which is  equivalent to $G \omega < \infty$,  is
 $\int_0^1 x (1-x) \, d \omega < \infty$.  For $\omega$ such
that
$\int_0^1 d \omega  = \infty$, the existence of positive
 solutions
depends
on the interplay between $\sigma$ and $\omega$ at the endpoints,
and is determined  by (7.19), or the equivalent inequalities (7.21)
and (7.22). Moreover, it is not difficult to see that
in this case
(7.21)$\Rightarrow$(7.22), and so the infinitesimal inequality
alone characterizes the solvability problem in this case.
\enddemo

Now we consider  the inhomogeneous Dirichlet problem
with  boundary data $\phi \not= 0$.
 We assume that $0 \le P \phi < \infty$.
(The estimates of the Poisson kernel for
elliptic  operators of second order can be found in \cite{5}, \cite{36},
 \cite{48}, \cite{56}.)

It follows from Theorem 7.4  that
(7.10) is solvable for some (small enough) $\epsilon> 0$
if and only
if both $G^\omega \bold 1 $ and $P \phi \in \Cal Z_{q, G}$.
Thus, applying  Theorems 7.4 and 7.5 established above,
we obtain the following statement.

\proclaim{Theorem 7.6} Under the assumptions of Theorem 7.4,
   the Dirichlet problem
$$\left \{ \aligned
- &L u = \sigma \, u^q + \epsilon \, \omega, \quad  u \ge 0 \quad
  \text {on} \quad
 \Omega, \\ & u = \epsilon \, \phi \quad \text {on} \quad
\partial \Omega, \endaligned \right.$$
 has a solution for some  $\epsilon > 0$ if and only if
both
$$G^\sigma (G^\omega \bold 1)^q
 \le C \, G^\omega \bold 1\quad
d \sigma\text{-a.e.}$$
 and
$$G^\sigma [G^\sigma (P \phi)^q]^q \le C \, G^\sigma (P \phi)^q
 \quad
d \sigma\text{-a.e.}$$
hold.
\endproclaim

Equivalent characterizations of solvability in terms of the
infinitesimal inequalities and testing inequalities (or
capacitary inequalities) follow as in
Theorems  7.4 and 7.5.

In conclusion we consider the special case of the nonlinear
Dirichlet problem
$$\left \{ \aligned
- &L u = u^q  + \epsilon \, \omega, \quad
 u \ge 0 \quad
  \text {on} \quad
 \Omega, \\  & u = 0 \quad \text {on} \quad
\partial \Omega, \endaligned \right. \tag 7.25$$
previously characterized by D. Adams and M. Pierre \cite{3}
in the case $\supp \omega \Subset \Omega$.
 We observe that our methods are totally different, and give
additional pointwise estimates of solutions with sharp constants
up to the boundary.

We introduce a weighted capacity associated with (7.25). To any
$E \subset \Omega$ we associate
$$\text{Cap}_{p, G}(E)  =  \inf \,
\left \{\int_\Omega  g^p \, \delta(x)^{1-p} \,
d x :    \  G g (x) \ge \delta (x) \chi_E (x),
\  g\ge  0 \right \}, \tag 7.26$$
where
$$G g (x) = \int_\Omega G(x, y) \, g(y) \, d y$$
is the Green potential of $g$.

Let
$$\text{Cap}_{p, I_2}(E)  =  \inf \, \left
\{\int_\Omega  g^p  \,
d x :    \  I_2 g (x) \ge  \chi_E (x),
\  g\ge  0  \right \} \tag 7.27$$
be the nonlinear Newtonian
capacity associated with the Sobolev space $W^{2,p} (\bold R^n)$
used in
\cite{3}.
It is easily seen that for any compact set $E \Subset \Omega$
one has $\text{Cap}_{p, G}(E) \asymp \text{Cap}_{p, I_2}(E)$
with  constants of equivalence which
depend on $\text{dist} \,
(E, \partial \Omega)$.

\proclaim{Theorem 7.7} Let $\Omega$  and $L$ be as in Theorem 7.4.
Let $1 < q < \infty$ and let $\omega \in M_+ (\Omega)$.
 Then the following statements are equivalent.
\newline (1) The Dirichlet problem
(7.25)
has a solution for
some $\epsilon > 0$.
\newline (2) $G \omega < \infty$ and
there exists a constant $C >0$ such that
$$G (G^\omega \bold 1)^q \le C \, G^\omega \bold 1. \tag 7.28$$
\newline (3) The weighted inequality
$$\int_\Omega (G h)^p \,  \delta^{1-p} \, d \omega
\le C \, \int_\Omega h^p \, \delta^{1-p} \, d \sigma, $$
holds.
\newline (4) The testing inequality
$$\int_{B\cap \Omega} (G^\omega \chi_B)^q \, \delta \,
d x \le C
\, \int_B  \delta \, d \omega$$
holds for all Euclidean balls $B$.
\newline (5) There exists a constant $C$ such that
$$ \int_E \delta(x) \, d \omega (x)
\le C \, \text{{\rm Cap}}_{p, G}(E) \tag 7.29$$
for all compact sets $E \subset \Omega$.
\endproclaim

\demo{Remarks} (a) The equivalence of (1) and (5) for compactly
 supported $\omega$ was proven in \cite{3}. As in our preceding
results, (7.28) with $C = p^{1-q} q^{-1}$ implies that equation
(7.25) with $\epsilon =1$ has a solution $u$ such that
$G \omega \le u \le p \, G \omega$.\newline
(b) Theorem 7.7 together with the estimate of solution
given above yields the following criterion for the  existence
of solutions in $L^r, \, 0 < r \le \infty$ (or any other ideal
function space): (7.25) has a solution for some  $\epsilon > 0$
which belongs  to $L^r$ if and only if
(7.28) holds and $G^\omega \bold 1 \in L^r$.
\enddemo

 \demo{Proof} The theorem follows from Theorem 7.6 if one
can show that the corresponding
infinitesimal inequality
is a consequence of the testing inequality. To prove this,
we apply Theorem 5.9 to the integral operator with the \Naim
kernel and then pass to the Green potential as above.

Recall that the \Naim
kernel $N(x,y)$ is quasi-metric and
$$ \rho(x,y) = N(x,y)^{-1}
\asymp |x-y|^{n-2} \, \max \{ |x-y|, \delta (x),
\delta (y)\}^2$$
(see the proof of Lemma 7.1).
It only remains to show that the assumptions
of Theorem 5.9 hold (see also Remark after Theorem 5.9), namely
that
$$|B_{2a}(x)|_{\sigma}\le   C \, a \, M(x,a)\tag  7.30$$
 and
$$ M(y,a)  \le C \, M(x,a)\tag  7.31$$
 whenever $x,y \in  X$
and $\rho(x,y)\le a.$ Here $B_a (x) = \{y \in \Omega: \,
\rho(x,y) < a$,
$d \sigma = \delta (x)^{1+q} \, dx$
and
$$ M(x,a) = \int_0^a \frac {|B_{t}(x)|_{\sigma}} {t^2} dt.$$

Denote Euclidean balls by $\Cal B = \Cal B_r(x)$.  Note that
if we set
$$\alpha = a^{1/n}, \quad \beta = a^{1/(n-2)} \,
\delta (x)^{-2/(n-2)},$$
 then  $B_a (x) \subset
\Cal B_{\min(\alpha, \beta)} (x)$. On the other hand
 if $|x-y| \le r$, then
$\rho (x, y) \le 4 \,
\max \{r^n, \, r^{n-2} \, \delta (x)^{2} \}$.
Hence
${\Cal B}_{4^{1/n} \min(\alpha, \beta)} (x) \subset B_a (x)$.

It is easy to see that if $a < C \, \delta(x)^n$ this leads
to an estimate
$$ |B_a (x)|_\sigma \sim \delta(x)^{1 + q} \, \beta^n =
a^{n/(n-2)} \, \delta(x)^{1 + q - 2 n/(n-2)}.$$
Let $\delta_0$ be such that if $ \delta(x) + r < \delta_0$ then
there exists $y \in  \Cal B_r(x)$ with
$\delta(x) = \delta(y) + r$.
If $a < 2^{-n} \, \delta_0^n$ and $2 \, \delta < a^{1/n}$
then pick
$y \in \Cal B_{\alpha/2} (x)$ with
 $\delta (y) \ge \alpha/2$. Then
$ \Cal B_{\alpha/2} (y) \subset B_a(x) \subset
\Cal B_{3 \alpha/2} (y)$. From this we get an estimate:
$$|B_a (x)|_\sigma \sim \alpha^{1+q} \, \alpha^n =
a^{1 + (1+q)/n}.$$
Combining these estimates we have
$$|B_a (x)|_\sigma \sim \max \{a^{n/(n-2)} \,
\delta(x)^{1 + q - 2n/(n-2)}, a^{1+ (1+q)/n} \}$$
for small enough $a < c \, \delta_0^n$.
Since $ |B_a (x)|_\sigma $ is bounded for $a \ge \delta_0$
now it is easy to verify that
 (7.30) and (7.31) hold.\qed\enddemo

\noindent {\bf Addendum.} After this paper was accepted for publication
we learned that H. Brezis and X. Cabre have been able to
modify their approach
in \cite{9} (see the Introduction) to find another proof
of the  necessity of our condition $G (\sigma (G \omega)^q) \le C \,
G \omega$ for the solvability of the Dirichlet problem for
$- \Delta u = \sigma \, u^q + \omega$.
Moreover, they showed that the
constant $C$ in the necessity part
can be chosen as  $C= p-1$, which complements our sufficiency result
with the sharp constant $C = q^{-1} p^{1-q}$. Note that the latter
constant is applicable  to very general superlinear operator equations,
while the former is established only for the Laplacian. For more general
superlinear  differential and integral
equations the constant $C = C(q, \kappa)$ in the necessity statements
 could be easily estimated from our proof; generally it depends on the
quasi-metric constant of the kernel.

\vfill
\eject

\bigskip
\enddocument

                  \end